\documentclass[12pt]{amsart}

\usepackage[utf8]{inputenc}
\usepackage[english]{babel}
\usepackage[margin=1.36in]{geometry} 
\usepackage{amsmath, amsfonts, amssymb, amscd}
\usepackage{enumitem}
\usepackage{hyperref}
\usepackage{mathtools}
\usepackage{xurl}

\usepackage{xcolor}

\makeatletter
\newsavebox{\@brx}
\newcommand{\llangle}[1][]{\savebox{\@brx}{\(\m@th{#1\langle}\)}%
  \mathopen{\copy\@brx\mkern2mu\kern-0.9\wd\@brx\usebox{\@brx}}}
\newcommand{\rrangle}[1][]{\savebox{\@brx}{\(\m@th{#1\rangle}\)}%
  \mathclose{\copy\@brx\mkern2mu\kern-0.9\wd\@brx\usebox{\@brx}}}
\makeatother

\let\le\leqslant
\let\ge\geqslant
\let\emptyset\varnothing
\let\epsilon\varepsilon
\let\phi\varphi

\DeclareMathOperator{\Diff}{Diff}
\DeclareMathOperator{\Dom}{Dom}
\DeclareMathOperator{\End}{End}

\DeclareMathOperator{\Hom}{Hom}

\DeclareMathOperator{\Fred}{Fred}
\DeclareMathOperator{\spn}{span}

\DeclareMathOperator{\Sym}{Sym}
\DeclareMathOperator{\codim}{codim}
\DeclareMathOperator{\coker}{coker}
\DeclareMathOperator{\gr}{gr}
\DeclareMathOperator{\im}{im}
\DeclareMathOperator{\ind}{ind}
\DeclareMathOperator{\lt}{lt}
\DeclareMathOperator{\ord}{ord}
\DeclareMathOperator{\pr}{pr}
\DeclareMathOperator{\rk}{rk}
\DeclareMathOperator{\tr}{tr}

\renewcommand{\d}{\mathrm{d}}
\renewcommand{\div}{\mathop{\mathrm{div}}\nolimits}

\newcommand{\const}{\mathrm{const}}
\newcommand{\fix}{\mathrm{fix}}

\newcommand{\iprod}{\mathbin{\lrcorner}}
\newcommand{\kere}{\mathop{\mathrm{ker}}\nolimits^e}
\newcommand{\reg}{\mathrm{reg}}
\newcommand{\rf}{\mathrm{ref}}
\newcommand{\sym}{\mathrm{sym}}
\newcommand{\loc}{\mathrm{loc}}

\numberwithin{equation}{section}

\theoremstyle{plain}
\newtheorem{proposition}{Proposition}[section]
\newtheorem{claim}[proposition]{Claim}

\newtheorem{corollary}[proposition]{Corollary}
\newtheorem{lemma}[proposition]{Lemma}
\newtheorem{theorem}{Theorem}

\theoremstyle{definition}
\newtheorem{definition}[proposition]{Definition}

\theoremstyle{remark}
\newtheorem{remark}[proposition]{Remark}
\newtheorem{example}[proposition]{Example}

\DeclareRobustCommand{\SkipTocEntry}[5]{}

\title[Generic equivariant simplicity]{Generic equivariant simplicity of the
spectrum of Laplace-type operators with finite symmetry}
\author{Egor Morozov\and Egor Shelukhin}

\address{D\'epartement de Math\'ematiques et de Statistique, Universit\'e de
Montr\'eal, C.P. 6128 Succ. Centre-Ville, Montreal (Québec), H3C 3J7, Canada} 

\email{egor.morozov@umontreal.ca} 
\email{egor.shelukhin@umontreal.ca}

\begin{document}

\begin{abstract}
We prove that for generic elements in three natural families of Laplace-type
differential operators on smooth manifolds, which are equivariant with respect
to a finite symmetry group, all the eigenspaces are irreducible representations.
In other words, the spectrum of these operators is generically equivariantly
simple. The three families of operators with symmetry are Schr\"odinger
operators parametrized by the potential, going back to the work of Wigner,
Laplace-Beltrami operators parametrized by invariant Riemannian metrics and
conformal Laplacians parametrized by invariant metrics in a given conformal
class. The latter two results extend prior work of Zelditch for certain free
group actions and resolve, for finite group actions, a well-known open question
raised by Guillemin and Yau. Our approach stems from methods developed by Wendl
and Doan-Walpuski to address equivariant transversality in symplectic topology.
\end{abstract}

\maketitle

\tableofcontents

\section{Introduction}

A celebrated result of Uhlenbeck~\cite{Uhl1976} states that for a generic, in a
suitable sense, Riemannian metric \(g\) on a closed connected manifold of
dimension at least \(2\), the associated Laplace-Beltrami operator \(\Delta_g\)
has a simple spectrum. That is, every eigenspace is of dimension one. Related
results, extensions and alternative approaches to this result appeared
in~\cite{Alb1975,Tan1979, BW1980,BU1983,Can2014,Col1988,Tey1999,GLS2014,CH2026}. 

In problems of geometric and physical origin manifolds often appear with
symmetries. It has therefore been of interest to establish whether or not a
similar statement holds in an equivariant sense when the manifold is endowed
with an action of a compact Lie group \(G\). Specifically, every eigenspace of
the Laplace-Beltrami operator of a \(G\)-invariant metric becomes a
\(G\)-representation, and one asks whether or not every such representation is
irreducible for a generic \(G\)-invariant metric. This question of \emph{generic
equivariant simplicity} of the spectrum or \emph{generic irreducibility} of
eigenspaces was asked by Guillemin~\cite{Sch2017} and by Yau~\cite[Problem
42]{Yau1993}. It is reinforced by results of Donnelly~\cite{Don1978,Don1976}
stating for instance that for finite \(G\) acting with trivial minimal isotropy
group, a real irreducible representation \(\rho\) asymptotically occupies
\(\deg(\rho)^2/(\alpha|G|)\) of the total dimension of the sum of the
eigenspaces up to any given eigenvalue (here \(\alpha=\dim_{\mathbb
R}\End_G(\rho)\in\{1,2,4\}\) depending on whether \(\rho\) is of real, complex,
or quaternionic type).


From the perspective of physics, this question corresponds to the classical
notion of \emph{accidental degeneracy} of energy levels of a quantum mechanical
system with symmetry. For instance, as observed by Wigner~\cite{Wig1959}, every
multi-particle system composed of \(N\) identical scalar particles admits a
\(G=S_N\) symmetry by permutations. This forces a necessary degeneracy of the
energy levels corresponding to irreducible \(G\)-representations, and generic
equivariant simplicity corresponds to the absence of further, accidental,
degeneracy. We refer to~\cite{FadYak2009} for further background on the physical
aspects of this question.

Generic irreducibility was studied in~\cite{Sch2017} for Lie groups, such as
\(SU(2)\) and \(SU(2)^m\times T^n\) acting on themselves, where it was shown to
hold. This result was extended in~\cite{MarGo2019}. Moreover, \cite{CJLS2022}
establishes generic equivariant simplicity for free actions of \(G=T^n\) on
manifolds (see also~\cite{JunZel2020} for \(n=1\)) of dimension \(>n\).
Finally,~\cite{PR2018} studies it for \(G\)-invariant metrics on homogeneous
spaces \(G/K\), establishing it in some cases, and proving that it \emph{fails}
whenever \(G/K\) is a compact irreducible symmetric space of rank at least two.
Thus, one does not expect generic irreducibility to hold for arbitrary actions
of compact Lie groups of positive dimension. However, one may still expect it to
hold in dimension zero, that is, for finite groups \(G\). This case was studied
by Zelditch~\cite{Zel1990} in the case of a free finite group action on a
manifold \(M\), namely a Riemannian cover, where equivariant simplicity was
established under the rather stringent ``HDLD" condition that every real
irreducible representation \(\rho\) of \(G\) has \(\deg(\rho)\le\dim(M)\). We
note that a similar low-degree property is essential for the proofs of the above
results regarding \(T^n\)-actions.

In this paper we study the generic equivariant simplicity of the spectrum of
differential operators on closed manifolds invariant under an arbitrary finite
group action. In particular, we impose no low-degree condition. We establish it
for three kinds of families of Laplace-type operators: the Laplace-Beltrami
operator, the conformal Laplacian, and Schr\"odinger operators with smooth
potentials. All manifolds are assumed to be connected.

\begin{theorem}\label{th:intro-yau}
Let \(G\) be a finite group acting on a closed manifold \(M\) of dimension at
least \(2\). For a residual set of smooth \(G\)-invariant Riemannian metrics
\(g\) on \(M\), every eigenspace of \(\Delta_g\) is an irreducible real
\(G\)-representation.
\end{theorem}
 
We note that there are many closed manifolds where every smooth action of a
compact Lie group factors through that of a finite group. Indeed, these are the
manifolds which do not admit smooth circle actions: for instance closed spin
\(4k\)-manifolds with non-vanishing \(\widehat A\)-genus~\cite{AH1970}, certain
complete intersections in complex projective spaces~\cite{DW2017}, and other
manifolds~\cite{FM1978, AB1981, Puppe1988}. In particular, we completely resolve
Guillemin's and Yau's question for such manifolds.

We generalize our first main result to the conformal Laplacian in arbitrary
dimensions in Theorem~\ref{th:mainres-confcov} below. Furthermore, we prove a
generic simplicity result for Schr\"odinger operators on manifolds varying the
\(G\)-invariant potential. This provides an equivariant generalization of a
classical theorem of Albert~\cite{Alb1975}.

\begin{theorem}\label{th:intro-sm}
Let \(G\) be a finite group acting on a closed manifold \(M\) and \(g\) be a
\(G\)-invariant Riemannian metric. For a residual set of smooth \(G\)-invariant
potentials \(f\in C^\infty(M)\) on \(M\), every eigenspace of the operator
\(L_f=\Delta_g+f\) is an irreducible real \(G\)-representation.
\end{theorem}

The example of multi-particle systems with \(S_N\)-symmetry on the sphere, and
having a smooth potential is of pertinence to systems of interest in physics
with soft interactions: for example the Maier-Saupe theory in liquid crystal
theory (see~\cite[Section 10.4]{Sel2016}) measuring alignment of orientations of
cylindrical mo\-le\-cu\-les, or to suitable versions of the quantum harmonic
oscillator. Note that the case of smooth potential does not directly apply to
systems with hard interactions, without regularization, as they tend to appear
more abundantly in theoretical physics. However, we expect that variations of
our approach could apply to such situations as well. We illustrate our
expectation by addressing an example case of a Coulomb gas on \(\mathbb R^3\) in
a confining potential as follows. 

\begin{theorem}\label{th:intro-co}
Let \(M=(\mathbb R^3)^N\) with \(G=S_N\) action by permutations of particles.
Let
\[
U=\sum_{i=1}^N v(x_i)+\sum_{1\le i<j\le N}C(|x_i-x_j|)
\]
for \(v\colon\mathbb R^3\to\mathbb R\) smooth, \(v(x)\to +\infty\) as \(|x|\to
\infty\), \(C\colon(0,\infty) \to \mathbb R\) smooth, \(C(y) \sim \frac 1y\)
near \(y=0\) and \(C \ge 0\) for large \(y\). For generic smooth compactly
supported \(G\)-invariant \(f\colon M\to\mathbb R\),
\[
L_f=\Delta+U+f
\]
has discrete \(G\)-simple spectrum, where \(\Delta\) is the standard Laplacian
on \(M\).
\end{theorem}

Here genericity should be understood as follows: fix an arbitrary \(S>0\). Then
for a residual set of \(G\)-invariant \(f\in C^\infty(M)\) supported in the
closed ball \(B_S\) of radius \(S\), the conclusion holds.

This theorem extends in a straightforward way to a number of cases: for
instance, to the case where \(M=X^N\) for a closed manifold \(X\) of dimension
at least~\(3\), to compact \(X\) with Lipschitz boundary where Dirichlet
boundary conditions are imposed, and to certain non-compact manifolds. Certain
conditions on \(v,C\), smoothness for example, can be weakened. However, the
requirement \(C(y)\sim\frac 1{y}\) near \(y=0\) appears to be quite important.
We do not pursue a maximal generalization of this result here as it is merely an
illustration of our methods in a more physically pertinent situation. 

\bigskip

\noindent\emph{Future directions.} In future work, we will return to such
generalizations as well as to generalizations of our other results, Theorems
\ref{th:intro-yau} and \ref{th:intro-sm} for instance, to the real-analytic
case: where all the metrics and the potentials are real-analytic. Furthermore,
we will study the case of elliptic operators on vector bundles, with
applications to the Dirac operator in particular. Finally, we will study
super-rigidity questions for minimal submanifolds and other distinguished
submanifolds in manifolds with special holonomy: for instance, associative
submanifolds of \(G_2\)-manifolds.


\bigskip

\noindent\emph{Methods of proof.} Let us now briefly discuss the methods of our
proofs. Previous works on the subject relied on more classical methods, such as
the Sard-Smale transversality theorem, which are in general hard to apply
because perturbations of \(G\)-invariant objects are not local. Indeed, this is
ultimately the source of the HDLD condition in Zelditch's work. Instead, we use,
and further develop in our particular situations of interest, the techniques
introduced by Wendl~\cite{Wen2023} to solve the super-rigidity conjecture for
\(J\)-holomorphic curves in symplectic topology. One key point of the argument
is a reduction of a natural condition ensuring generic simplicity of the
spectrum, the Petri condition, to a local question depending only on the
situation on jets at a point, which readily extends to the equivariant case.
That the resulting condition on jets itself holds generically boils down to a
natural rank estimate of certain differential operators, determined by symbols
of the initial operators, acting on polynomials up to a given degree. This
so-called polynomial Wendl condition, was previously established by Wendl (see
also~\cite{DW2023}) for Cauchy-Riemann operators on Riemann surfaces, involving
polynomials of two variables. The second key point in our argument is
establishing the polynomial Wendl condition for the symbol of the Laplacian
acting on polynomials in arbitrarily many variables. This requires a new
technique, which is a certain version of the Fourier transform for polynomials,
related to the inner product on polynomials encountered in the study of
spherical harmonics. Using the same technique, we also provide first examples
where the polynomial Wendl condition fails, for instance for the symbol of the
bi-Laplacian. 

\subsection{Rigorous statements of main results}\label{sec:mainres}

In this section we restate our main results in more precise terms and fix some
notation. Throughout this paper \((M,g_\fix)\) is an \(n\)-dimensional closed
smooth Riemannian manifold and \(G\) is a finite group acting on \(M\) by
\(g_\fix\)-isometries. We do not suppose that the action is effective, that is,
the maximal subgroup \(G_0\subset G\) acting trivially on \(M\) may be
nontrivial. A point \(x\) is called \textbf{regular} if the stabilizer of \(x\)
is \(G_0\). It follows easily from the principal orbit type theorem that the set
of regular points is open and dense. If \(G\) acts on a set \(V\) (that is,
\(V\) is a \textbf{\(G\)-space}), then the subset of \(G\)-invariants is denoted
by \(V^G\). If \(V\) and \(W\) are \(G\)-representations, we denote
\(\Hom_G(V,W):=\Hom(V,W)^G\). A similar
notation is used for all \(G\)-equivariant mapping spaces.

A scalar differential operator \(\mathbf L\colon C^\infty(M)\to C^\infty(M)\) is
called \textbf{\(G\)-equivariant} if
\[
\mathbf L(gu)=g(\mathbf Lu)
\]
for all \(u\in C^\infty(M)\) and \(g\in G\). The action of \(G\) on
\(C^\infty(M)\) turns each eigenspace of a \(G\)-equivariant operator \(\mathbf
L\) into a \(G\)-representation. If \(\mathbf L\) is elliptic, then all its
eigenspaces are finite-dimensional. It is a standard result of representation
theory that \(G\) has finitely many irreducible representations
\(W_1,\ldots,W_N\) and any other finite-dimensional \(G\)-representation \(V\)
uniquely decomposes in a finite direct sum of irreducible ones, that is,
\[
V\cong W_1^{\oplus m_1}\oplus\ldots\oplus W_N^{\oplus m_N}
\]
for some \(\overline m=(m_1,\ldots,m_N)\in\mathbb N_0^N\). We just write
\([V]=\overline m\) in this case. For each \(j=1,\ldots,N\) the number \(m_j\)
is called the \textbf{multiplicity} of \(W_j\) in \(V\).

A subset of a topological space is called \textbf{Baire} (or \emph{comeagre} or
\emph{residual}) if it is a countable intersection of open and dense subsets.

\begin{definition}
Let \(\mathcal P\) be a topological space and \((\mathbf L_p)_{p\in\mathcal P}\)
be a family of linear elliptic formally self-adjoint \(G\)-equivariant
differential operators on \(M\) parametrized by \(\mathcal P\). We say that the
spectrum of \((\mathbf L_p)_{p\in\mathcal P}\) is \textbf{equivariantly
generically simple} if there exists a Baire subset \(\mathcal
P_\reg\subset\mathcal P\) such that for each \(p\in\mathcal P_\reg\) all the
eigenspaces of \(\mathbf L_p\) are irreducible \(G\)-representations. Similarly,
we say that the spectrum of \((\mathbf L_p)_{p\in\mathcal P}\) is equivariantly
generically simple \textbf{at nonzero eigenvalues} if for each \(p\in\mathcal
P_\reg\) all the eigenspaces of \(\mathbf L_p\) \emph{corresponding to nonzero
eigenvalues} are irreducible \(G\)-representations.
\end{definition}

Let \(\mathcal G\) be the space of all smooth Riemannian metrics on \(M\). Note
that \(\mathcal G\) both carries a canonical topology and a Fr\'echet manifold
structure (but not a Banach manifold structure). For each \(g\in\mathcal G\),
let \(\Delta_g\) denote the Laplace-Beltrami operator on \((M,g)\) given by
\(\Delta_g u=-\div_g \nabla_g u\). The first main result of the paper is the
generic equivariant simplicity of spectrum of the Laplace-Beltrami operator with
respect to variations of the metric.

\begin{theorem}\label{th:mainres-laplace}
Let \(M\) be a closed manifold of dimension at least \(2\) and \(G\) be a finite
group acting on \(M\) by diffeomorphisms. Let \(\mathcal P\) be an open subset
of \(\mathcal G^G\). Then the spectrum of the family of operators
\((\Delta_g)_{g\in\mathcal P}\) is equivariantly generically simple.
\end{theorem}

Note that Theorem~\ref{th:mainres-laplace}, as well as Uhlenbeck's celebrated
result, are easily seen to fail for the circle \(M=S^1\). Therefore the
dimensional restriction is necessary. Our second main result is about families
of Schr\"odinger operators of order two with varying potential. Recall that
\(C_G^\infty(M)\) denotes \(G\)-invariant smooth functions on \(M\).

\begin{theorem}\label{th:mainres-schrodlap}
Let \(M\) be a closed Riemannian manifold and \(G\) be a finite group acting on
\(M\) by isometries. Let \(\mathbf L\) be a linear, elliptic, formally
self-adjoint \(G\)-equivariant differential operator on \(M\) of order two and
\(\mathcal P\) be an open subset in {\(C_G^\infty(M)\).} Then the spectrum of
the family of operators \((\mathbf L+f)_{f\in\mathcal P}\) is equivariantly
generically simple.
\end{theorem}

Let \(\mathcal C\subset\mathcal G\) stand for an arbitrary conformal class on a
closed manifold \(M\) of dimension \(n\ge 2\). For a Riemannian metric
\(g\in\mathcal G\) the corresponding \textbf{conformal Laplacian} is given by
\[
\mathbf L_g:=\Delta_g+\frac{n-2}{4(n-1)}R_g,
\]
where \(R_g\) is the scalar curvature of \(g\). Note that \(\mathbf L_g\) is a
\(G\)-equivariant differential operator as soon as \(g\in\mathcal G^G\). Our
last main result is the generic equivariant simplicity of spectrum for the
conformal Laplacian with respect to conformal variations of the metric.

\begin{theorem}\label{th:mainres-confcov}
Let \(M\) be a closed manifold of dimension at least \(2\) and \(G\) be a finite
group acting on \(M\) by diffeomorphisms. Let \(\mathcal P\) be an open subset
of either \(\mathcal G^G\) or \(\mathcal C^G\) where \(\mathcal C\) is the
conformal class of a \(G\)-invariant metric. Then the spectrum of the family of
operators \((\mathbf L_g)_{g\in\mathcal P}\), where \(\mathbf L_g\) is the
conformal Laplacian on \((M,g)\), is equivariantly generically simple at nonzero
eigenvalues.
\end{theorem}

On surfaces, the conformal Laplacian coincides with the Laplace-Beltrami
operator so Theorems~\ref{th:mainres-laplace} and~\ref{th:mainres-confcov} are
essentially the same in this case. Moreover, our proof of
Theorem~\ref{th:mainres-laplace} combines Theorem~\ref{th:mainres-confcov} for
surfaces and a different argument for \(\dim M>2\).

%
%

\subsection{Outline of the proof}

Theorem~\ref{th:mainres-schrodlap} follows from the following general
technical result.

\begin{theorem}\label{th:outline-schrod} Let \(M\) be a closed Riemannian
manifold and \(G\) be a finite group acting on \(M\) by isometries. Let
\(\mathbf L\) be a linear, elliptic, formally self-adjoint \(G\)-equivariant
differential operator on \(M\) and \(\mathcal P\) be an open subset in
\(C_G^\infty(M)\). Suppose that there exists a regular
point \(x\in M\) such that for any \(f\in\mathcal P\) and \(\lambda\in\mathbb
R\) the operator \(\mathbf L+f-\lambda\) has the strong unique continuation
property at \(x\) and the principal symbol \(\sigma_x(\mathbf L)\) satisfies the
polynomial Wendl condition. Then the spectrum of the family of operators
\((\mathbf L+f)_{f\in\mathcal P}\) is equivariantly generically simple.
\end{theorem}

Similarly, Theorem~\ref{th:mainres-confcov} is a consequence of the following
result.

\begin{theorem}\label{th:outline-confcov}
Let \(M\) be a closed manifold and \(G\) be a finite group acting on \(M\) by
diffeomorphisms. Let \(\mathbf L\) be a \(G\)-equivariant conformally covariant
differential operator on \(M\) with smooth coefficients and \(\mathcal P\) be an
open subset in either \(\mathcal G^G\) or \(\mathcal C^G\). Suppose that there
exists a regular point \(x\in M\) such that for any \(g\in\mathcal P\) and
nonzero \(\lambda\in\mathbb R\) the operator \(\mathbf L_g-\lambda\) has the
strong unique continuation property at \(x\) and the principal symbol
\(\sigma_x(\mathbf L_g)\) satisfies the polynomial Wendl condition. Then the
spectrum of the family of operators \((\mathbf L_g)_{g\in\mathcal P}\) is
equivariantly generically simple at nonzero eigenvalues.
\end{theorem}

Theorem~\ref{th:mainres-laplace} is not a direct consequence of either one of
these results but its proof is based on similar ideas.

For the definition of a \(G\)-equivariant conformally covariant differential
operator with smooth coefficients see Section~\ref{sec:fam} and
Example~\ref{ex:fam-cc}. The strong unique continuation property is discussed in
Section~\ref{sec:ucon}. The \emph{polynomial Wendl condition} is a certain
linear algebraic yet nontrivial condition on the principal symbol of an elliptic
differential operator. See the discussion below and Section~\ref{sec:wendl} for
more details.

It may be possible to apply Theorems~\ref{th:outline-schrod}
and~\ref{th:outline-confcov} to other families of operators. However, the strong
unique continuation property and the polynomial Wendl condition are subtle
points for such generalizations. Indeed, both properties are not easy to verify
and sometimes fail. For example, in Section~\ref{sec:wendl} we show that the
polynomial Wendl condition fails for certain powers of Laplacian. For higher
order operators, the strong unique continuation property also becomes an issue:
for example, see~\cite[Section~2.3]{HPP2025} for a discussion concerning the
GJMS operators.

\smallskip

It was already Uhlenbeck who observed that generic simplicity of spectrum is
closely related to parametric transversality. However, a direct generalization
of Uhlenbeck's original approach to the equivariant situation encounters
significant difficulties, forcing one to impose additional assumptions such as
Zelditch's HDLD condition.
In this paper we avoid these problems by implementing a new technique with
origins in symplectic geometry. This technique was developed by C.\,Wendl in his
proof of the super-rigidity conjecture for \(J\)-holomorphic
curves~\cite{Wen2023}. Here we use a refined version inspired by the work of
A.\,Doan and T.\,Walpuski~\cite{DW2023}. We give a sketch of the whole argument
here.

Suppose that we want to show that the spectrum of a family of (scalar, linear,
elliptic, formally-self adjoint) \(G\)-equivariant operators \((\mathbf
L_p)_{p\in\mathcal P}\) on \(M\) is equivariantly generically simple. Although
in all natural settings including Theorems~\ref{th:outline-schrod}
and~\ref{th:outline-confcov}, \(\mathcal P\) is a \emph{Fr\'echet} manifold, for
the sake of simplicity and intuition let us now assume that \(\mathcal P\) is a
\emph{Banach} manifold. In fact, Floer's method of \(C_\epsilon\)-spaces (see
Section~\ref{sec:floer}) allows one to essentially reduce considerations to this
case. 

One can think of \(\mathbf L\) as a map \(\mathbf
L\colon\mathcal P\to\Fred_G^\sym(\mathcal D,\mathcal H)\), where \(\mathcal
D,\mathcal H\) are appropriate functional spaces and \(\Fred_G^\sym(\mathcal
D,\mathcal H)\) is the corresponding space of \(G\)-equivariant symmetric
Fredholm operators. We assume that \(\mathbf L\) interpreted this way is a
smooth map of Banach manifolds.

Now for each \(K\in\mathbb N\), let \(\mathcal P_\reg^K\subset\mathcal P\) be
such that for each \(p\in\mathcal P_\reg^K\) and eigenvalue \(\lambda\in
[-K,K]\) of \(\mathbf L_p\), the eigenspace \(\ker(\mathbf L_p-\lambda)\) is an
irreducible \(G\)-representation. We put \(\mathcal P_\reg:=\bigcap_{K\in\mathbb
N}\mathcal P_\reg^K\) and prove that it is a Baire subspace. It suffices to show
that \(\mathcal P_\reg^K\) is open and dense in \(\mathcal P\) for each
\(K\in\mathbb N\). That \(\mathcal P_\reg^K\) is open is easy and follows from
an appropriate version of the continuous dependence of eigenvalues on the
parameter. Showing that \(\mathcal P_\reg^K\) is dense is much more challenging
and is known as \emph{splitting the eigenvalues}. This is the point where the
parametric transversality comes into play. Extend \(\mathbf L\) to a map
\[
\mathbf T\colon\mathcal P\times\mathbb R\to\Fred_G^\sym(\mathcal D,\mathcal
H),\qquad
\mathbf T_{p,\lambda}:=\mathbf L_p-\lambda.
\]
It turns out that \(\Fred_G^\sym(\mathcal D,\mathcal H)\) is stratified by
submanifolds according to the decomposition of the kernel of the operator into
irreducible \(G\)-representations. Moreover, the codimension of each stratum can
be computed using Schur's lemma and it turns out that this codimension is equal
to one if and only if the kernel is an irreducible representation. Now let
\(p\notin\mathcal P_\reg^K\) and suppose that \(\lambda\in [-K,K]\) is such that
\(\ker\mathbf T_{p,\lambda}\) is not an irreducible representation. By
continuity of eigenvalues, for any point \(p_0\) near \(p\) there exists an
eigenvalue \(\lambda_0\) near \(\lambda\). If we are able to choose
\((p_0,\lambda_0)\) so that \(\mathbf T\) is regular at \((p_0,\lambda_0)\),
then we can split the eigenvalue \(\lambda\) by perturbing \(p\). Indeed, if
\(\ker\mathbf T_{p_0,\lambda_0}\) is an irreducible \(G\)-representation, then
we just take \(p_0\). Otherwise, since \(\mathbf T\) is regular at
\((p_0,\lambda_0)\), the preimage under the map \(\mathbf T\) of the stratum
containing \(\mathbf T_{p_0,\lambda_0}\) has codimension at least two in
\(\mathcal P\times\mathbb R\). The possibility to split the eigenvalue
\(\lambda\) now follows from the Sard-Smale theorem applied to the projection
\(\mathcal P\times\mathbb R\to\mathcal P\). Repeating the procedure for the
remaining eigenvalues inside \([-K,K]\), we can hope to split them all.

Thus, splitting the eigenvalues essentially reduces to showing that the map
\(\mathbf T\) is regular for a generic \((p,\lambda)\). It turns out that in
some situations the regularity condition at \((p,\lambda)\) reduces to an easily
formulated condition on the corresponding eigenspace \(\ker\mathbf
T_{p,\lambda}\). This is the so-called \emph{\(G\)-equivariant Petri condition}
considered in~\cite{DW2023}. Unfortunately, this condition is difficult to
verify and, in general, it is not satisfied everywhere (see
Example~\ref{ex:ucon-petrifail}). However, an insight by Wendl~\cite{Wen2019}
provides a method to show that Petri condition is satisfied generically. In
summary, the method consists of two steps. First, using the strong unique
continuation property, Petri condition is reduced to a similar condition on
\(\infty\)-jets, the \emph{\(\infty\)-jet Petri condition}. Although this
condition is still difficult to verify, it is local and allows one to deal with
issues of equivariance. Secondly, truncating \(\infty\)-jets to finite jets and
applying Sard-Smale again, one reduces the \(\infty\)-jet Petri condition to a
purely linear algebraic condition on the principal symbol of \(\mathbf L_p\).
This is exactly the polynomial Wendl condition of
Theorems~\ref{th:outline-schrod} and~\ref{th:outline-confcov}.

In case of Theorem~\ref{th:mainres-laplace}, the regularity condition reduces not
to \(G\)-equivariant Petri condition but to its variant, which we call
\emph{traceless \(\Delta\)-Petri condition}. Although this condition is more
tricky, its genericity can be proved by similar methods.

The argument goes almost as described above modulo one important technical
detail. As mentioned above, in all natural settings including
Theorems~\ref{th:mainres-laplace}--\ref{th:mainres-confcov}, the parameter space
\(\mathcal P\) is only a Fr\'echet manifold, not a Banach manifold. This means
that we cannot apply Sard-Smale theorem directly. The issue is fixed by another
pair of techniques from symplectic geometry called \emph{Floer
\(C_\epsilon\)-spaces} and \emph{Taubes' trick}. Roughly speaking, the idea is
to replace a Fr\'echet manifold \(\mathcal P\) by a Banach manifold \(\mathcal
P_\epsilon\), which embeds continuously into \(\mathcal P\) and approximates
\(\mathcal P\) better and better as \(\epsilon\to 0\), see
Section~\ref{sec:floer} for the details. Although it is a technical issue, we
consider it to be quite important since it goes through the whole exposition and
should be always kept in mind.

\subsection{Organization of the paper}

We start with numerous preliminaries, which we divide into two parts. The first
part, Section~\ref{sec:prealg}, contains algebraic (in fact, \emph{linear}
algebraic) preliminaries. Those include the language of filtered vector spaces
(Section~\ref{sec:filt}) and jet spaces (Section~\ref{sec:jets}), as well as
some basic facts on the scalar elliptic symbols. The second part,
Section~\ref{sec:preanal}, contains analytic preliminaries. Those are more
specific and include the stratification of the space of \(G\)-equivariant
Fredholm operators (Section~\ref{sec:walls}), Floer \(\epsilon\)-spaces
(Section~\ref{sec:floer}), \(C^\infty\)-subvarieties
(Section~\ref{sec:cinfsub}), and some technical details on the families of
\(G\)-equivariant differential operators such as used in the paper
(Section~\ref{sec:fam}) and an appropriate version of the continuity of
eigenvalues with respect to the parameter (Section~\ref{sec:cont}).

Theorem~\ref{th:outline-schrod} is proved in Section~\ref{sec:schrod}. Each
subsection is devoted to a particular step of the proof. In
Section~\ref{sec:petri} we show how splitting the eigenvalues is related to the
\(G\)-equivariant Petri condition. In Section~\ref{sec:ucon} we apply the strong
unique continuation property to reduce the \(G\)-equivariant Petri condition to
the \(\infty\)-jet Petri condition. In Section~\ref{sec:petrigen} we use the
polynomial Wendl condition and Taubes trick to show that \(\infty\)-jet Petri
condition is generic. In Section~\ref{sec:split} we gather everything together
and use Taubes trick once again to split the eigenvalues. In the end of the
section we deduce Theorem~\ref{th:mainres-schrodlap} from
Theorem~\ref{th:outline-schrod} modulo the verification of the polynomial Wendl
condition for the Laplace operator, which is postponed to
Section~\ref{sec:wendl}.

Theorems~\ref{th:outline-confcov} and~\ref{th:mainres-laplace} are discussed in
Sections~\ref{sec:confcov} and~\ref{sec:laplace} respectively. Since their
proofs are mostly similar to that of Theorem~\ref{th:outline-schrod}, we do not
give all the details and discuss only the points where the proof of
Theorem~\ref{th:outline-schrod} must be modified. Again, for
Theorem~\ref{th:mainres-laplace}, the verification of the corresponding Wendl
condition is postponed to Section~\ref{sec:wendl}.

All variants of the Wendl condition are discussed in Section~\ref{sec:wendl}.
This section is pure algebra. The jet Wendl condition and the polynomial Wendl
condition are covered in Sections~\ref{sec:jetwendl} and~\ref{sec:polywendl}. In
Sections~\ref{sec:coord}--\ref{sec:wendliff} we state and prove a criterion for
the polynomial Wendl condition and then use this criterion to verify the
polynomial Wendl condition for the symbol of the Laplace operator. In
Section~\ref{sec:tdwendl} we cover the version of the Wendl conditions that we
use in the proof of Theorem~\ref{th:mainres-laplace}.

Finally Appendix~\ref{sec:coulomb} describes the proof of Theorem
\ref{th:intro-co}, the new arguments being solely of functional-analytic
nature.

\addtocontents{toc}{\SkipTocEntry}
\subsection*{AI disclosure} 
Some standard but technical results in Appendix~\ref{sec:coulomb} were obtained
with the assistance of ChatGPT Astra. It has also helped us fix a few minor gaps
and typos in the initial version of the paper, which was written without AI
assistance. We have verified all arguments in the appendix and take
responsibility for their correctness.

\addtocontents{toc}{\SkipTocEntry}
\subsection*{Acknowledgments}
This paper is a part of the first author's Ph.D. thesis under the supervision of
the second author and Iosif Polterovich. The authors are grateful to A.\,Doan
for useful discussions and to Dylan Cant for helpful comments on a previous
version of this project. We also thank Alexandre Girouard, Bernard Helffer, and
Zeev Rudnick for their interest in this project. Finally, we thank the
organizers of the spectral geometry seminar at the Isaac Newton Institute and
the organizers of the Oxford topological data analysis seminar for the
opportunity to present a preliminary version of our work. EM is partially
supported by the ISM Graduate Scholarship and the FRQNT graduate scholarship. ES
is partially supported by an NSERC Discovery grant, by an FRQNT Teams grant, by
a Simons Targeted Research Group in Hamiltonian Dynamics and Beyond, by the
Courtois chair in fundamental research, and by an anonymous donor at the
Institute for Advanced Study during the Fall semester of 2026. 

\section{Algebraic preliminaries}\label{sec:prealg}

\subsection{Some general notation and conventions}

The sets of all nonnegative integers and nonzero reals are denoted by \(\mathbb
N_0\) and \(\mathbb R^*\) respectively. The cardinality of a set \(S\) is
denoted by \(|S|\). The identity linear map on any vector space \(V\) is denoted
by \(\mathbf 1\) (or \(\mathbf 1_V\) if we want to avoid possible confusion).
Throughout the text and especially in Section~\ref{sec:wendl}, we extensively
use the multi-index notation: if \(\alpha=(\alpha_1,\ldots,\alpha_n)\in\mathbb
N_0^n\), then \(X^\alpha:=X_1^{\alpha_1}\ldots
X_n^{\alpha_n},|\alpha|:=\alpha_1+\ldots+\alpha_n,\alpha!:=\alpha_1!\ldots
\alpha_n!,(-1)^\alpha:=(-1)^{|\alpha|}\), etc. If \(\mathcal A\) is an algebra
and \(f\in\mathcal A\), then we often denote the multiplication operator
\(\mathcal A\to\mathcal A,a\mapsto fa\) simply by \(f\). A similar notational
convention is used for group actions. For any vector space \(V\) its symmetric
square \(S^2 V\) is always considered as a subspace of the tensor square
\(V^{\otimes 2}=V\otimes V\), and we denote \(v \odot w := \frac12(v\otimes w+w
\otimes v)\). For example, if \(f,g\colon V\to W\) are linear maps, then
\[
(f\odot g)(u\odot v)=\frac{1}4\bigl(f(u)\otimes g(v)+f(v)\otimes g(u)+g(u)
\otimes f(v)+g(v)\otimes f(u)\bigr).
\]

\subsection{Filtered vector spaces}\label{sec:filt}

Let \(V\) be a vector space. A \textbf{filtration} \(V^\bullet\) on \(V\) is a
descending chain of subspaces
\[
V=V^0\supset V^1\supset\ldots\supset V^j\supset\ldots,\qquad j\in\mathbb N_0.
\]
Sometimes it is useful to assume that \(V^j:=V\) for all integers \(j<0\). A
filtration is called \textbf{finite} if there exists \(j_0\) such that
\(V^j=\{0\}\) for all \(j>j_0\) and \textbf{infinite} otherwise. For an infinite
filtration we put \(V^\infty:=\bigcap_{j\in\mathbb N_0} V^j\). The
\textbf{order} \(\ord\) with respect to the filtration \(V^\bullet\) is defined
by
\[
\ord\colon V\to\mathbb N_0\cup\{\infty\},\qquad
\ord v:=\sup\{j\in\mathbb N_0\colon v\in V^j\}.
\]
In particular, \(\ord 0=\infty\). The \textbf{associated graded space} \(\gr V\)
is defined by
\[
\gr V:=\bigoplus_{j\in\mathbb N_0}\gr^j V,\quad\text{where}\quad
\gr^j V:=V^j/V^{j+1}.
\]
The equivalence class of a vector \(v\in V^j\) in \(\gr^j V\) is denoted by
\(v^j\). The \textbf{lowest order term map} \(\lt\) is given by
\[
\lt\colon\ord^{-1}(\mathbb N_0)\to\gr V,\qquad
\lt v:=v^j\in\gr^j V,\quad\text{where}\quad j=\ord v.
\]
Let us discuss some ways to construct new filtrations from given ones. If
\(W\subset V\) is a subspace, then the induced filtration \(W^\bullet\) on \(W\)
is defined by
\[
W^j:=W\cap V^j,\qquad j\in\mathbb N_0,
\]
and the induced filtration \((V/W)^\bullet\) on \(V/W\) is defined by
\[
(V/W)^j:=V^j/(V^j\cap W).
\]
An important particular case of the latter construction is when \(W=V^{j_0}\)
for some \(j_0\in\mathbb N_0\cup\{\infty\}\). In this case the resulting
filtration on \(V/V^{j_0}\) is finite whenever \(j_0<\infty\).

If \(V,W\) are two filtered vector spaces with filtrations
\(V^\bullet,W^\bullet\), then there is a natural filtration on their tensor
product \(V\otimes W\) defined by
\[
(V\otimes W)^j:=\sum_{j'=0}^j V^{j'}\otimes W^{j-j'}.
\]
Combining these constructions, one can define natural filtrations on other
vector spaces. For example, given a filtration on \(V\), a natural filtration on
\(S^2 V\) is defined by restricting the tensor product filtration on \(V\otimes
V\) onto symmetric tensors.

Finally, let us discuss maps between filtered vector spaces. Let \(V,W\) be two
filtered vector spaces with filtrations \(V^\bullet,W^\bullet\). A linear map
\(f\colon V\to W\) is of \textbf{order at least} \(k\) if \(f(V^j)\subset
W^{j+k}\) for all \(j\in\mathbb N_0\). In this case \(f\) induces a map
\begin{equation}\label{eq:filt-grkf}
\gr_k f\colon\gr V\to\gr W,\qquad
\gr_k f(v^j):=f(v)^{j+k},\quad j\in\mathbb N_0,v\in V^j.
\end{equation}
This map is nonzero provided that \(f\) has order \emph{exactly} \(k\), that is,
\(f(V^j)\not\subset W^{j+k+1}\) for some \(j\in\mathbb N_0\). It is easy to
check directly that the map \(\gr_k f\) comes with the natural inclusion
\(\gr\ker f\hookrightarrow\ker\gr_k f\) and the natural projection \(\coker\gr_k
f\twoheadrightarrow\gr\coker f\).

\subsection{Jet spaces}\label{sec:jets}

Let \(E\) be a vector bundle of rank \(r\) over \(M\). Fix a point \(x\in M\).
The space \(C_x^\infty(E)\) of germs of smooth sections of \(E\) at \(x\)
carries the natural \textbf{vanishing order filtration}
\(C_x^\infty(E)^\bullet\) defined for \(j\in\mathbb N_0\) by 
\[
C_x^\infty(E)^j:=\{u\in C_x^\infty(E)\colon\text{all the derivatives of \(u\) at
\(x\) up to order \(j-1\) vanish}\}.
\]

For any \(l\in\mathbb N_0\) the \textbf{\(l\)-jet space} at \(x\) is the space
\[
J_x^l(E):=C_x^\infty(E)/C_x^\infty(E)^{l+1}
\]
equipped with a finite filtration \(J_x^l(E)^\bullet\) induced from the
vanishing order filtration. As usual, if \(E=M\times\mathbb R\), then we write
just \(J_x^l(E)=J_x^l(M)\). Let
\[
j_x^l\colon C_x^\infty(E)\to J_x^l(E)
\]
denote the quotient map.

Once a local coordinate system \(x_1,\ldots,x_n\) in a neighbourhood of \(x\)
and a local trivialization of \(E\) are chosen, the Taylor expansion identifies
the space \(J_x^l(E)\) with the space \(\mathbb R_{\le l}[x_1,\ldots,x_n]^r\),
where \(\mathbb R_{\le l}[x_1,\ldots,x_n]\) is the space of polynomials of
degree at most \(l\). Although the isomorphism depends on the choice of local
coordinates, for any nonzero \(u\in J_x^l(E)\) the lowest order term of the
Taylor polynomial of \(u\) gives a well-defined element in \(S^j T_x^*M\otimes
E_x\) with \(j=\ord u\). In other words, there are natural isomorphisms
\[
\gr^j J_x^l(E)\cong S^j T_x^*M\otimes E_x,\qquad
\gr J_x^l(E)\cong S^{\le l}T_x^*M\otimes E_x.
\]
These isomorphisms are used implicitly in the sequel.

The multiplication of sections of \(E\) by functions descends to jets. More
precisely, there is a well-defined product map
\[
J_x^{l_1}(M)\times J_x^{l_2}(E)\to J_x^{\min(l_1,l_2)}(E),\qquad
l_1,l_2\in\mathbb N_0.
\]
Of course, in local coordinates this operation is given by the multiplication of
Taylor polynomials followed by the truncation.

It is useful to know that a \(G\)-equivariant section can have any jet at a
regular point.

\begin{claim}\label{cl:jets-regpoint}
Let \(x\in M\) be regular, \(E\) be a \(G\)-equivariant vector bundle, and
\(u\in J_x^l(E)\) for some \(l\in\mathbb N_0\). Suppose that \(G_0\) acts
trivially on \(E\). Then there exists \(\tilde u\in C_G^\infty(E)\) such that
\(j_x^l\tilde u=u\).
\end{claim}

\begin{proof}
Since \(x\) is regular, there exists a coordinate neighbourhood \(U\) of
\(x\) such that for every \(h\in G\setminus G_0\) we have \(h(U)\cap
U=\emptyset\). Let \(x_1,\ldots,x_n\) be local coordinates in \(U\) centered at
\(x\) and let \(p(x_1,\ldots,x_n)\) be the Taylor polynomial of the \(l\)-jet
\(u\). Take a cutoff function \(\eta\in C_c^\infty(U)\), that is, \(\eta\equiv
0\) on \(M\setminus U\) and \(\eta\equiv 1\) on some smaller neighbourhood
\(V\subset U\) of \(x\). Then \(\eta p\in C^\infty(E)\) and the section given
by 
\[
\tilde u:={\frac 1{|G_0|}} \sum_{h\in G} h^*(\eta p)
\]
is \(G\)-equivariant and has \(l\)-jet at \(x\) equal to \(u\) as desired.
\end{proof}

\smallskip

Let \(x\in M\) and \(k\in\mathbb N_0\). A \textbf{symbol} of order \(k\) at
\(x\) is an element \(\sigma\in S^k T_x M\). Recall that any tangent vector
\(\xi\in T_x M\) canonically defines a derivation \(\partial_\xi\) on
\(S^\bullet T_x^*M\). Hence any symbol \(\sigma\in S^k T_x M\) (also
canonically) gives rise to a formal differential operator \(\hat\sigma\) on
\(S^\bullet T_x^*M\). {For each \(l\in\mathbb N_0\) we put
\(\hat\sigma^l:=\hat\sigma|_{S^{\le l}T_x^*M}\). This is a linear map
\(\hat\sigma^l\colon S^{\le l}T_x^*M\to S^{\le l-k}T_x^*M\).} The symbol
\(\sigma\) is called \textbf{elliptic} if \(\sigma(\zeta^{\otimes k})\ne 0\) for
all \(\zeta\in T_x^*M\setminus\{0\}\). When \(\sigma\) is a symbol of a real
scalar formally self-adjoint differential operator, this is possible only for
\(k\) even, in which case \(\sigma(\zeta^{\otimes k})\) has the same sign for
all \(\zeta\ne 0\). In this paper we adhere to the convention
\((-1)^{k/2}\sigma(\zeta^{\otimes k})>0\). This assures that the spectrum of a
differential operator having such symbol is bounded from below, not from above.

Let \(l\in\mathbb N_0\). Any differential operator \(\mathbf P\colon
C^\infty(M)\to C^\infty(M)\) of order \(k\) descends to a linear map
\begin{equation}\label{eq:jets-Pl}
\mathbf P^l\colon J_x^{l+k}(M)\to J_x^l(M)
\end{equation}
of order \(-k\). Moreover, it is easy to see that the induced map
\[
\gr_{-k}\mathbf P^l\colon S^{\le l+k}T_x^*M\to S^{\le l}T_x^*M
\]
coincides with \(\hat\sigma^{l+k}\), where \(\sigma=\sigma_x(\mathbf P)\in S^k
T_x M\) is the principal symbol of \(\mathbf P\) at \(x\). Suppose that
\(\mathbf P\) is elliptic. By definition, this means that the symbol \(\sigma\)
is elliptic. In this case it can be shown
(see~\cite[Proposition~1.6.20]{DW2023}) that both \(\hat\sigma^{l+k}\) and
\(\mathbf P^l\) are surjective for all \(l\in\mathbb N_0\). In particular,
\begin{equation}\label{eq:jets-dimkersigma}
\dim\ker\hat\sigma^{l+k}=\dim S^{l+k}T_x^* M-
\dim S^l T_x^* M=\binom{n+l+k}n-\binom{n+l}n,
\end{equation}
as well as
\begin{equation}\label{eq:jets-dimkerP}
\dim\ker\mathbf P^l=\dim J_x^{l+k}(M)-\dim J_x^l(M)=\binom{n+l+k}n-\binom{n+l}n.
\end{equation}
The map \(\hat\sigma^{l+k}\) is surjective for all \(l\in\mathbb N_0\) if and
only if the map \(\hat\sigma\) is surjective or, equivalently, it has a right
inverse \(\widehat R\). Since \(\hat\sigma\) is a map of degree \(-k\), we can
always choose the map \(\widehat R\) to be of degree \(k\) and in this paper we
consider only such right inverses (in other words, we require \(\widehat R\) to
be a morphism of graded vector spaces). We will need the fact that the map
\(\mathbf P^l\) has a right inverse that is `compatible' with any given right
inverse \(\widehat R\) of \(\hat\sigma\).

\begin{claim}\label{cl:jets-compR}
For any right inverse \(\widehat R\) to \(\hat\sigma\) and \(l\in\mathbb N_0\)
there exists a right inverse \(R\) to \(\mathbf P^l\) such that \(R\) has order
\(k\) and \(\gr_k R=\widehat R^l\).
\end{claim}

\begin{proof}
We define \(R\) on \(J_x^l(M)^j\) by descending induction in \(j\). Since
\(J_x^l(M)^{l+1}=\{0\}\), for \(j=l+1\) there is nothing to do. Suppose that
\(R\) is already defined, and satisfies \(\gr_k R=\widehat R^l\), on
\(J_x^l(M)^{j+1}\) and let us extend it to \(J_x^l(M)^j\). Let
\(v_1,\ldots,v_p\in J_x^l(M)^j\) be such that \(v_1^j,\ldots v_p^j\) form a
basis of \(S^j T_x^*M\). For each \(i=1,\ldots,p\) choose \(u_i\in
J_x^{l+k}(M)^{j+k}\) such that \(u_i^{j+k}=\widehat R(v_i^j)\in S^{j+k}T_x^*M\).
Then \(\hat\sigma^{l+k}(u_i^{j+k})=v_i^j\) and
\[
(\mathbf P^l(u_i)-v_i)^j=\mathbf P^l(u_i)^j-v_i^j=
\gr_{-k}\mathbf P^l(u_i^{j+k})-v_i^j=\hat\sigma^{l+k}(u_i^{j+k})-v_i^j=0,
\]
where we used the definition~\eqref{eq:filt-grkf} of \(\gr_{-k}\) and that
\(\gr_{-k}\mathbf P^l=\hat\sigma^{l+k}\). It follows that \(\mathbf
P^l(u_i)-v_i\in J_x^l(M)^{j+1}\) and \(R(\mathbf P^l(u_i)-v_i)\) is already
defined. Put
\[
R(v_i):=u_i-R(\mathbf P^l(u_i)-v_i)\in J_x^{l+k}(M)^{j+k}.
\]
It is easy to check that defined this way, \(R\) is exactly what we need.
\end{proof}

\section{Analytic preliminaries}\label{sec:preanal}

\subsection{Walls in the space of symmetric Fredholm operators}\label{sec:walls}

In this section we follow~\cite[Section~3.2.2]{Wen2016}. Let \(\mathcal H\) be a
real Hilbert space with inner product \(\langle\cdot,\cdot\rangle_{\mathcal H}\)
equipped with an action of a finite group \(G\) by isometries. Let \(\mathcal
D\subset\mathcal H\) be a \(G\)-equivariant linear subspace with a Banach space
structure such that the inclusion \(\mathcal D\hookrightarrow\mathcal H\) is a
continuous 
operator. In the space \(\mathcal L_G(\mathcal D,\mathcal H)\) of bounded linear
\(G\)-equivariant operators there is a closed subspace of symmetric operators
defined by
\[
\mathcal L_G^\sym(\mathcal D,\mathcal H):=\{\mathbf T\in\mathcal L_G(\mathcal
D,\mathcal H)\colon\langle x,\mathbf Ty\rangle_{\mathcal
H}=\langle\mathbf Tx,y\rangle_{\mathcal H},\quad x,y\in\mathcal D\}.
\]
Let \(\Fred_G(\mathcal D,\mathcal H)\subset\mathcal L_G(\mathcal D,\mathcal H)\)
denote the open subset of \(G\)-equivariant Fredholm operators. For any
\(\mathbf T\in\Fred_G(\mathcal D,\mathcal H)\) both \(\ker\mathbf T\) and
\(\coker\mathbf T\) are finite-dimensional representations of \(G\). The
\textbf{\(G\)-index} of the operator \(\mathbf T\) is defined by
\[
\ind_G\mathbf T:=[\ker\mathbf T]-[\coker\mathbf T]\in\mathbb Z^N,
\]
where \(N\) is the number of irreducible \(G\)-representations as in
Section~\ref{sec:mainres}. Now for any \(\bar\imath\in\mathbb Z^N,\bar
k\in\mathbb N_0^N\) we put
\begin{align*}
\Fred_G^{\bar\imath}(\mathcal D,\mathcal H)&:=
\{\mathbf T\in\Fred_G(\mathcal D,\mathcal H)\colon\ind_G\mathbf
T=\bar\imath\},\\
\Fred_G^{\bar\imath,\bar k}(\mathcal D,\mathcal H)&:=
\{\mathbf T\in\Fred_G(\mathcal D,\mathcal H)\colon\ind_G\mathbf
T=\bar\imath,[\ker\mathbf T]=\bar k\}
\end{align*}
and
\[
\Fred_G^\sym(\mathcal D,\mathcal H):=
\Fred_G^{\bar 0}(\mathcal D,\mathcal H)\cap\mathcal L_G^\sym(\mathcal D,\mathcal
H),\] \[
\Fred_G^{\sym,\bar k}(\mathcal D,\mathcal H):=
\Fred_G^{\bar 0,\bar k}(\mathcal D,\mathcal H)\cap\mathcal L_G^\sym(\mathcal
D,\mathcal H),
\]
where \(\bar 0:=(0,\ldots,0)\in\mathbb N_0^N\).

Fix an operator \(\mathbf T_0\in\Fred_G^\sym(\mathcal D,\mathcal H)\). Applying
the same argument as in~\cite[Lemma~3.8]{Wen2016} it is easy to see that there
exist splittings of \(\mathcal D\) and \(\mathcal H\) into closed linear
\(G\)-invariant subspaces
\[
\mathcal D=V\oplus K,\qquad\mathcal H=W\oplus K,
\]
where \(K=\ker\mathbf T_0,W=\im\mathbf T_0\), and \(V=W\cap\mathcal D\). Using
these splittings, every operator \(\mathbf T\in\mathcal L_G(\mathcal D,\mathcal
H)\) can be written in block form as
\[
\mathbf T=\begin{pmatrix}
\mathbf A & \mathbf B\\
\mathbf C & \mathbf D
\end{pmatrix},
\]
with \(\mathbf A,\mathbf B,\mathbf C,\mathbf D\) being \(G\)-equivariant and
\(\mathbf T_0\) in this form looks like
\[
\mathbf T_0=\begin{pmatrix}
\mathbf A_0 & 0\\
0 & 0
\end{pmatrix}
\]
for some \(G\)-equivariant Banach space isomorphism \(\mathbf A_0\colon V\to
W\). Let \(\mathcal O\) be a neighbourhood of \(\mathbf T_0\) in
\(\Fred_G^\sym(\mathcal D,\mathcal H)\) such that for all \(\mathbf T\in\mathcal
O\) the block \(\mathbf A\) is invertible. Then there are well-defined
\(G\)-equivariant maps
\begin{align*}
\Phi&\colon\mathcal O\to\Sym_G(K),&\qquad
\Phi(\mathbf T)&:=\mathbf D-\mathbf{CA}^{-1}\mathbf B,\\
\Psi&\colon\mathcal O\to\mathcal L_G(\mathcal D),&\quad
\Psi(\mathbf T)&:=\begin{pmatrix}
\mathbf 1 & -\mathbf A^{-1}\mathbf B\\
0 & \mathbf 1
\end{pmatrix},
\end{align*}
where \(\Sym(K)\) is the space of symmetric linear maps on
\((K,\langle\cdot,\cdot\rangle_{\mathcal H})\). Both maps are smooth and the
derivative of \(\Phi\) at \(\mathbf T_0\) is given by
\begin{equation}\label{eq:walls-dPhi}
\d_{\mathbf T_0}\Phi\colon\mathcal L_G^\sym(\mathcal D,\mathcal
H)\to\Sym_G(K),\qquad
\d_{\mathbf T_0}\Phi(H)=\pr_K\circ(H|_K),
\end{equation}
where \(\pr_K\) is the orthogonal projection onto \(K\). Note that for any
\(\mathbf T\in\mathcal O\) the map \(\Psi(\mathbf T)\) is invertible with the
inverse given explicitly by
\[
\Psi(\mathbf T)^{-1}=\begin{pmatrix}
\mathbf 1 & \mathbf A^{-1}\mathbf B\\
0 & \mathbf 1
\end{pmatrix}.
\]
Since
\[
\mathbf T\Psi(\mathbf T)=\begin{pmatrix}
\mathbf A & 0\\
\mathbf C & \Phi(\mathbf T)
\end{pmatrix},
\]
where the map \(\mathbf A\) is invertible, it follows that \(\Psi(\mathbf T)\)
maps \(\ker\Phi(\mathbf T)\subset K=\ker\mathbf T_0\) \(G\)-equivariantly and
isomorphically onto \(\ker\mathbf T\). In particular, Schur's lemma implies that
\begin{equation}\label{eq:walls-kerusc}
[\ker\mathbf T]\le [\ker\mathbf T_0]\quad\text{with equality if and only
if}\quad\Phi(\mathbf T)=\mathbf 0,
\end{equation}
where `\(\le\)' is the product partial order on \(\mathbb N_0^N\), that is, for
\(\overline m,\overline n\in\mathbb N_0^N\) we write \(\overline m\le\overline
n\) if \(m_j\le n_j\) for all \(j=1,\ldots,N\).

It is a consequence of Schur's lemma that the endomorphism algebra of any
irreducible representation is a division algebra over \(\mathbb R\), that is, it
is isomorphic to \(\mathbb R,\mathbb C\), or \(\mathbb H\). In particular, if
\(W_1,\ldots,W_N\) are all irreducible \(G\)-representations, then
\(\alpha_j:=\dim\End_G(W_j)\in\{1,2,4\}\) for each \(j=1,\ldots,N\). This
corresponds to irreducible representations of \textbf{real}, \textbf{complex},
and \textbf{quaternionic} types respectively. The following lemma
extends~\cite[Proposition~3.9]{Wen2016} to an equivariant setting,
cf.~\cite[Theorem~1.4.6 and Section~1.A]{DW2023}

\begin{lemma}\label{lem:walls-main}
Let \(Z\) be a Banach manifold and
\[
\mathbf T\colon Z\to\Fred_G^\sym(\mathcal D,\mathcal H),\qquad
z\mapsto\mathbf T_z
\]
be a smooth map. Let \(z\in Z,K:=\ker\mathbf T_z,\bar k:=[\ker\mathbf
T_z]\), and
\[
Z^{\bar k}:=\{z'\in Z\colon\mathbf T_{z'}\in\Fred_G^{\sym,\bar
k}(\mathcal D,\mathcal H)\}.
\] 
Suppose that the map
\begin{equation}\label{eq:walls-Lambda}
\Lambda_z\colon T_z Z\to\Sym_G(K),\qquad
\Lambda_z(\zeta):=\pr_K\circ (\d_z\mathbf T(\zeta)|_K)
\end{equation}
is surjective. Then there exists a neighbourhood \(\mathcal U\) of \(z\) in
\(Z\) such that \(\mathcal U\cap Z^{\bar k}\) is a Banach submanifold in
\(\mathcal U\) of codimension
\begin{equation}\label{eq:walls-codim}
\codim(\mathcal U\cap Z^{\bar k})=\dim\Sym_G(K)=
\sum_{j=1}^N k_j+\alpha_j\cdot\frac{k_j(k_j-1)}2.
\end{equation}
In particular, \(\codim(\mathcal U\cap Z^{\bar k})=1\) if and only if
\(\ker\mathbf T_z\) is an irreducible \(G\)-representation.
\end{lemma}

\begin{proof}
Since by~\eqref{eq:walls-dPhi} and chain rule the map \(\Lambda_z\) is nothing
but \(\d_z(\Phi\circ\mathbf T)\) for \(\Phi\) defined as above (with \(\mathbf
T_z\) in place of \(\mathbf T_0\)), it follows from the implicit function
theorem that there exists a neighbourhood \(\mathcal U\) of \(z\) in \(Z\) such
that \(\mathcal U\cap (\Phi\circ\mathbf T)^{-1}(0)\) is a submanifold in
\(\mathcal U\) of codimension \(\dim\Sym_G(K)\). Since the explicit formula for
\(\dim\Sym_G(K)\) is another consequence of Schur's lemma, it only remains to
show that \((\Phi\circ\mathbf T)^{-1}(0)=Z^{\bar k}\). For this observe that
by~\eqref{eq:walls-kerusc} we may assume that for every \(z'\in\mathcal U\) we
have
\[
[\ker\mathbf T_{z'}]\le[\ker\mathbf T_z]=\bar k
\]
with equality if and only if \(\Phi(\mathbf T_{z'})=0\), that is, if \({z'}\in
(\Phi\circ\mathbf T)^{-1}(0)\). This implies the claim.
\end{proof}

\subsection{Floer \texorpdfstring{\(\epsilon\)}{ε}-spaces}\label{sec:floer}

Floer \(\epsilon\)-spaces provide a device designed to work with Fr\'e\-chet
manifolds in the same manner as if they were Banach manifolds. Our exposition
follows~\cite[Appendix~B]{Wen2016},~\cite{Wen2021},
and~\cite[Section~5.4]{Wen2023} but is slightly more general. More precisely,
the underlying Fr\'echet space of our Fr\'echet manifold is not necessarily the
space of \(C^\infty\) sections of some vector bundle (for this reason, we speak
about Floer \(\epsilon\)-spaces, not Floer \(C_\epsilon\)-spaces). Although we
do not use this extension, i.e. all Fr\'echet manifolds we encounter in this
text are modeled on the spaces of \(C^\infty\) sections of some vector bundle,
we hope that it simplifies the exposition and may be useful for other
applications.

Let \(P\) be a Fr\'echet space. Suppose that the topology on \(P\) can be
defined by a countable family of \emph{norms} \((\|\cdot\|_k)_{k\in\mathbb
N_0}\) such that for each \(k\), the completion \(P_k\) of \(P\) w.r.t. the norm
\(\|\cdot\|_k\) is a separable Banach space. Let
\(\epsilon=(\epsilon_k)_{k\in\mathbb N_0}\) be a sequence of reals such that
\(\epsilon_k>0\) for all \(k\) and \(\lim\limits_{k\to\infty}\epsilon_k=0\).
Then for any \(\eta\in P\) the \textbf{(Floer) \(\epsilon\)-norm} of \(\eta\) is
defined by
\[
\|\eta\|_\epsilon:=\sum_{k=0}^\infty\epsilon_k\|\eta\|_k.
\]
Put
\[
P_\epsilon:=\{\eta\in P\colon \|\eta\|_\epsilon<\infty\}.
\]

\begin{proposition}
For each \(\epsilon\in\mathcal E\), \(P_\epsilon\) is a separable Banach space.
\end{proposition}

\begin{proof}
That \(P_\epsilon\) is complete can be proved by a standard argument,
see~\cite[Theorem~B.2]{Wen2016} for the details. The separability of
\(P_\epsilon\) can be proved by the same trick as
in~\cite[Theorem~B.5]{Wen2016}. Let \(P_{(\epsilon)}\) be the space of sequences
\[
\xi=(\xi_0,\xi_1,\ldots)\in\prod_{k=0}^\infty P_k
\]
satisfying
\[
\|\xi\|_{(\epsilon)}:=\sum_{k=0}^\infty\epsilon_k\|\xi_k\|_k<\infty.
\]
Since \((P_k,\|\cdot\|_k)\) is a separable Banach space for each \(k\), it
contains a countable dense subset \(A_k\). Then it is easy to check that
\[
A:=\{(\xi_0,\xi_1,\ldots,\xi_N,0,0,\ldots)\in P_{(\epsilon)}\colon
N\in\mathbb N_0,\xi_k\in A_k\}
\]
is countable and dense in \((P_{(\epsilon)},\|\cdot\|_{(\epsilon)})\). Since the
map
\[
P_\epsilon\hookrightarrow P_{(\epsilon)},\qquad\eta\mapsto (\eta,\eta,\ldots)
\]
is an isometric embedding of metric spaces and a subspace of a separable metric
space is separable, we get the claim.
\end{proof}

Now let \(\mathcal P\) be a Fr\'echet manifold modeled on the space \(P\). Take
any point \(p_\rf\in\mathcal P\) (called a \textbf{reference point}). Let
\(\mathcal U\) be a coordinate neighbourhood of \(p_\rf\). Since
\(T_{p_\rf}\mathcal P\cong P\), there exists a neighbourhood \(\mathcal O\) of
the origin in \(T_{p_\rf}\mathcal P\) diffeomorphic to \(\mathcal U\). Fix a
diffeomorphism and denote it by \(\exp_{p_\rf}\colon\mathcal O\to\mathcal U\).
Then the \textbf{(Floer) \(\epsilon\)-space} is defined by
\[
\mathcal P_\epsilon:=\{\exp_{p_\rf}(\eta)\mid\eta\in\mathcal O,
\|\eta\|_\epsilon<\infty\}.
\]
Defined this way, \(\mathcal P_\epsilon\) is a smooth Banach manifold modeled on
\(P_\epsilon\) and covered by a single chart. Since the embedding
\(P_\epsilon\hookrightarrow P\) is continuous, the natural inclusion \(\mathcal
P_\epsilon\hookrightarrow\mathcal P\) is continuous as well.

Let \(\mathcal E\) be the set of all sequences \((\epsilon_k)_{k\in\mathbb
N_0}\) such that \(\epsilon_k>0\) for all \(k\) and
\(\lim\limits_{k\to\infty}\epsilon_k=0\). A preorder \(\prec\) on \(\mathcal E\)
is defined by
\[
\epsilon\prec\epsilon'\qquad\Leftrightarrow\qquad
\limsup_{k\to\infty}(\epsilon_k/\epsilon_k')<\infty.
\]

\begin{definition}[{\cite[Definition~5.28]{Wen2023}}]
A statement depending on an element \(\epsilon\in\mathcal E\) is said to be true
\textbf{for all \(\epsilon\) with sufficiently rapid decay} if there exists
\(\epsilon_0\in\mathcal E\) such that the statement holds for each
\(\epsilon\in\mathcal E\) with \(\epsilon\prec\epsilon_0\).
\end{definition}

The following lemma is proved in the same way as~\cite[Lemma~5.29]{Wen2023}.

\begin{lemma}\label{lem:floer-srd}
The norms \(\|\cdot\|_\epsilon\) enjoy the following properties:
\begin{enumerate}[label=(\arabic*)]
\item{if \(\epsilon\prec\epsilon'\), then the norm \(\|\cdot\|_\epsilon\)
is dominated by \(\|\cdot\|_{\epsilon'}\);}
\item{for every \(\eta\in P\) we have \(\|\eta\|_\epsilon<\infty\) for all
\(\epsilon\) with sufficiently rapid decay;}\label{it:floer-chooseeps}
\item{every at most countable subset of \(\mathcal E\) has a lower bound in
\(\mathcal E\) with respect to the preorder \(\prec\).}\label{it:floer-lb}
\end{enumerate}
\end{lemma}

\subsection{\texorpdfstring{\(C^\infty\)}{C∞}-subvarieties}\label{sec:cinfsub}

Of course, the core of the proof of Theorem~\ref{th:outline-schrod} is some kind
of Sard-Smale theorem. The appropriate version is the \emph{Sard-Smale theorem
for \(C^\infty\)-subvarieties}, see~\cite[Appendix~C]{Wen2023}.

\begin{definition}\label{def:cinfsub-subj}
Let \(X\) be a Banach manifold, \(Y\subset X\), and \(c\in\mathbb N\). We say
that \(Y\) is a \textbf{\(C^\infty\)-subvariety} of \textbf{codimension at least
\(c\)} if for any \(x\in Y\) there exists a neighbourhood \(U\) of \(x\) in
\(X\), a finite-dimensional vector space \(W\), and a smooth map \(f\colon U\to
W\) such that \(Y\cap U=f^{-1}(0)\) and \(\rk\d_x f\ge c\).
\end{definition}

The following lemma is a version of the Sard-Smale theorem used in the sequel,
see~\cite[Proposition C.3]{Wen2023}.

\begin{lemma}[Sard-Smale theorem for \(C^\infty\)-subvarieties]
\label{lem:cinfsub-sardsmale}
Let \(X\) and \(Z\) be separable and metrizable Banach manifolds, \(\pi\colon
X\to Z\) be a smooth Fredholm map, and \(Y\subset X\) be a
\(C^\infty\)-subvariety of codimension at least \(c\in\mathbb N\). Let
\(Z_\pi^\reg\subset Z\) be the set of regular values of \(\pi\). Then there
exists a Baire subset \(Z_Y^\reg\subset Z\) such that for all \(z\in
Z_\pi^\reg\cap Z_Y^\reg\) the set \(Y\cap\pi^{-1}(z)\) is a
\(C^\infty\)-subvariety of codimension at least \(c\) in \(\pi^{-1}(z)\).
\end{lemma}

We record the following important special case of
Lemma~\ref{lem:cinfsub-sardsmale}, see~\cite[Corollary C.4]{Wen2023}.

\begin{corollary}
In the setting of Lemma~\ref{lem:cinfsub-sardsmale}, if \(\ind\d_x\pi<c\) for
all \(x\in Y\), then \(Y\cap\pi^{-1}(z)=\emptyset\) for generic \(z\in Z\).
\end{corollary}

\subsection{Smooth families of \texorpdfstring{\(G\)}{G}-equivariant
differential operators}\label{sec:fam}

Let \(\Diff_G^k(M)\) denote the space of \(G\)-equivariant scalar linear
differential operators on \(M\) of order at most \(k\).

\begin{definition}\label{def:fam-subj}
Let \(\mathcal P\) be a Fr\'echet manifold modeled on a Fr\'echet space \(P\) as
in Section~\ref{sec:floer} and \(k\in\mathbb N\). A \textbf{smooth family of
\(G\)-equivariant scalar differential operators} on \(M\) of order \(k\)
parametrized by \(\mathcal P\) is a map
\[
\mathbf D\colon\mathcal P\to\Diff_G^k(M),
\]
where for each \(p\in\mathcal P\), \(\mathbf D_p:=\mathbf D(p)\) is a linear
elliptic formally self-adjoint (with respect to the Riemannian metric \(g_\fix\)
on \(M\)) \(G\)-equivariant differential operator of order \(k\) on \(M\) such
that any point \(x\in M\) has a neighbourhood \(U\subset M\) with the property
that the coefficients of \(\mathbf D_p\) are smooth when considered as functions
on the Fr\'echet manifold \(U\times\mathcal P\).
\end{definition}

\begin{remark}
Note that \(\mathbf D\) restricts to a smooth family on \(U \times P_\epsilon\)
for an \(\epsilon\)-space \(P_\epsilon\) (see Section~\ref{sec:floer}) contained
in \(P\).
\end{remark}

\begin{remark}[{cf.~\cite[footnote on p.~4]{DW2023}}]
Clearly, any family of \(G\)-equivariant differential operators \(\mathbf
D\colon\mathcal P\to\Diff_G^k(M)\) extends to a smooth map
\[
\mathbf D\colon\mathcal P\to\Fred_G^\sym(\mathcal D,\mathcal H),
\]
where \(\mathcal D=H^k(M,\d v_{g_\fix}),\mathcal H=L^2(M,\d v_{g_\fix})\). One
may ask if it is safe to replace the smoothness assumption in
Definition~\ref{def:fam-subj} by this weaker one. In fact, part of the argument
(namely, Section~\ref{sec:walls} and
Propositions~\ref{pr:petri-main},\ref{pr:petri2-main} relying on it) still works
in this case. However, in
Sections~\ref{sec:petrigen},\ref{sec:rkest2},\ref{sec:rkest3} we use that
\begin{equation}\label{eq:fam-clairaut}
\d_p J_x^l\mathbf D=J_x^l\d_p\mathbf D,\qquad p\in\mathcal P,\quad x\in M,\quad
l\in\mathbb N_0.
\end{equation}
In general, this fails with a weaker smoothness assumption as above.
Indeed, what equation~\eqref{eq:fam-clairaut} essentially says is that for any
coefficient of \(\mathbf D\), its partial derivatives in \(p\in\mathcal P\) and
\(x\in M\) commute. This is known to be not true in general and, in fact, any
counterexample (e.g., Peano function) can be used to construct a family
\(\mathbf D\) violating~\eqref{eq:fam-clairaut}. This does not happen if the
coefficients of \(\mathbf D_p\) are smooth \emph{jointly} in \(x\) and \(p\), as
it is required in Definition~\ref{def:fam-subj}.
\end{remark}

In this paper we deal with the following smooth families of \(G\)-equivariant
differential operators.

\begin{example}[Schr\"odinger operators]
Let \(\mathbf L\) be a linear elliptic formally self-adjoint \(G\)-equivariant
differential operator on \(M\) of order \(k\) and \(\mathcal P\) be an open
subset in \(C_G^\infty(M)\). Then
\[
\mathbf D\colon\mathcal P\to\Diff_G^k(M),\qquad
\mathbf D_f:=\mathbf L+f
\]
is a smooth family of \(G\)-equivariant differential operators on \(M\).
\end{example}

\begin{example}[metric-dependent operators]\label{ex:fam-metdep}
A (scalar) \textbf{metric-dependent operator} of order \(k\) on \(M\) is a
family of differential operators of the form \((\mathbf L_g)_{g\in\mathcal G}\)
such that \({\bf L}_g\) is \(g\)-self-adjoint. A metric dependent operator
\(\mathbf L\) is called \textbf{\(G\)-equivariant} if for any \(g\in\mathcal
G^G\), the operator \(\mathbf L_g\) is \(G\)-equivariant.

Now let \(\mathbf L\) be a \(G\)-equivariant metric-dependent differential
operator and \(\mathcal P\subset\mathcal G^G\) be an open subset. Then \(\mathbf
L\) itself is \emph{not} a smooth family of \(G\)-equivariant differential
operators parametrized by \(\mathcal P\) in the sense of
Definition~\ref{def:fam-subj}. Indeed, the operator \(\mathbf L_g\) is
self-adjoint with respect to a \emph{varying} Riemannian metric \(g\) on \(M\)
while Definition~\ref{def:fam-subj} requires it to be self-adjoint with respect
to some \emph{fixed} Riemannian metric \(g_\fix\).

However, if we put
\[
\mathbf D\colon\mathcal P\to\Diff_G^k(M),\qquad
\mathbf D_g:=\rho(g)^{\frac12}\circ\mathbf L_g\circ\rho(g)^{-\frac12},
\]
where \(\rho(g):=\d v_g/\d v_{g_\fix}\in C_G^\infty(M)\) is the density of the
Riemannian volume measure \(\d v_g\) with respect to the Riemannian volume
measure \(\d v_{g_\fix}\), then \(\mathbf D\) is a smooth family of
\(G\)-equivariant differential operators on \(M\) in the sense of
Definition~\ref{def:fam-subj} as soon as the smoothness assumption is satisfied.
The Laplace-Beltrami operator is the classical example of a metric-dependent
operator.
\end{example}

\begin{example}[conformally covariant operators]\label{ex:fam-cc}
A metric dependent operator \(\mathbf L\) is called \textbf{conformally
covariant} of \textbf{bidegree} \((a,b)\) for \(a,b \in \mathbb R\) if for any
\(g_1,g_2\in\mathcal G\) and \(f\in C^\infty(M)\) with \(g_2=e^{2f}g_1\) we have
\begin{equation}\label{eq:fam-conf}
\mathbf L_{g_2}=e^{-bf}\circ\mathbf L_{g_1}\circ e^{af}.
\end{equation}

We assume that \(a\ne b\) throughout the paper. This holds for the conformal
Laplacian as well as, to the best of our knowledge, for all known interesting
examples of such operators, see~\cite{Can2014},\cite{CG2026}.

A conformally covariant operator is called \textbf{\(G\)-equivariant} if it is
\(G\)-equivariant as a metric-dependent operator. Let \(\mathcal
C\subset\mathcal G\) be an arbitrary conformal class. As in
Example~\ref{ex:fam-metdep}, any \(G\)-equivariant conformally covariant
operator \(\mathbf L\) can be turned into a smooth family of \(G\)-equivariant
differential operators \(\mathbf D\) parametrized by \(\mathcal G^G\) or
\(\mathcal C^G\) as soon as the smoothness assumption is satisfied, in which
case we say that \(\mathbf L\) \textbf{has smooth coefficients}. As we will see
in Section~\ref{sec:confcov}, for these operators the variational formulas take
a particularly nice form. The conformal Laplacian, Paneitz operator, and, more
generally, GJMS operators are the classical examples of conformally covariant
operators (see~\cite{Can2014} for a review of the literature).
\end{example}

\subsection{Continuity of eigenvalues}\label{sec:cont}

Let \(\mathcal P,\mathbf D\) be as in Definition~\ref{def:fam-subj}. For the
purposes of this paper we need to incorporate the dependence on the eigenvalue
into the family. To this end, we put \(\mathcal M:=\mathcal P\times\mathbb R\)
and
\begin{equation}\label{eq:cont-T}
\mathbf T\colon\mathcal M\to\Diff_G(M),\qquad
\mathbf T_{p,\lambda}:=\mathbf D_p-\lambda.
\end{equation}
Then \(\mathbf T\) extends to a smooth map
\[
\mathbf T\colon\mathcal M\to\Fred_G^\sym(\mathcal D,\mathcal H),
\]
where \(\mathcal D:=W^{k,2}(M,\d v_{g_\fix}),\mathcal H:=L^2(M,\d v_{g_\fix})\).

Let \(p\in\mathcal P\). For any \(\lambda\in\mathbb R\) the vector \(\mathfrak
m(p,\lambda):= [\ker\mathbf T_{p,\lambda}]\) in \(\mathbb N_0^N\) is called the
\textbf{\(G\)-multiplicity vector} of \(\lambda\) and the number
\(m(p,\lambda):= | \mathfrak m(p,\lambda)|\) is called the
\textbf{\(G\)-multiplicity} of \(\lambda\). In particular, \(\lambda\) is an
eigenvalue of \(\mathbf D_p\) whenever \(m(p,\lambda)\ne 0\). Since \(\mathbf
D_p\) is a linear, elliptic, formally self-adjoint scalar operator, its spectrum
is real, discrete, and has the form
\[
\lambda_1(p)\le\lambda_2(p)\le\ldots\le\lambda_\alpha(p)\le\ldots\nearrow
+\infty,
\]
where each eigenvalue is counted with its \(G\)-multiplicity. The aim of this
section is to prove the continuity of eigenvalues in the \(G\)-equivariant
setting.

\begin{proposition}\label{pr:cont-main}
For each \(\alpha\in\mathbb N\)
\begin{enumerate}[label=(\roman*)]
\item
the eigenvalue \(\lambda_\alpha(p)\) is a continuous function in \(p\in\mathcal
P\);\label{it:cont-eigenval}
\item
the \(G\)-multiplicity vector \(\mathfrak m(p,\lambda_\alpha(p))\) is an upper
semicontinuous \(\mathbb N_0^N\)-valued function in \(p\in\mathcal P\) with
respect to the product partial order on \(\mathbb N_0^N\).
\end{enumerate}
\end{proposition}

\begin{proof}
The proof proceeds by standard methods of~\cite{Kat1966}. Fix \(p_0\in\mathcal
P\). First let us show that \(\mathbf D_p\) is uniformly bounded from below near
\(p_0\), i.e. there exists a neighbourhood \(\mathcal U_0\) of \(p_0\) in
\(\mathcal P\) and \(\mu_0\in\mathbb R\) such that \(\mu_0<\lambda_1(p)\) for
each \(p\in\mathcal U_0\). Shifting \(\mathbf D_{p_0}\) if necessary, we may
assume that \(\lambda_1(p_0)>0\). The elliptic estimate implies that there
exists a constant \(C>0\) such that
\[
\|u\|_{H^k}\le C(\|u\|+\|\mathbf D_{p_0}u\|),\qquad
u\in H^k(M,\d v_{g_\fix}),
\]
where \(\|\cdot\|\) stands for the \(L^2\)-norm. On the other hand, it is easy
to see that the smoothness condition of Definition~\ref{def:fam-subj} implies
that for any \(\epsilon>0\) there exists a neighbourhood \(\mathcal U_\epsilon\)
of \(p_0\) in \(\mathcal P\) such that
\[
\|(\mathbf D_p-\mathbf D_{p_0})u\|\le
\frac{\epsilon}C\|u\|_{H^k},\qquad
p\in\mathcal U_\epsilon,\quad u\in H^k(M,\d v_{g_\fix}).
\]
Combining the last two inequalities, we obtain
\[
\|(\mathbf D_p-\mathbf D_{p_0})u\|\le
\epsilon(\|u\|+\|\mathbf D_{p_0}u\|),\qquad
p\in\mathcal U_\epsilon,\quad u\in H^k(M,\d v_{g_\fix}).
\]
For \(\epsilon\) such that \(q:=\epsilon(1+\|\mathbf D_{p_0}^{-1}\|)<1\)
this gives
\[
\|(\mathbf D_p-\mathbf D_{p_0})u\|\le q\|\mathbf D_{p_0}u\|<
\|\mathbf D_{p_0}u\|,\qquad u\in H^k(M,\d v_{g_\fix})\setminus\{0\}.
\]
On the other hand, if \(u\) is an eigenfunction of \(\mathbf D_p\) with
eigenvalue \(\lambda\le 0\), then
\begin{multline*}
\|(\mathbf D_p-\mathbf D_{p_0})u\|^2=
\|(\mathbf D_{p_0}-\lambda)u\|^2=
\|\mathbf D_{p_0}u\|^2-2\lambda\langle\mathbf D_{p_0}u,u\rangle+\\
\lambda^2\|u\|^2\ge\|\mathbf D_{p_0}u\|^2,
\end{multline*}
a contradiction. Hence \(\lambda>0\), i.e. all eigenvalues of \(\mathbf D_p\)
are positive. Moreover,
\[
\lambda\|u\|=\|\mathbf D_p u\|\ge (1-q)\|\mathbf D_{p_0}u\|\ge
\frac{1-q}{\|\mathbf D_{p_0}^{-1}\|}\|u\|,
\]
so \(\mathcal U_0:=\mathcal U_\epsilon\) does the job (see
also~\cite[Theorem~4.11 in Chapter~V]{Kat1966}).


Now fix \(\alpha\in\mathbb N\) and let
\[
\mu_1<\mu_2<\ldots<\mu_\beta<\ldots\nearrow +\infty
\]
be the eigenvalues of \(\mathbf D_{p_0}\) counting \emph{without}
\(G\)-multiplicity. Note that \(\mu_0\le\mu_1\). Let \(\alpha,\beta\in\mathbb
N\) be such that \(\lambda_\alpha(p_0)=\mu_\beta\). For each
\(i=1,\ldots,\beta\), let \(\Gamma_i\subset\mathbb C\) be a circle centered at
\(\mu_i\) so small that \(\mu_i\) is the only eigenvalue inside (the disk
bounded by) \(\Gamma_i\). Also, let \(\Gamma_0\subset\mathbb C\) be a circle
such that \(\mu_0,\mu_1,\ldots,\mu_\beta\) lie inside \(\Gamma_0\) and \(\mu_i\)
lie outside \(\Gamma_0\) for \(i>\beta\).

Since \(\mathbf T_{p,\lambda}\) depends continuously on \(p\in\mathcal P\),
there exists a neighbourhood \(\mathcal U\subset\mathcal U_0\) of \(p_0\) in
\(\mathcal P\) such that for all \(p\in\mathcal U,i=0,\ldots,\beta\), and
\(\lambda\in\Gamma_i\) we have \(\mathbf T_{p,\lambda}^{-1}\in\mathcal
L_G(\mathcal H)\). One can show then that
\[
\pi_p:=-\frac 1{2\pi\sqrt{-1}}\int_{\Gamma_i}\mathbf T_{p,\lambda}^{-1}\d\lambda
\]
is the projection onto the direct sum of eigenspaces corresponding to the
eigenvalues of \(\mathbf D_p\) lying inside \(\Gamma_i\). The right-hand side of
this equality depends continuously on \(p\in\mathcal U\) and therefore so does
the projection \(\pi_p\). A straightforward extension of~\cite[Chapter~I,
Lemma~4.10]{Kat1966} to the equivariant setting now shows that for each
\(p\in\mathcal U,i=1,\ldots,\beta\) there are exactly \(m(p_0,\mu_i)=:m_i\)
eigenvalues \(\mu_i^1\le\ldots\le\mu_i^{m_i}\) of \(\mathbf D_p\) inside the
circle \(\Gamma_i\) counting with \(G\)-multiplicity. Clearly,
\begin{equation}\label{eq:cont-multvecs}
\mathfrak m(p,\mu_i^j)\le\mathfrak m(p_0,\mu_i),\qquad i=1,\ldots,\beta,\quad
j=1,\ldots,m_i.
\end{equation}
Similarly, there are exactly \(m_1+\ldots+m_\beta\) eigenvalues of \(\mathbf
D_p\) inside the circle \(\Gamma_0\) counting with \(G\)-multiplicity. Since
\(\mu_0\) also lies inside \(\Gamma_0\), we obtain that
\[
\mu_i^j=\lambda_{m_1+\ldots+m_{i-1}+j}(p),\qquad
i=1,\ldots,\beta,\quad j=1,\ldots,m_i.
\]
In particular, taking \(i:=\beta\) and \(j:=\alpha-m_1-\ldots-m_{\beta-1}\), we
obtain that \(\mu_i^j=\lambda_\alpha(p)\) lies inside the circle
\(\Gamma_\beta\), and since this circle can be taken arbitrarily small, the
continuity of \(\lambda_\alpha(p)\) follows. Finally, it follows
from~\eqref{eq:cont-multvecs} that
\[
\mathfrak m(p,\lambda_\alpha(p))=\mathfrak m(p,\mu_i^j)\le
\mathfrak m(p_0,\mu_\beta)=\mathfrak m(p_0,\lambda_\alpha(p_0)),
\]
which shows the upper semicontinuity of \(\mathfrak m(p,\lambda_\alpha(p))\).
\end{proof}

\section{Generic equivariant simplicity of the spectrum for
Schr\texorpdfstring{\"o}{ö}dinger operators}\label{sec:schrod}

In this section we prove Theorem~\ref{th:outline-schrod} and therefore
Theorem~\ref{th:intro-sm} in particular.

\subsection{Petri condition}\label{sec:petri}

Let \(\mathcal P,\mathbf D\) be as in Definition~\ref{def:fam-subj} and
\(\mathbf T\) be given by~\eqref{eq:cont-T}. For \((p,\lambda)\in\mathcal
M=\mathcal P\times\mathbb R\) and \(K:=\ker\mathbf T_{p,\lambda}\), consider the
map
\[
\Lambda_{p,\lambda}\colon T_{p,\lambda}\mathcal M\to\Sym_G(K),\qquad
\Lambda_{p,\lambda}(q,\mu):=\pr_K\circ\,(\d_{p,\lambda}\mathbf T(q,\mu)|_K).
\]
Clearly, for any \(\epsilon\in\mathcal E\) the restriction
\(\Lambda_{p,\lambda;\epsilon}\) of \(\Lambda_{p,\lambda}\) onto \(\mathcal
M_\epsilon:=\mathcal P_\epsilon\times\mathbb R\) is a special case of the
map~\eqref{eq:walls-Lambda}. Therefore Lemma~\ref{lem:walls-main} suggests to
study the surjectivity of the map \(\Lambda_{p,\lambda;\epsilon}\). 

From now on we work under the assumptions of Theorem~\ref{th:outline-schrod},
i.e. \(\mathbf D\) is a family of \(G\)-equivariant Schr\"odinger operators that
comes from a linear, elliptic, formally self-adjoint \(G\)-equivariant
differential operator \(\mathbf L\) on \(M\). In this case we have
\[
\Lambda_{f,\lambda}(h,\mu)=\pr_K\circ\,(h-\mu)|_K,
\]
where \(f\in\mathcal P\subset C_G^\infty(M),\,h\in T_f\mathcal P\cong
C_G^\infty(M)\), and \(\lambda,\mu\in\mathbb R\). It turns out that the
surjectivity of the map \(\Lambda_{f,\lambda;\epsilon}\) can be restated in
terms of the \emph{injectivity} of the \textbf{Petri map} defined by
\[
\Pi\colon S^2 C^\infty(M)\to C^\infty(M),\quad
\Pi(u\odot v)(x):=u(x)v(x),\qquad x\in M.
\]

\begin{definition}
A linear elliptic formally self-adjoint \(G\)-equivariant differential
operator \(\mathbf P\colon C^\infty(M)\to C^\infty(M)\) satisfies the
\textbf{\(G\)-equivariant Petri condition} if the restriction of the map \(\Pi\)
onto \((S^2\ker\mathbf P)^G\) is injective.
\end{definition}

\begin{proposition}\label{pr:petri-main}
Fix \(f_\rf \in \mathcal P\). Suppose that \(\mathbf T_{f_\rf,\lambda}\)
satisfies the \(G\)-equivariant Petri condition. Then the map
\(\Lambda_{f_\rf,\lambda;\epsilon}\) is surjective for all \(\epsilon\) with
sufficiently rapid decay.
\end{proposition}

\begin{proof}
Observe that the \(G\)-invariant Riemannian metric \(g_\fix\) induces a perfect
pairing \(\llangle\cdot,\cdot\rrangle\) between \(\Sym(K)\) and \(S^2 K\)
defined on generators by
\[
\llangle A,u\odot v\rrangle:=\int_M v\,Au\,\d v_{g_\fix},\qquad
A\in\Sym(K),\quad u,v\in K.
\]
This pairing induces a perfect pairing between the \(G\)-invariants
\(\Sym_G(K)\) and \((S^2 K)^G\). We have
\begin{multline*}
\llangle\Lambda_{f_\rf,\lambda}(h,0),u\odot v\rrangle=
\int_M v\Lambda_{f_\rf,\lambda}(h,0)u\,\d v_{g_\fix}=
\int_M v\pr_K(hu)\d v_{g_\fix}=\\
=\int_M huv\,\d v_{g_\fix}=
\langle h,\Pi(u\odot v)\rangle_{L^2},\qquad
h\in C_G^\infty(M),\quad u,v\in K.
\end{multline*}
It follows that
\[
\llangle\Lambda_{f_\rf,\lambda;\epsilon}(h,0),B\rrangle=
\langle h,\Pi(B)\rangle_{L^2},\qquad
h\in C_G^\infty(M)_\epsilon,\quad B\in (S^2 K)^G.
\]
From parts~\ref{it:floer-chooseeps} and~\ref{it:floer-lb} of
Lemma~\ref{lem:floer-srd} we obtain that \(\Pi((S^2 K)^G)\subset
C_G^\infty(M)_\epsilon\) for all \(\epsilon\) with sufficiently rapid
decay. Choose such an \(\epsilon\) and let \(B\in (S^2 K)^G\) be orthogonal to
\(\im\Lambda_{f_\rf,\lambda;\epsilon}\) with respect to the pairing
\(\llangle\cdot,\cdot\rrangle\). Then for \(h=\Pi(B)\) we have
\[
0=\llangle\Lambda_{f_\rf,\lambda;\epsilon}(h,0),B\rrangle=\|\Pi(B)\|_{L^2}^2.
\]
Therefore \(\Pi(B)=0\), and since \(\mathbf T_{f_\rf,\lambda}\) satisfies the
\(G\)-equivariant Petri condition, we get \(B=0\). Since the pairing
\(\llangle\cdot,\cdot\rrangle\) is perfect, this concludes the proof.
\end{proof}

\subsection{The strong unique continuation property}\label{sec:ucon}

Unfortunately, there are easy examples where the Petri condition fails.

\begin{example}\label{ex:ucon-petrifail}
Let \(M=\mathbb R^2/\mathbb Z^2\) be a torus and
\(\Delta=-(\partial_x^2+\partial_y^2)\) be the flat Laplacian on \(M\). Then
\(f(x)=\sin 2\pi x,\,g(x)=\cos 2\pi x,\,f(y)=\sin 2\pi y,\,g(y)=\cos 2\pi y\)
are eigenfunctions of \(\Delta\) with eigenvalue \(4\pi^2\) and
\[
\Pi\big(f(x)f(y) \odot g(x)g(y)-
f(x)g(y) \odot g(x)f(y)\big)=0,
\]
the four different products being eigenfunctions with eigenvalue \(8\pi^2\).

\end{example}

The aim of this and the next sections is to show that Petri condition is
satisfied \emph{generically}. The proof consists of two ingredients: the strong
unique continuation property and the polynomial Wendl condition. In this section
we discuss the former ingredient. Let \(\mathbf P\) be a linear elliptic
formally self-adjoint differential operator on \(M\). 

\begin{definition}
We say that \(\mathbf P\) has the \textbf{strong unique continuation property}
at \(x\in M\) if the restriction of the map
\[
j_x^\infty\colon C^\infty(M)\to J_x^\infty(M),\qquad u\mapsto j_x^\infty u
\]
onto \(\ker\mathbf P\) is injective.
\end{definition}

Observe that for every \(l\in\mathbb N_0\cup\{\infty\}\) the Petri map induces
the map \[\Pi^l\colon S^2 J_x^l(M)\to J_x^l(M).\]

\begin{definition}
We say that \(\mathbf P\) satisfies the \textbf{\(\infty\)-jet Petri condition}
at \(x\in M\) if the restriction of the map \(\Pi^\infty\) onto \(S^2\ker\mathbf
P^\infty\) is injective.
\end{definition}

\begin{proposition}\label{pr:ucon-main}
Suppose that the operator \(\mathbf P\) is \(G\)-equivariant and has the strong
unique continuation property at \(x\in M\). If \(\mathbf P\) satisfies the
\(\infty\)-jet Petri condition at \(x\), then \(\mathbf P\) satisfies the
\(G\)-equivariant Petri condition.
\end{proposition}

\begin{proof}
The following diagram commutes:
\[
\begin{CD}
(S^2\ker\mathbf P)^G @>\Pi>> C_G^\infty(M)\\
@VVS^2 j_x^\infty V @VVj_x^\infty V\\
S^2\ker\mathbf P^\infty @>\Pi^\infty>> J_x^\infty(M).
\end{CD}
\]
The bottom map is injective by the assumption that \(\mathbf P\) satisfies the
\(\infty\)-jet Petri condition. Moreover, it is a standard fact of linear
algebra that the symmetric square of an injective linear map is injective.
Hence, the left map is injective by the assumption that \(\mathbf P\) has the
strong unique continuation property at \(x\in M\). The proposition now follows
immediately.
\end{proof}

\subsection{The Petri condition is generic}\label{sec:petrigen}

The aim of this section is to prove the following statement, assuming the
notions and results from Section~\ref{sec:wendl} regarding Wendl's conditions.

\begin{proposition}\label{pr:petrigen-main}
Let \(x\in M\) be regular and such that \(\sigma_x(\mathbf L)\) satisfies the
polynomial Wendl condition. There exists a Baire subset \(\mathcal
P_\Pi\subset\mathcal P\) such that for all \((f,\lambda)\in\mathcal M\) with
\(f\in\mathcal P_\Pi\) the operator \(\mathbf T_{f,\lambda}\) satisfies the
\(\infty\)-jet Petri condition at \(x\).
\end{proposition}

We follow Wendl's original approach (see~\cite{Wen2019} for an intuitive
explanation of the main ideas). The proof relies on the construction of the
\emph{universal Petri space}. It follows from~\eqref{eq:jets-dimkerP} that for
each \(l\in\mathbb N_0\),
\[
\mathcal K^l:=\{(f,\lambda,u)\in\mathcal M\times J_x^{l+k}(M)\colon
u\in\ker\mathbf T_{f,\lambda}^l\}
\]
is a vector bundle of rank
\begin{equation}\label{eq:petrigen-m}
m:=\binom{n+l+k}n-\binom{n+l}n.
\end{equation}
over \(\mathcal P\times\mathbb R\).

To proceed we need an extra restriction on a related space coming from the
following notion. Recall that \(\mathbf P\) induces a linear map \(\mathbf P^l\)
on jets, see~\eqref{eq:jets-Pl}.

By the \textbf{rank} \(\rk B\) of an element \(B\in S^2\ker\mathbf P^l\) we mean
its rank as an element of \((\ker\mathbf P^l)^{\otimes 2}\), i.e. the minimal
number \(\tau\in\mathbb N_0\) such that \(B=\sum_{i=1}^\tau u_i\otimes v_i\) for
some \(u_1,\ldots,u_\tau,v_1,\ldots,v_\tau\in\ker\mathbf P^l\). Equivalently, it
is the minimal number \(\tau\in\mathbb N_0\) such that \(B=\sum_{i=1}^\tau
\epsilon_i u_i\odot u_i\) for some \(u_1,\ldots,u_\tau\in\ker\mathbf P^l\) and
\(\epsilon_1,\ldots,\epsilon_\tau=\pm 1\). This coincides also with the rank of
the symmetric matrix represented by \(B\) w.r.t. any basis of \(\ker\mathbf
P^l\).

Now for all \(d\in\mathbb N_0,\tau\in\mathbb N\),
\begin{equation}\label{eq:petrigen-B}
\mathcal J_{d,\tau}^l:=\{(f,\lambda,B)\in S^2\mathcal K^l\colon\ord B\le d,\rk
B=\tau\}
\end{equation}
is a fiber bundle over \(\mathcal P\) with fibers of dimension
\begin{equation}\label{eq:petrigen-dimfib}
\tau m-\frac{\tau(\tau-1)}2+1.
\end{equation}
It comes with the Petri map
\[
\Pi_{d,\tau}^l\colon\mathcal J_{d,\tau}^l\to J_x^{l+k}(M),\qquad
\Pi_{d,\tau}^l(f,\lambda,B):=\Pi^{l+k}(B).
\]
The \textbf{universal Petri space} (of order at most \(d\) and rank \(\tau\)) is
\begin{equation}\label{eq:petrigen-petri}
\mathcal Z_{d,\tau}^l:=(\Pi_{d,\tau}^l)^{-1}(0).
\end{equation}
The following proposition uses the notion of polynomial Wendl condition, defined
in Section~\ref{sec:polywendl}. However, its use is purely formal at this stage. 

\begin{proposition}[the rank estimate]\label{pr:petrigen-rkest}
Suppose that \((f_\rf,\lambda)\in\mathcal M\) is such that \(\mathbf
T_{f_\rf,\lambda}\) satisfies the {polynomial Wendl condition at \(x\).} Then
for every \(d\in\mathbb N_0,\tau\in\mathbb N\) there exist constants
\(C(d,\tau),L(d,\tau)\) such that for all \(l\ge L(d,\tau)\) and
\((f_\rf,\lambda,B)\in\mathcal Z_{d,\tau}^l\) we have
\[
\rk\d_{(f_\rf,\lambda,B)}\Pi_{d,\tau}^l\ge C(d,\tau)l^n.
\]
\end{proposition}

\begin{proof}
For simplicity, we write just \(\mathbf T\) instead of \(\mathbf
T_{f_\rf,\lambda}\). Let \(l\in\mathbb N_0\). Since \(x\) is regular, by
Claim~\ref{cl:jets-regpoint} for any \(f\in J_x^l(M)\) there exists \(\tilde
f\in C_G^\infty(M)\) such that \(j_x^l\tilde f=f\). Then \(f(t):=f_\rf+t\tilde
f,t\in\mathbb R\) is a path in \(C_G^\infty(M)\). Set \(\mathbf T(t):=\mathbf
T_{f(t),\lambda}\) so that \(\mathbf T(0)=\mathbf T\), and let \(R\) be a right
inverse to \(\mathbf T^l\). For small \(t\), the operator \(\mathbf T(t)^l R\)
is close to the identity. In particular, it is invertible and \(R(t):=R(\mathbf
T(t)^l R)^{-1}\) is a right inverse to \(\mathbf T(t)^l\). For any
\(u\in\ker\mathbf T^l\) we define
\[
u(t):=u-R(t)\mathbf T(t)^l u.
\]
Then \(u(t)\in\ker\mathbf T(t)^l,u(0)=u\). Since \((\mathbf T^l)'(0)=\mathbf
T'(0)^l=f\), where \(f\) stands for the corresponding multiplication map
\(J_x^{l+k}(M)\to J_x^l(M)\) (note that we use~\eqref{eq:fam-clairaut} here), it
follows that \(-u'(0)=R((\mathbf T^l)'(0)u)=R(fu)\).

Now let \(d\in\mathbb N_0,\tau\in\mathbb N\) and \(B=\sum_{i=1}^\tau\epsilon_i
u_i\odot u_i\in S^2\ker\mathbf T^l\cap\ker\Pi^{l+k},\,\epsilon_i=\pm 1\) satisfy
\(\ord B\le d,\rk B=\tau\). Let the paths \(u_i(t)\) in \(J_x^{l+k}(M)\) be
constructed as in the previous paragraph. Put \(B(t):=\sum_{i=1}^\tau\epsilon_i
u_i(t)\odot u_i(t)\). Then for \(t\) small enough we still have \(\ord B(t)\le
d\) and \(\rk B(t)=\tau\), whence \(( f'(0),0,B'(0))\in
T_{(f,\lambda,B)}\mathcal J_{d,\tau}^l\). Moreover,
\begin{multline*}
-\frac12\d_{(f_\rf,\lambda,B)}\Pi_{d,\tau}^l( f'(0),0,B'(0))
=-\frac12\left.\frac{\d}{\d t}\right|_{t=0}\Pi^{l+k}(B(t))
=-\frac12\Pi^{l+k}(B'(0))=
\\=\Pi^{l+k}((Rf\odot\mathbf 1)B)
=\Pi^{l+k}((Rf\otimes\mathbf 1)B)
=W_B^l(f),
\end{multline*}
where \(W_B^l\) is the \emph{jet Wendl map} defined in
Section~\ref{sec:jetwendl}. Since \(\sigma_x(\mathbf L)\) satisfies the
polynomial Wendl condition, the proposition follows from
Claim~\ref{cl:polywendl-tojets}.
\end{proof}

Fix \(f_\rf\in\mathcal P\). For any \(d,l\in\mathbb N_0,\tau\in\mathbb N\), and
\(\epsilon\in\mathcal E\) consider the \(\epsilon\)-versions of the spaces
\eqref{eq:petrigen-B},\eqref{eq:petrigen-petri}:
\[
\mathcal J_{d,\tau;\epsilon}^l:=\{(f,\lambda,B)\in\mathcal J_{d,\tau}^l\colon
f\in\mathcal P_\epsilon\},\qquad
\mathcal Z_{d,\tau;\epsilon}^l:=\{(f,\lambda,B)\in\mathcal Z_{d,\tau}^l\colon
f\in\mathcal P_\epsilon\}.
\]
Moreover, let \(\Pi_{d,\tau;\epsilon}^l\) be the restriction of
\(\Pi_{d,\tau}^l\) to \(\mathcal J_{d,\tau;\epsilon}^l\). Also, for a subset
\(\mathcal V\subset S^2\mathcal K^l\) and \(f\in\mathcal P\), the intersection
of \(\mathcal V\) with the fiber of \(S^2\mathcal K^l\) over \(f\) is denoted by
\(\mathcal V(f)\). Recall that \(\mathcal K^l\) is a vector bundle over
\(\mathcal P\times\mathbb R\) and hence so is \(S^2\mathcal K^l\).
 
Recall the constants \(C(d,\tau), L(d,\tau)\) from
Proposition~\ref{pr:petrigen-rkest}.

\begin{definition}
Let \(\epsilon \in \mathcal E\). An element \((f,\lambda,B)\in\mathcal
Z_{d,\tau}^l\) is called \textbf{\(\epsilon\)-clean} if \(f\in\mathcal
P_\epsilon\) and there exists a neighbourhood \(\mathcal O\) of
\((f,\lambda,B)\) in \(\mathcal J_{d,\tau;\epsilon}^l\) such that \(\mathcal
O\cap\mathcal Z_{d,\tau;\epsilon}^l\) is a \(C^\infty\)-subvariety in \(\mathcal
J_{d,\tau;\epsilon}^l\) of codimension at least \(C(d,\tau)l^n\).
\end{definition}

\begin{proposition}\label{pr:petrigen-epsreg}
Let \(d\in\mathbb N_0,\tau\in\mathbb N\). For every \(\epsilon\) with
sufficiently rapid decay, all elements in \(\mathcal Z_{d,\tau}^l(f_\rf)\) for
\(l\ge L(d,\tau)\) are \(\epsilon\)-clean.
\end{proposition}

\begin{proof} 
Let \(z\in\mathcal Z_{d,\tau}^l(f_\rf)\) and let \(F\subset T_z\mathcal
J_{d,\tau}^l\) be a finite set such that \(\d_z\Pi_{d,\tau}^l(F)\) is a basis of
\(\im\d_z\Pi_{d,\tau}^l\). Then by Lemma~\ref{lem:floer-srd}\ref{it:floer-lb}
for every \(\epsilon\) with sufficiently rapid decay we have \(F\subset
T_z\mathcal J_{d,\tau;\epsilon}^l\). Then it follows from
Proposition~\ref{pr:petrigen-rkest} that
\[
\rk\d_z\Pi_{d,\tau;\epsilon}^l=\rk\d_z\Pi_{d,\tau}^l \ge C(d,\tau)l^n.\]
 In particular, there exists a neighbourhood
\(\mathcal O\) of \(z\) in \(\mathcal J_{d,\tau;\epsilon}^l\) such that
\[
\rk\d_{z'}\Pi_{d,\tau;\epsilon}^l\ge C(d,\tau)l^n,
\] for all \(z'\in\mathcal O\).
Then \(\mathcal O\cap\mathcal Z_{d,\tau;\epsilon}^l\) satisfies
Definition~\ref{def:cinfsub-subj} with \(W=J_x^{l+k}(M)\) and
\(f=\Pi_{d,\tau;\epsilon}^l\). This shows that the proposition is true with the
order of quantifiers exchanged. More precisely, for every \(z\in\mathcal
Z_{d,\tau}^l(f_\rf)\) there exists \(\epsilon_z\in\mathcal E\) such that \(z\)
is \(\epsilon_z\)-clean. Now observe that \(\epsilon_z\)-cleanness is an
open condition, i.e. \(z\) has a neighbourhood \(\mathcal U_z\) in \(\mathcal
Z_{d,\tau}^l(f_\rf)\) consisting of \(\epsilon_z\)-clean elements. As a
topological space, \(\mathcal Z_{d,\tau}^l(f_\rf)\) is just a subspace of the
total space of a vector bundle over \(\mathbb R\). In particular, it is
Lindel\"of, and the open cover \(\{\mathcal U_z\}_{z\in\mathcal
Z_{d,\tau}^l(f_\rf)}\) has an at most countable subcover. The proposition now
follows from Lemma~\ref{lem:floer-srd}\ref{it:floer-lb}.
\end{proof}

\begin{proposition}\label{pr:petrigen-subset}
For every \(\epsilon\in\mathcal E\) there exists a Baire subset \(\mathcal
P_{\Pi;\epsilon}\subset\mathcal P_\epsilon\) such that for any \(f\in\mathcal
P_{\Pi;\epsilon}\) and \(d\in\mathbb N_0,\tau\in\mathbb N\) the following holds:
there exists a constant \(L'(d,\tau)\) such that for all \(l\ge L'(d,\tau)\) the
set of \(\epsilon\)-clean elements in \(\mathcal Z_{d,\tau}^l(f)\) is empty.
\end{proposition}

\begin{proof}
Fix \(\epsilon\in\mathcal E\). Applying Lemma~\ref{lem:cinfsub-sardsmale} to the
canonical projection
\[
\mathcal J_{d,\tau;\epsilon}^l\to\mathcal P_\epsilon,\qquad
(f,\lambda,B)\mapsto f,
\]
we obtain a Baire subset \(\mathcal P_{\Pi;\epsilon}(d,\tau,l)\subset\mathcal
P_\epsilon\) such that for every \(f\in\mathcal P_{\Pi;\epsilon}(d,\tau,l)\) the
set of \(\epsilon\)-clean points is a \(C^\infty\)-subvariety in \(\mathcal
J_{d,\tau}^l(f)\) of codimension at least \(C(d,\tau)l^n\). Recall that
\(\mathcal J_{d,\tau}^l(f)\) is a smooth manifold of dimension given
by~\eqref{eq:petrigen-dimfib}. It follows that the set of \(\epsilon\)-clean
points in \(\mathcal Z_{d,\tau}^l(f)\) is empty as soon as
\begin{equation}\label{eq:petrigen-ineq}
\tau m-\frac{\tau(\tau-1)}2+1-C(d,\tau)l^n<0,
\end{equation}
where \(m\) is given by~\eqref{eq:petrigen-m}. Since \(m\) grows like
\(l^{n-1}\) in \(l\), there exists a constant \(L'(d,\tau)\) such that for all
\(l\ge L'(d,\tau)\) the inequality~\eqref{eq:petrigen-ineq} is satisfied. Hence
for all such \(l\) the set of \(\epsilon\)-clean points in \(\mathcal
Z_{d,\tau}^l(f)\) is empty for all \(f\in\mathcal P_{\Pi;\epsilon}(d,\tau,l)\).
It remains to put
\[
\mathcal P_{\Pi;\epsilon}:=\bigcap_{d\in\mathbb N_0,\tau\in\mathbb N}
\bigcap_{l\ge L'(d,\tau)}\mathcal P_{\Pi;\epsilon}(d,\tau,l).
\]
\end{proof}

Finally we prove Proposition~\ref{pr:petrigen-main} by passing from
\(\epsilon\)-spaces back to smooth functions.

\begin{proof}[Proof of Proposition~\ref{pr:petrigen-main}]
The proof is based on the so-called \emph{Taubes trick}. 

\textbf{Step 1: defining the sets \(\mathcal P_\reg^{K_1,K_2}\) and \(\mathcal
P_\reg\).} For each \(l\in\mathbb N_0\) choose a metric \(\gamma^l\) on \(S^2
J_x^l(M)\). Let \(K_1,K_2\in\mathbb N\). Fix \(l=l(K_1)\) so large that for all
\(\tau\in\mathbb N\) satisfying \(\tau\le K_1\) we have
\[
l=l(K_1)\ge\max\{K_1,L(K_1,\tau),L'(K_1,\tau)\}.
\]
Define the set
\[
\mathcal N^{K_1,K_2}\subset\bigcup_{\tau=1}^{K_1}\mathcal Z_{K_1,\tau}^l\subset
S^2\mathcal K^l
\]
in the following way: an element \((f,\lambda,B)\in\mathcal Z_{K_1,\tau}^l\)
belongs to \(\mathcal N^{K_1,K_2}\) whenever the following compactness
conditions are satisfied:
\begin{enumerate}[label=(\arabic*)]
\item{\emph{\(B\) is normalized:} \(\|B\|_{\gamma^{l+k}}=1\);}
\item{\emph{\(\lambda\) is not large:} \(|\lambda|\le K_2\);}
\item{\emph{the order of \(B\) does not blow up:} the distance from \(B\) to
the subspace \((S^2 J_x^{k+l}(M))^{K_1+1}\) in the metric \(\gamma^{l+k}\) is at
least \(1/K_2\);}
\item{\emph{the rank of \(B\) does not shrink suddenly:} the distance from \(B\)
to the subset of elements of rank at most \(\rk B-1\) in \(S^2 J_x^{l+k}(M)\) in
the metric \(\gamma^{l+k}\) is at least \(1/K_2\).}
\end{enumerate}
Finally, put
\[
\mathcal P_\reg^{K_1,K_2}:=\{f\in\mathcal P\colon\mathcal N^{K_1,K_2}(f)=
\emptyset\},\qquad
\mathcal P_\reg:=\bigcap_{K_1,K_2\in\mathbb N}
\mathcal P_\reg^{K_1,K_2}\subset\mathcal P.
\]

\textbf{Step 2: \(\mathcal P_\reg^{K_1,K_2}\) is open.} We prove that the
complement is closed. Let \(f_\nu\) be a sequence in \(\mathcal
P\setminus\mathcal P_\reg^{K_1,K_2}\) converging to some \(f\in\mathcal P\).
Then for each \(\nu\) there exists an element \(z_\nu\in\mathcal
N^{K_1,K_2}(f_\nu)\). The compactness conditions imposed on \(\mathcal
N^{K_1,K_2}\) provide a subsequence converging to an element \(z\in\mathcal
N^{K_1,K_2}(f)\). This shows that \(\mathcal N^{K_1,K_2}(f)\ne\emptyset\), that
is, \(f\in \mathcal P \setminus \mathcal P_\reg^{K_1,K_2}\).

\textbf{Step 3: \(\mathcal P_\reg^{K_1,K_2}\) is dense.} Fix
any \(f_\rf\in\mathcal P\). By Lemma~\ref{lem:floer-srd}\ref{it:floer-lb} there
exists \(\epsilon\in\mathcal E\) such that the condition of
Proposition~\ref{pr:petrigen-epsreg} is satisfied for \(d=K_1\) and any
\(\tau\le K_1\). Let \(\mathcal P_{\Pi;\epsilon}\subset\mathcal P_\epsilon\) be
the corresponding Baire subset given by Proposition~\ref{pr:petrigen-subset}.
Then there exists a sequence \(f_\nu\) in \(\mathcal P_{\Pi;\epsilon}\)
converging to \(f_\rf\). It suffices to show that there exists \(\nu\)
arbitrarily large such that \(f_\nu\in\mathcal P_\reg^{K_1,K_2}\). Assume the
converse. Again, by the compactness conditions imposed on \(\mathcal
N^{K_1,K_2}\), we may assume (after passing to a subsequence if necessary) that
there is some \(\tau\le K_1\) and a sequence \(z_\nu\in\mathcal
N^{K_1,K_2}(f_\nu)\cap\mathcal Z_{K_1,\tau}^l\) converging to an element
\(z\in\mathcal N^{K_1,K_2}(f_\rf)\cap\mathcal Z_{K_1,\tau}^l\), which is
\(\epsilon\)-clean by our choice of \(\epsilon\). Since
\(\epsilon\)-cleanness is an open condition, it follows that \(z_\nu\) is
\(\epsilon\)-clean for all \(\nu\) sufficiently large. With our choice of
\(l\), this contradicts the fact that \(f_\nu\in\mathcal P_{\Pi;\epsilon}\) due
to Proposition~\ref{pr:petrigen-subset}.

\textbf{Step 4: \(\infty\)-jet Petri condition is satisfied on \(\mathcal
P_\reg\).}
By Steps~2 and~3, \(\mathcal P_\reg\) is a Baire subset. Let
\((f,\lambda)\in\mathcal M\) with \(f\in\mathcal P_\reg\). Suppose that there is
a nontrivial element \(\widetilde B\in S^2\ker\mathbf T_{f,\lambda}^\infty\)
violating the \(\infty\)-jet Petri condition for \(\mathbf T_{f,\lambda}\). Let
\(d:=\ord\widetilde B,\tau:=\rk\widetilde B\). Choose \(K_1>d,\tau\) and let
\(B\) be the projection of \(\widetilde B\) onto \(S^2 J_x^{l+k}(M)\) normalized
so that \(\|B\|_{\gamma^{l+k}}=1\). Then \(B\) lies in \(\ker\Pi^{l+k}\) and
satisfies \(\ord B=d\) (in particular, \(B\) is nonzero) and \(\rk B\le\tau\).
Choosing \(K_2\) so large that \((f,\lambda,B)\in\mathcal
N^{K_1,K_2}\), we obtain a contradiction to the construction of \(\mathcal
P_\reg\).
\end{proof}

\subsection{Splitting the eigenvalues}\label{sec:split}

Let \(f_\rf\in\mathcal P\) be fixed.

\begin{definition}
A point \((f,\lambda)\in\mathcal M\) is called \textbf{\(\epsilon\)-regular} 
if \(f\in\mathcal P_\epsilon\) and the map \(\Lambda_{f,\lambda;\epsilon}\) is
surjective.
\end{definition}

For each \(\bar k\in\mathbb N_0^N\) and \(\epsilon\in\mathcal E\) we define
\[
\mathcal M^{\bar k}:=\{(f,\lambda)\in\mathcal M\colon [\ker\mathbf
T_{f,\lambda}]=\bar k\},\qquad
\mathcal M_\epsilon^{\bar k}:=\mathcal M^{\bar k}\cap\mathcal M_\epsilon.
\]
It follows from Lemma~\ref{lem:walls-main} that all \(\epsilon\)-regular
elements lying in \(\mathcal M_\epsilon^{\bar k}\) form a Banach submanifold in
\(\mathcal M_\epsilon\) of codimension given by~\eqref{eq:walls-codim}. Applying
the Sard-Smale theorem to the canonical projection \(\mathcal
M_\epsilon\to\mathcal P_\epsilon\), we obtain a Baire subset \(\mathcal
P_{\epsilon,\reg}^{\bar k}\subset\mathcal P_\epsilon\) such that for all
\(f\in\mathcal P_{\epsilon,\reg}^{\bar k}\) all \(\epsilon\)-regular elements
lying in \(\mathcal M_\epsilon^{\bar k}(f)\) form a Banach submanifold in
\(\mathcal M_\epsilon(f)=\mathbb R\) of the same codimension. In particular, if
this codimension is greater than one, then the fiber \(\mathcal M_\epsilon^{\bar
k}(f)\) does not contain any \(\epsilon\)-regular elements. This motivates the
following definition:
\begin{definition}
We call an element \(\bar k\in\mathbb N_0^N\) \textbf{simple} if \(|\bar{k}| =
1\) and \textbf{non-simple} if \(|\bar{k}|>1\).
\end{definition}
In other words, \(\bar k\) is simple if and only if the codimension
formula~\eqref{eq:walls-codim} gives one. Intersecting the subsets \(\mathcal
P_{\epsilon,\reg}^{\bar k}\) over all \(\bar k\) non-simple, we obtain the
following

\begin{proposition}\label{pr:split-subset}
For every \(\epsilon\in\mathcal E\) there exists a Baire subset \(\mathcal
P_{\epsilon,\reg}\subset\mathcal P_\epsilon\) such that for all \(f\in\mathcal
P_{\epsilon,\reg}\) and \(\bar k\in\mathbb N_0^N\) non-simple the fiber
\(\mathcal M_\epsilon^{\bar k}(f)\) does not contain any \(\epsilon\)-regular
elements.\qed
\end{proposition}

We are now ready to prove Theorems~\ref{th:outline-schrod}
and~\ref{th:mainres-schrodlap}.

\begin{proof}[Proof of Theorem~\ref{th:outline-schrod}]
For each \(K\in\mathbb N\) set \(I_K:=[-K,K]\) and using the notations of
Section~\ref{sec:cont} define
\[
\mathcal P_\reg^{K}:=\{f\in\mathcal P\colon
\mathfrak m(f,\lambda_\alpha(f))\text{ is simple }
\forall\alpha\in\mathbb N\text{ s.t. }
\lambda_\alpha(f)\in I_K\},\]\[
\mathcal P_\reg:=\bigcap_{K\in\mathbb N}\mathcal P_\reg^K.
\]
It suffices to show that for each \(K\in\mathbb N\) the set \(\mathcal
P_\reg^K\) is open and dense in \(\mathcal P\).

\textbf{\(\mathcal P_\reg^K\) is open.} We prove that the complement is
closed. Put
\[
\Xi^K(f):=\{i\in\mathbb N\colon\lambda_i(f)\in I_K\},\qquad
f\in\mathcal P.
\]
Let \(f_\nu\) be a sequence in \(\mathcal P\setminus\mathcal P_\reg^K\)
converging to some \(f\in\mathcal P\). It follows from
Proposition~\ref{pr:cont-main}\ref{it:cont-eigenval} that
\(\Xi^K(f_\nu)\subset\Xi^K(f)\) for all \(\nu\) sufficiently large. Hence, up to
a subsequence, we may assume that there exists \(\alpha\in\Xi^K(f)\) such that
for all \(\nu\), \(\lambda_\alpha(f_\nu)\in I_K\) and \(\mathfrak
m(f_\nu,\lambda_\alpha(f_\nu))=\bar k\) for some \(\bar k\in\mathbb N_0^N\) not
simple. Applying Proposition~\ref{pr:cont-main} again, we see that
\(\lambda_\alpha(f)=\lim\lambda_\alpha(f_\nu)\in I_K\) and \(\mathfrak
m(f,\lambda_\alpha(f))\ge\bar k\) is non-simple, i.e. \(f\in\mathcal P
\setminus \mathcal P_\reg^K\). Hence \(\mathcal P_\reg^K\) is open.

\textbf{\(\mathcal P_\reg^K\) is dense.} Fix any
\(f_\rf\in\mathcal P\). Perturbing \(f_\rf\) if necessary, we may
assume by Propositions~\ref{pr:ucon-main} and~\ref{pr:petrigen-main} that the
operator \(\mathbf T_{f_\rf,\lambda}\) satisfies the \(G\)-equivariant
Petri condition for all \(\lambda\in\mathbb R\). Then by
Proposition~\ref{pr:petri-main} and Lemma~\ref{lem:floer-srd}\ref{it:floer-lb}
there exists \(\epsilon\in\mathcal E\) such that the point
\((f_\rf,\lambda)\in\mathcal M\) is \(\epsilon\)-regular for all \(\lambda\in
I_K\). Let \(\mathcal P_{\epsilon,\reg}\subset\mathcal P_\epsilon\) be the Baire
subset given by Proposition~\ref{pr:split-subset}. Then there exists a sequence
\(f_\nu\) in \(\mathcal P_{\epsilon,\reg}\) converging to \(f_\rf\). It suffices
to show that there exists \(\nu\) arbitrarily large such that \(f_\nu\in\mathcal
P_\reg^K\). Assume the converse. Again, up to a subsequence we may assume that
there exists \(\alpha\in\Xi^K(f_\rf)\) such that for all \(\nu\),
\(\lambda_\alpha(f_\nu)\in I_K\) and \(\mathfrak
m(f_\nu,\lambda_\alpha(f_\nu))=\bar k\) for some \(\bar k\in\mathbb N_0^N\) not
simple. The point \((f_\rf,\lambda_\alpha(f_\rf))\) is \(\epsilon\)-regular by
our choice of \(\epsilon\). Since \(\epsilon\)-regularity is an open condition,
it follows that for all \(\nu\) sufficiently large,
\((f_\nu,\lambda_\alpha(f_\nu))\) is an \(\epsilon\)-regular element in the
fiber \(\mathcal M_\epsilon^{\bar k}\), which contradicts the construction of
\(\mathcal P_{\epsilon,\reg}\).
\end{proof}

\begin{proof}[Proof of Theorem~\ref{th:mainres-schrodlap}]
Let \(x\in M\) be regular. By~\cite{Aro1957},~\cite[Theorem~1.8]{Kaz1988}, the
operator \(\mathbf L+f-\lambda\) has the strong unique continuation property at
\(x\). By Proposition~\ref{pr:wendliff-lap}, the principal symbol of \(\mathbf
L\) at \(x\) satisfies the polynomial Wendl condition. The result now follows
from Theorem~\ref{th:outline-schrod}.
\end{proof}

\section{Generic equivariant simplicity for conformally
covariant operators}\label{sec:confcov}

In this section we work under the assumptions of
Theorem~\ref{th:outline-confcov}. The proof goes along the general lines of the
proof of Theorem~\ref{th:outline-schrod}. For this reason, here we only indicate
two points where the proof of Theorem~\ref{th:outline-schrod} requires a
modification. All other changes are cosmetic and are left to the reader.

\subsection{Petri condition}

Fix a metric \(g_\fix\in\mathcal P\). As explained in
Example~\ref{ex:fam-metdep}, \(\mathbf L\) gives rise to a smooth family
\(\mathbf D\) of \(G\)-equivariant operators on \(M\) parametrized by \(\mathcal
P\). Let \(\mathbf T\) be as in~\eqref{eq:cont-T} and let the maps
\(\Lambda_{g,\lambda},\Lambda_{g,\lambda;\epsilon}\) be as in
Section~\ref{sec:petri}. Then for all \(g\in\mathcal P\) and \(h\in T_g\mathcal
P\) we have
\begin{align*}
\d_{g,\lambda}\mathbf T(h,0)&=
\d_g(\rho^{\frac12}\circ\mathbf L\circ\rho^{-\frac12})(h)\\
&=\begin{multlined}[t]
\frac12\rho'(h)\rho(g)^{-1}\circ\mathbf D_g-
\frac12\mathbf D_g\circ\rho'(h)\rho(g)^{-1}\,+\\
\rho(g)^{\frac12}\circ\d_g\mathbf L(h)\circ\rho(g)^{-\frac12}
\end{multlined}\\
&=\begin{multlined}[t]\frac12\rho'(h)\rho(g)^{-1}\circ\mathbf T_{g,\lambda}-
\frac12\mathbf T_{g,\lambda}\circ\rho'(h)\rho(g)^{-1}\,+\\
\rho(g)^{\frac12}\circ\d_g\mathbf L(h)\circ\rho(g)^{-\frac12},
\end{multlined}\\
\Lambda_{g,\lambda}(h,0)&=
\pr_K\circ(\rho(g)^{\frac12}\circ\d_g\mathbf
L(h)\circ\rho(g)^{-\frac12})|_K.
\end{align*}
Now, using~\eqref{eq:fam-conf}, for any \(h\in T_g\mathcal P\) of the form
\(h=2fg\in C_G^\infty(S^2 T^*M)\), where \(f\in C_G^\infty(M)\), we have
\begin{align}
\d_g\mathbf L(h)
&=\left.\frac{\d}{\d t}\right|_{t=0}\mathbf L_{e^{2ft}g}
=\left.\frac{\d}{\d t}\right|_{t=0}(e^{-bft}\circ\mathbf L_g\circ
e^{aft})\label{eq:petri2-varT}\\
&=-bf\circ\mathbf L_g+\mathbf L_g\circ af,\notag
\\
\rho(g)^{\frac12}\circ\d_g\mathbf L(h)\circ\rho(g)^{-\frac12}&=
-bf\circ\mathbf D_g+\mathbf D_g\circ af\notag\\
&=-bf\circ\mathbf T_{g,\lambda}+\mathbf T_{g,\lambda}\circ af+
\lambda(a-b)f.\notag
\end{align}
Hence,
\begin{equation}\label{eq:petri2-Lambda}
\Lambda_{g,\lambda}(h,0)=\pr_K\circ(\rho(g)^{\frac12}\circ\d_g\mathbf
L(h)\circ\rho(g)^{-\frac12})|_K=\pr_K\circ\,(\lambda(a-b)f|_K).
\end{equation}

With~\eqref{eq:petri2-Lambda} in hand, the following proposition is proved
in the same way as Proposition~\ref{pr:petri-main}.

\begin{proposition}\label{pr:petri2-main}
Suppose that \(\mathbf T_{g_\rf,\lambda}\) with \(\lambda\ne 0\) satisfies the
\(G\)-equivariant Petri condition. Then the map
\(\Lambda_{g_\rf,\lambda;\epsilon}\) is surjective for all \(\epsilon\) with
sufficiently rapid decay.
\end{proposition}

\subsection{The rank estimate}\label{sec:rkest2}

Another point where the proof of Theorem~\ref{th:outline-schrod} has to be
modified is the proof of Proposition~\ref{pr:petrigen-main}. Put
\(\widetilde{\mathbf T}_{g,\lambda}:=\mathbf L_g-\lambda\). Observe that the map
\begin{equation}\label{eq:rkest2-tilde}
\ker\mathbf T_{g,\lambda}\to
\ker\widetilde{\mathbf T}_{g,\lambda},\quad
\tilde u:=\rho(g)^{-\frac12}u,
\end{equation}
is an isomorphism. Moreover, the induced map on \(S^2\ker\mathbf T_{g,\lambda}\)
preserves \(\ker\Pi\). It follows that we can replace \(\mathbf
T_{g,\lambda}\) by \(\widetilde{\mathbf T}_{g,\lambda}\) while proving that the
\(G\)-equivariant Petri condition is generic.

\begin{proposition}[the rank estimate]\label{pr:rkest2-rkest}
Suppose that \((g,\lambda)\in\mathcal M^*:=\mathcal P\times\mathbb R^*\) is such
that \(\widetilde{\mathbf T}_{g,\lambda}\) satisfies the jet Wendl condition at
\(x\). For every \(d\in\mathbb N_0,\tau\in\mathbb N\) there exist constants
\(C(d,\tau),L(d,\tau)\) such that for all \(l\ge L(d,\tau)\) the following
holds: for all \((g,\lambda,B)\in\mathcal Z_{d,\tau}^l\) we have
\[
\rk\d_{(g,\lambda,B)}\Pi_{d,\tau}^l\ge C(d,\tau)l^n.
\]
\end{proposition}

To prove this, we need the following auxiliary claim (see also
Remark~\ref{rem:jetwendl-R}).

\begin{claim}\label{cl:rkest2-addK}
Suppose that a linear elliptic formally self-adjoint operator \(\mathbf P\)
satisfies the jet Wendl condition at \(x\) and let
\[
K^l\colon J_x^l(M)\times\ker\mathbf P^l\to\ker\mathbf P^l
\]
be a family of bilinear maps depending on \(l\in\mathbb N_0\). Put
\[
\widetilde W_B^l\colon J_x^l(M)\to J_x^{l+k}(M),\qquad
\widetilde W_B^l(f):=W_B^l(f)+
\Pi^{l+k}((K^l(f,\cdot)\otimes\mathbf 1)B).
\]
Then, up to a change of the constants \(C\) and \(L\), the jet Wendl
condition is satisfied with \(\widetilde W_B^l\) instead of \(W_B^l\).
More precisely, for all \(d\in\mathbb N_0,\tau\in\mathbb N\) there exist
constants \(\widetilde C=\widetilde C(d,\tau),\widetilde L=\widetilde
L(d,\tau)\) such that for all \(l\ge\widetilde L\) the following holds: for
every \(B\in\ker\Pi^{l+k}\cap S^2\ker\mathbf P^l\) satisfying \(\ord B\le d,\rk
B=\tau\) the map \(\widetilde W_B^l\) satisfies
\[
\rk\widetilde W_B^l\ge\widetilde Cl^n.
\]
\end{claim}

\begin{proof}
Let \(l\in\mathbb N_0\). Take any \(B\in\ker\Pi^{l+k}\cap S^2\ker\mathbf P^l\)
satisfying \(\ord B\le d\) and \(\rk B=\tau\) and write
\[
B=\sum_{i=1}^\tau u_i\odot v_i,\qquad u_i,v_i\in\ker\mathbf P^l.
\]
For each \(i=1,\ldots,\tau\) consider the map \(w_i\colon J_x^l(M)\to
J_x^{l+k}(M)\),
\[
w_i(f):=\Pi^{l+k}((K^l(f,\cdot)\otimes\mathbf 1)(u_i\odot v_i))
=\frac12(K^l(f,u_i)v_i+u_i K^l(f,v_i)).
\]
It follows from~\eqref{eq:jets-dimkerP} that there exist constants
\(C'=C'(n,k),L'=L'(n,k)\) such that
\[
\dim\ker\mathbf P^l\le C'l^{n-1},\qquad l\ge L'.
\]
In particular, for each \(i=1,\ldots,\tau\) and \(l\ge L'\) we have
\[
\rk w_i\le 2\dim\ker\mathbf P^l\le 2C'l^{n-1}.
\]
Thus we get the estimate
\[
\rk\widetilde W_B^l
\ge\rk W_B^l-\rk (W_B^l-\widetilde W_B^l)
\ge\rk W_B^l-\sum_{i=1}^\tau\rk w_i
\ge Cl^n-2\tau C'l^{n-1},\quad l\ge L'.
\]
It follows that for some new constants \(\widetilde C,\widetilde L\) we have
\[
\rk\widetilde W_B^l\ge Cl^n-2\tau C'l^{n-1}\ge\widetilde Cl^n,\qquad
l\ge\widetilde L,
\]
as desired.
\end{proof}

\begin{proof}[Proof of Proposition~\ref{pr:rkest2-rkest}]
As in the proof of Proposition~\ref{pr:petrigen-rkest}, we write just \(\mathbf
L\) and \(\widetilde{\mathbf T}\) instead of \(\mathbf L_g\) and
\(\widetilde{\mathbf T}_{g,\lambda}\) for simplicity. Let \(l\in\mathbb N_0\)
and
\[
{}^*\colon J_x^l(M)\to J_x^{l+k}(M),\qquad f\mapsto f^*
\]
be an arbitrary right inverse to the natural projection \(J_x^{l+k}(M)\to
J_x^l(M)\). Since \(x\) is regular, by Claim~\ref{cl:jets-regpoint} for
every \(f\in J_x^l(M)\) there exists \(\tilde f\in C_G^\infty(M)\) such that
\(j_x^{l+k}\tilde f=f^*\) (and hence \(j_x^l\tilde f=f\)). Then
\(g(t):=e^{2\tilde ft}g\) is a path of metrics in \(\mathcal C^G\). Let
\(u\in\ker\widetilde{\mathbf T}^l\) and \(\widetilde{\mathbf T}(t),R(t),u(t)\)
be defined as in the proof of Proposition~\ref{pr:petrigen-rkest}.
Using~\eqref{eq:petri2-varT} it is easy to see that
\[
(\widetilde{\mathbf T}^l)'(0)=\widetilde{\mathbf T}'(0)^l
=(-b\tilde f\circ\mathbf L+\mathbf L\circ a\tilde f)^l
=-bf\circ\mathbf L^l+\mathbf L^l\circ af^*.
\]
Hence
\begin{equation}\label{eq:rkest2-lift}
-u'(0)=R(\widetilde{\mathbf T}^l)'(0)u=R(-bf\mathbf L^l u+\mathbf L^l(af^* u))
=\lambda(a-b)R(fu)+R\widetilde{\mathbf T}^l(af^* u).
\end{equation}
Now let \(d\in\mathbb N_0,\tau\in\mathbb N\) and \(B\in
S^2\ker\widetilde{\mathbf T}^l\cap\ker\Pi^{l+k}\) satisfy \(\ord B\le d,\rk
B=\tau\). Let \(B(t)\) be defined as in the proof of
Proposition~\ref{pr:petrigen-rkest}. It follows from~\eqref{eq:rkest2-lift} that
\begin{align*}
-\frac12\d_{(g,\lambda,B)}\Pi_{d,\tau}^l(g'(0),0,B'(0))
&=-\frac12\left.\frac{\d}{\d t}\right|_{t=0}\Pi^{l+k}(B(t))
=-\frac12\Pi^{l+k}(B'(0))\\
&=\Pi^{l+k}(((\lambda(a-b)Rf+R\widetilde{\mathbf T}^l af^*)\odot\mathbf 1)B)\\
&=\Pi^{l+k}(((\lambda(a-b)Rf+R\widetilde{\mathbf T}^l af^*)\otimes\mathbf 1)B)\\
&=\Pi^{l+k}(((\lambda(a-b)Rf+(R\widetilde{\mathbf T}^l-\mathbf 1)af^*)\otimes
\mathbf 1)B)\\
&=\lambda(a-b)W_B^l(f)+\Pi^{l+k}(((R\widetilde{\mathbf T}^l-\mathbf
1)af^*\otimes\mathbf 1)B),
\end{align*}
where we used that
\[
\Pi^{l+k}((af^*\otimes\mathbf 1)B)=af^*\Pi^{l+k}(B)=0.
\]
The proposition now follows from Claim~\ref{cl:rkest2-addK} with
\[
K^l(f,\cdot):=\frac 1{\lambda(a-b)}(R\widetilde{\mathbf T}^l-\mathbf 1)af^*.
\]
\end{proof}

\begin{proof}[{Proof of Theorems~\ref{th:outline-confcov}
and~\ref{th:mainres-confcov}}]
The proofs are similar to the proofs of Theorems~\ref{th:outline-schrod}
and~\ref{th:mainres-schrodlap}. The set \(I_K=[-K,K]\) from the proof of
Theorem~\ref{th:outline-schrod} should be replaced by \([-K,-1/K]\cup [1/K,K]\)
in order to avoid the zero eigenvalue.
\end{proof}

\section{Generic equivariant simplicity for the Laplace-Beltrami
operator}\label{sec:laplace}

\subsection{\texorpdfstring{The traceless \(\Delta\)}{Δ}-Petri condition}

We assume \(n=\dim M>2\) throughout this section. Given a smooth \(G\)-action on
\(M\), consider the Laplace-Beltrami operator \(\Delta\) as a \(G\)-equivariant
metric dependent differential operator of order two on \(M\). As explained in
Example~\ref{ex:fam-metdep}, \(\Delta\) gives rise to a smooth family \(\mathbf
D\) of \(G\)-equivariant operators on \(M\) parametrized by \(\mathcal P\). Let
\(\mathbf T_{g,\lambda}\) be as in~\eqref{eq:cont-T}, \(\widetilde{\mathbf
T}_{g,\lambda}:=\Delta_g-\lambda\), and let the maps
\(\Lambda_{g,\lambda},\Lambda_{g,\lambda;\epsilon}\) be as in
Section~\ref{sec:petri}. For all \(h\in T_g\mathcal P\cong C_G^\infty(S^2
T^*M)\) we have
\[
\d_g\Delta(h)\tilde u=
\div_g (d\tilde u\iprod h)-\frac12\langle\d(\tr_g h),\d\tilde u\rangle_g,\qquad
\tilde u\in C^\infty(M).
\]
Applying the divergence theorem, it is easy to see that for all \(\tilde
u,\tilde v\in\ker\widetilde{\mathbf T}_{g,\lambda},
\)
\[
\int_M\tilde v\,\d_g\Delta(h)\tilde u\,\d v_{g}=
\int_M\langle h,-\Pi_\d(\d\tilde u\odot\d\tilde v)+
\frac12\tr_g(\Pi_\d(\d\tilde u\odot\d\tilde v))g-
\frac{\lambda}2\Pi(\tilde u\odot\tilde v)g\rangle_g\d v_g,\]

where
\[
\Pi_\d\colon S^2 C^\infty(T^*M)\to C^\infty(S^2 T^*M),\qquad
\Pi_\d(\theta\odot\omega)(x):=\theta(x)\odot\omega(x).
\]
Note that the first equality in~\eqref{eq:petri2-Lambda} holds for any metric
dependent operator, not necessarily conformally covariant. Since in our case
\(\mathbf L_g=\Delta_g\), we have
\[
\Lambda_{g,\lambda}(h,0)=
\pr_K\circ\,(\rho(g)^{\frac12}\circ\d_g\Delta(h)\circ
\rho(g)^{-\frac12})|_K,\qquad K:=\ker\mathbf T_{g,\lambda},
\]
and finally
\begin{align*}
\int_M v\Lambda_{g,\lambda}(h,0)u\,\d v_{g_\fix}
&=\int_M \rho(g)^{\frac12}v\,\d_g\Delta(h)(\rho(g)^{-\frac12}u)\d v_{g_\fix}
\\&=\int_M \tilde v\,\d_g\Delta(h)\tilde u\,\d v_g
\\&=-\langle h,\Pi_{g,\lambda}^\Delta(\tilde u\odot\tilde v)
\rangle_{L^2(S^2 T^*M)},\qquad
u,v\in K,
\end{align*}
where \(\tilde u:=\rho(g)^{-\frac12}u,\tilde v:=\rho(g)^{-\frac12}v\) are as
in~\eqref{eq:rkest2-tilde}, the inner product
\(\langle\cdot,\cdot\rangle_{L^2(S^2 T^*M)}\) is induced by the metric \(g\),
and
\[
\Pi_{g,\lambda}^\Delta\colon S^2 C^\infty(M)\to C^\infty(S^2 T^*M),\]\[
\Pi_{g,\lambda}^\Delta(B):=\Pi_\d(B_\d)-\frac12\tr_g(\Pi_\d(B_\d))g+
\frac{\lambda}2\Pi(B)g,\quad B_\d:=(\d\otimes\d)B
\]
is the \textbf{\(\Delta\)-Petri map}. In what follows we use only the traceless
part of \(\Pi_{g,\lambda}^\Delta\), which is given by
\[
\Pi_g^\perp\colon S^2 C^\infty(T^*M)\to C^\infty(S^2 T^*M),\qquad
\Pi_g^\perp(B_\d):=\Pi_\d(B_\d)-\frac{1}n\tr_g(\Pi_\d(B_\d))g.
\]
This comes with the following definition.

\begin{definition}
We say that the \textbf{\(G\)-equivariant traceless \(\Delta\)-Petri condition}
is satisfied at \((g,\lambda)\in\mathcal M^*\) if the restriction of the
traceless \(\Delta\)-Petri map \(\Pi_g^\perp\) onto
\((S^2(\d\ker\widetilde{\mathbf T}_{g,\lambda}))^G\) is injective.
\end{definition}

The following proposition is proved similarly to
Propositions~\ref{pr:petri-main} and~\ref{pr:petri2-main}.

\begin{proposition}
Suppose that the \(G\)-equivariant traceless \(\Delta\)-Petri condition is
satisfied at \((g,\lambda)\in\mathcal M^*\). Then the map
\(\Lambda_{g,\lambda;\epsilon}\) is surjective for all \(\epsilon\) with
sufficiently rapid decay.\qed
\end{proposition}

\subsection{The rank estimate}\label{sec:rkest3}

Let \(x\in M\) be regular. For each \(l\in\mathbb N_0\cup\{\infty\}\) and
\((g,\lambda)\in\mathcal M^*\) the map \(\Pi_g^\perp\) induces the map
\[
\Pi_g^{\perp;l}\colon S^2 J_x^l(T^*M)\to J_x^l(S^2 T^*M),\qquad
\Pi_g^{\perp;l}(B_\d)=\Pi_\d^l(B_\d)-\frac 1n\tr_g^l(\Pi_\d^l(B_\d))j_x^l g.
\]
where the map \(\Pi_\d^l\colon S^2 J_x^l(T^*M)\to J_x^l(S^2 T^*M)\) is induced
by \(\Pi_\d\).

We say that the \textbf{\(\infty\)-jet traceless \(\Delta\)-Petri condition} is
satisfied at \((g,\lambda)\in\mathcal M^*\) if the restriction of the map
\(\Pi_g^{\perp;\infty}\) onto \(S^2(\d\ker\widetilde{\mathbf
T}_{g,\lambda}^\infty)\) is injective.

\begin{proposition}
There exists a Baire subset \(\mathcal P_{\Pi^\perp}\subset\mathcal P\) such
that for all \((g,\lambda)\in\mathcal M^*\) with \(g\in\mathcal P_{\Pi^\perp}\)
the \(\infty\)-jet traceless \(\Delta\)-Petri condition is satisfied at
\((g,\lambda)\).
\end{proposition}

For each \(d\in\mathbb N_0,\tau\in\mathbb N\) we put
\begin{align*}
\mathcal K^l&:=\{(g,\lambda,u)\in\mathcal M^*\times J_x^{l+k}(M)\colon
u\in\ker\widetilde{\mathbf T}_{g,\lambda}^l\},\\
\mathcal J_{d,\tau}^l&:=\{(g,\lambda,B)\in S^2\mathcal K^l\colon
\ord B_\d\le d,\rk B_\d=\tau\},
\end{align*}
where now \(B_\d:=(\d^{l+1}\otimes\d^{l+1})B\). Note that \(\d^{l+1}\) is
injective on \(\ker\widetilde{\mathbf T}_{g,\lambda}^l\) since the latter does
not contain constants. In particular, \(\mathcal J_{d,\tau}^l\) is a fiber
bundle as in Section~\ref{sec:petrigen}. Put
\[
\Pi_{d,\tau}^{\perp;l}\colon\mathcal J_{d,\tau}^l\to J_x^{l+1}(S^2 T^*M),
\qquad
\Pi_{d,\tau}^{\perp;l}(g,\lambda,B):=\Pi_g^{\perp;l+1}(B_\d).
\]
The \textbf{universal traceless \(\Delta\)-Petri space} is given by \(\mathcal
Z_{d,\tau}^{\perp;l}:=(\Pi_{d,\tau}^{\perp;l})^{-1}(0)\).

\begin{proposition}[the rank estimate]
For every \(d\in\mathbb N_0,\tau\in\mathbb N\) there exist constants
\(C(d,\tau),L(d,\tau)\) such that for all \(l\ge L(d,\tau)\) the following
holds: for all \((g,\lambda,B)\in\mathcal Z_{d,\tau}^{\perp;l}\) we have
\[
\rk\d_{(g,\lambda,B)}\Pi_{d,\tau}^{\perp;l}\ge C(d,\tau)l^n.
\]
\end{proposition}

\begin{proof}
As in the proof of Propositions~\ref{pr:petrigen-rkest}
and~\ref{pr:rkest2-rkest}, we write just \(\widetilde{\mathbf T}\) and
\(\Delta\) instead of \(\widetilde{\mathbf T}_{g,\lambda}\) and \(\Delta_g\) for
simplicity. Also, since in this proof we work mostly with jets, we often
identify the maps of sections with the induced maps on jets, e.g. we write
simply \(\Pi_\d\) instead of \(\Pi_\d^l\), \(\tr_g\) instead of \(\tr_g^l\). 

Let \(l\in\mathbb N_0\) and \(h\in J_x^{l+1}(S^2 T^*M)\) satisfy \(\tr_g
h=0\). Since \(x\) is regular, by Claim~\ref{cl:jets-regpoint} there exists
\(\tilde h\in C_G^\infty(S^2 T^*M)\) such that \(j_x^{l+1}\tilde h=h\). Then,
for some small \(\delta>0\), \(g(t)=g+t\tilde h,t\in (-\delta,\delta)\) is a
path of metrics in \(\mathcal P\). Let \(u\in\ker\widetilde{\mathbf T}\) and
\(\mathbf T(t),R(t),u(t)\) be defined as in the proof of
Proposition~\ref{pr:petrigen-rkest}. We have
\[
(\widetilde{\mathbf T}^l)'(0)u=\widetilde{\mathbf T}'(0)^l u=
\div_g(\d u\iprod h).
\]
Hence,
\begin{equation}\label{eq:rkest3-lift}
-\d u'(0)=\d R(\widetilde{\mathbf T}^l)'(0)u=\d R\div_g(\d u\iprod h).
\end{equation}

Now let \(d\in\mathbb N_0,\tau\in\mathbb N\), and \((g,\lambda,B)\in\mathcal
Z_{d,\tau}^{\perp;l}\). Note that since \(\Pi_g^{\perp;l+1}(B_\d)=0\), we have
\[
\langle\Pi_\d(B_\d),h\rangle_g=
\frac 1n\langle\tr_g(\Pi_\d(B_\d))j_x^{l+1}g,h\rangle_g=
\frac 1n\tr_g(\Pi_\d(B_\d))\tr_g h=0.
\]
Let \(B=\sum_{i=1}^\tau\epsilon_i u_i\odot u_i,\,\epsilon_i=\pm1\), the paths
\(u_i(t)\) in \(J_x^{l+2}(M)\) be as above, and
\(B(t):=\sum_{i=1}^\tau\epsilon_i u_i(t)\odot u_i(t)\).
Using~\eqref{eq:rkest3-lift}, we obtain
\begin{align*}
-\d_{(g,\lambda,B)}\Pi_{d,\tau}^{\perp;l}(g'(0),0,B'(0))&=
-\left.\frac{\d}{\d t}\right|_{t=0}\Pi_{g(t)}^{\perp;l+1}(B_\d(t))\\&=
\begin{multlined}[t]
-\Pi_g^{\perp;l+1}(B_\d'(0))+
\frac{1}n\bigl(\tr_g(\Pi_\d(B_\d))h-\\
\langle\Pi_\d(B_\d),h\rangle_g\,j_x^{l+1}g\bigr)
\end{multlined}\\&=
-\Pi_g^{\perp;l+1}(B_\d'(0))+\frac{1}n\tr_g(\Pi_\d(B_\d))h=W_B^{\perp;l}(h),
\end{align*}
where \(W_B^{\perp;l}\) is the \emph{jet traceless \(\Delta\)-Wendl map} defined
in Section~\ref{sec:tdwendl}. The proposition now follows from
Claim~\ref{cl:tdwendl-tojets} and Proposition~\ref{pr:tdwendl-main}.
\end{proof}

\begin{proof}[Proof of Theorem~\ref{th:mainres-laplace}]
For \(\dim M=2\) the result follows from Theorem~\ref{th:mainres-confcov} so we
may assume that \(\dim M>2\). For any \(g\in\mathcal G\) there is a single
eigenspace corresponding to the eigenvalue \(\lambda=0\). This eigenspace
is formed by constant eigenfunctions and is a trivial \(G\)-representation,
hence we may assume that \(\lambda\ne 0\). The rest of the argument is similar
to the proof of Theorems~\ref{th:outline-schrod} and~\ref{th:mainres-schrodlap}
with Petri condition replaced by the traceless \(\Delta\)-Petri condition.
\end{proof}

\section{Wendl conditions}\label{sec:wendl}

\subsection{The jet Wendl condition}\label{sec:jetwendl}

Let an elliptic formally self-adjoint operator \(\mathbf P\) on \(M\) and a
point \(x\in M\) be fixed. Let \(l\in\mathbb N_0\). Recall that any \(f\in
J_x^l(M)\) induces a multiplication map \(J_x^{l+k}(M)\to J_x^l(M)\) also
denoted by \(f\). For \(B\in S^2 J_x^{l+k}(M)\) and a right inverse \(R\) to
\(\mathbf P^l\) define the \textbf{jet Wendl map} by
\[
W_B^l\colon J_x^l(M)\to J_x^{l+k}(M),\qquad
W_B^l(f):=\Pi^{l+k}((Rf\otimes\mathbf 1)B).
\]

\begin{definition}\label{def:jetwendl-subj}
We say that \(\mathbf P\) satisfies the \textbf{jet Wendl condition} at \(x\) if
for any choice of the map \(R\) and every \(d\in\mathbb N_0,\tau\in\mathbb N\)
there exist constants \(C=C(d,\tau),L=L(d,\tau)\) such that for all \(l\ge L\)
the following holds: for all \(B\in\ker\Pi^{l+k}\cap S^2\ker\mathbf P^l\)
satisfying \(\ord B\le d,\rk B=\tau\) the map \(W_B^l\) satisfies
\[
\rk W_B^l\ge Cl^n.
\]
\end{definition}

\begin{remark}\label{rem:jetwendl-R}
The choice of the map \(R\) does not really matter, i.e., if
Definition~\ref{def:jetwendl-subj} is satisfied for \emph{some} choice of the
map \(R\), then, up to a change of the constants \(C\) and \(L\), it is
satisfied for any other choice of the map \(R\). This is because for any other
choice \(\widetilde R\) we have \(\im(R-\widetilde R)\subset\ker\mathbf P^l\),
while the dimension of the latter space grows like \(l^{n-1}\)
by~\eqref{eq:jets-dimkerP} (see also the proof of Claim~\ref{cl:rkest2-addK}).
The same remark applies to the choice of \(\widehat R\) in the definition of the
polynomial Wendl map \(\widehat W_{\widehat B}\) below.
\end{remark}

\subsection{The polynomial Wendl condition}\label{sec:polywendl}

The jet Wendl condition is what we really need. The problem is that it is
difficult to show directly that it is satisfied. In fact, usually there is not
even a satisfying description of \(\ker\mathbf P^l\). However, if one looks
only at the lowest order terms, then the operator reduces to its principal
symbol, which makes the rank estimation problem much more tractable. In this
form the idea appears in~\cite[proof of Theorem~1.6.17]{DW2023}, cf. Wendl's
\emph{rescaling trick}~\cite[Section~5.2.4]{Wen2023}.

Let \(x\in M\) and let \(\sigma\in S^k T_x M\) be an elliptic self-adjoint
symbol of order \(k\). The \textbf{polynomial Petri map} is defined by
\[
\widehat\Pi\colon S^2(S^\bullet T_x^*M)\to S^\bullet T_x^*M,\qquad
\widehat\Pi(u\odot v):=uv.
\]
Take an element \(\widehat B\in S^2(S^\bullet T_x^*M)\). We say that \(\widehat
B\) is \textbf{homogeneous} of degree \(\hat d\) (and write \(\deg\widehat
B=\hat d\)) if it can be written in the form \(\widehat B=\sum_{i=1}^\tau
u_i\odot v_i\) for some homogeneous \(u_i,v_i\in S^\bullet T_x^*M\) with
\(\deg u_i+\deg v_i=\hat d\) for each \(i=1,\ldots,\tau\). Also, note that for
every \(l\in\mathbb N_0\) an element \(\hat f\in S^l T_x^*M\) induces a map
\(S^\bullet T_x^*M\to S^{\bullet+l}T_x^*M\) also denoted by \(\hat f\). For a
right inverse \(\widehat R\) to \(\hat\sigma\) and a homogeneous \(\widehat B\in
S^2(S^\bullet T_x^*M)\) of degree \(\hat d\) define the \textbf{polynomial Wendl
map} by
\[
\widehat W_{\widehat B}\colon S^\bullet T_x^*M\to
S^{\bullet+k+\hat d}T_x^*M,\qquad
\widehat W_{\widehat B}(\hat f):=
\widehat\Pi((\widehat R\hat f\otimes\mathbf 1)\widehat B).
\]

\begin{definition}\label{def:polywendl-subj}
We say that the symbol \(\sigma\) satisfies the \textbf{polynomial Wendl
condition} if for every choice of the map \(\widehat R\) and every \(\hat
d\in\mathbb N_0\) there exist constants \(\widehat C=\widehat C(\hat d),\widehat
L=\widehat L(\hat d)\) such that for all \(l\ge\widehat L\) the following holds:
for all nonzero \(\widehat B\in\ker\widehat\Pi\cap S^2\ker\hat\sigma\)
homogeneous of degree \(\hat d\) the map \(\widehat W_{\widehat B}^l\) satisfies
\[
\rk\widehat W_{\widehat B}^l\ge\widehat Cl^n.
\]
\end{definition}

Note that in view of Remark~\ref{rem:jetwendl-R} it is enough to check
Definition~\ref{def:polywendl-subj} for a particular choice of \(\widehat R\).

\begin{claim}\label{cl:polywendl-tojets}
Let \(\mathbf P\) be an elliptic differential operator on \(M\) and \(x\in M\).
If \(\sigma_x(\mathbf P)\) satisfies the polynomial Wendl condition, then
\(\mathbf P\) satisfies the jet Wendl condition at \(x\).
\end{claim}

\begin{proof}
Put \(\sigma:=\sigma_x(\mathbf P)\) and fix \(d\in\mathbb N_0,\tau\in\mathbb
N\). Let \(l\ge L\) and \(B\in\ker\Pi^{l+k}\cap S^2\ker\mathbf P^l\) be such
that \(\ord B\le d,\rk B=\tau\). Put \(\widehat B:=\lt B\in S^2(S^\bullet
T_x^*M)\). Then \(\widehat B\in\ker\widehat\Pi\), it is homogeneous, and its
degree \(\hat d\) satisfies \(\hat d\le d\). Furthermore,
\[
(\hat\sigma\otimes\mathbf 1)\widehat B=
((\mathbf P^l\otimes\mathbf 1)B)^{\hat d-k}=0,
\]
which implies that \(\widehat B\in S^2\ker\hat\sigma\). Let \(\widehat R\) be
any right inverse to \(\hat\sigma\). By Claim~\ref{cl:jets-compR} there exists a
right inverse \(R\) to \(\mathbf P^l\) such that \(R\) is of order \(k\) and
\(\gr_k R=\widehat R\). Let \(i\colon J_x^{l-\hat d}(M)\to J_x^l(M)\) be an
arbitrary right inverse to the natural projection \(J_x^l(M)\to J_x^{l-\hat
d}(M)\). Then the map \(W_B^l\circ i\) is of order \(k+\widehat d\) and it is
easy to see that the map \(\gr_{k+\hat d}(W_B^l\circ i)\) is identified with
\(\widehat W_{\widehat B}^{l-\hat d}\) under natural isomorphisms. In
particular, using the canonical inclusion
\[
\gr\ker (W_B^l\circ i)\hookrightarrow\ker\gr_{k+\hat d}(W_B^l\circ i)=
\ker\widehat W_{\widehat B}^{l-\hat d}
\]
and observing that the domains of both maps are of the same dimension, we get
\[
\rk W_B^l\ge\rk(W_B^l\circ i)\ge
\rk\widehat W_{\widehat B}^{l-\hat d}\ge
\widehat C(l-\hat d)^n\ge Cl^n,\qquad l\ge L.
\]
for some constants \(C,L\) depending only on \(d\) (note that
\(C,L\) are independent of \(\tau\) in this case).
\end{proof}

\subsection{Introducing coordinates}\label{sec:coord}


Let \(x\in M\) and let \(\sigma\in S^\bullet T_x M\) be an elliptic symbol of
order \(k\). Since in Sections~\ref{sec:coord}--\ref{sec:wendliff} we deal only
with polynomial Petri and Wendl maps, to simplify the notation, we omit all the
hats, i.e. we write simply \(\sigma,\Pi,B\) instead of
\(\hat\sigma,\widehat\Pi,\widehat B\), etc. till the end of
Section~\ref{sec:wendliff}.

Consider the {domain \(S^2(S^\bullet T_x^*M)\) of the map
\(\Pi\) as a subspace of \((S^\bullet T_x^*M)^{\otimes 2}\).}
It is convenient to extend \(\Pi\) to a larger domain by
\[
\Pi\colon (S^\bullet T_x^*M)^{\otimes 2}\to S^\bullet T_x^*M,\qquad
\Pi(u\otimes v)=uv.
\]
In particular, the map \(W_B\) can be defined for any \(B\in (S^\bullet
T_x^*M)^{\otimes 2}\).

Now let \(X=(X_1,\ldots,X_n)\) be an \(n\)-tuple of coordinate functions on
\(T_x M\) given by a basis of \(T_x^*M\). Then \(S^\bullet T_x^*M\cong\mathbb
R[X_1,\ldots,X_n]\), which we denote by \(\mathbb R[X]\). If
\(Y=(Y_1,\ldots,Y_n)\) is another \(n\)-tuple of coordinate functions on \(T_x
M\) obtained in the same way as \(X\), then we have a natural isomorphism
\[
(S^\bullet T_x^*M)^{\otimes 2}\cong\mathbb R[X_1,\ldots,X_n,Y_1,\ldots,Y_n]=:
\mathbb R[X,Y],
\]
so that \(B \in \mathbb R[X,Y]\). Employing the isomorphisms above, the map
\(\Pi\) (extended to a larger domain) becomes
\[
\Pi\colon\mathbb R[X,Y]\to\mathbb R[X],\qquad\Pi(B(X,Y))=B(X,X),
\] and the map \(W_B\) becomes
\[
W_B\colon\mathbb R[X]\to\mathbb R[X],\qquad W_B(f)=\Pi(R(fB)).
\]
Finally, consider the involution
\[
{}^\star\colon\mathbb R[X,Y]\to\mathbb R[X,Y],\qquad B^\star(X,Y):=B(Y,X).
\]
Clearly,
\[
S^2(S^\bullet T_x^*M)\cong\{B\in\mathbb R[X,Y]\colon B^\star=B\}.
\]
It is convenient to introduce the following notation. If \(C\in \mathbb
R[X,Y]\), then \(\partial_C\in\End(\mathbb R[X,Y])\) denotes the formal
differential operator on \(\mathbb R[X,Y]\) obtained from \(C\) by replacing all
\(X_i\)'s and \(Y_i\)'s by \(\partial_{X_i}\)'s and \(\partial_{Y_i}\)'s
respectively. In particular, there exists \(p\in\mathbb R[X]\) such that
\({\partial_p}=\sigma\). Tracking the isomorphisms used in this section, we see
that
\[
(\ker\sigma)^{\otimes 2}\cong\{B\in\mathbb R[X,Y]\colon\partial_p B=
\partial_{p^\star}B=0\}.
\]

\subsection{Weyl algebras}\label{sec:weyl}

Recall that the \textbf{(standard) Weyl algebra} \(A_n\) over \(\mathbb R\) is
the algebra over \(\mathbb R\) given by generators
\(X_1,\ldots,X_n,\partial_1,\ldots,\partial_n\) and relations
\begin{equation}\label{eq:weyl-strel}
[X_i,X_j]=0,\quad [\partial_i,\partial_j]=0,\quad
[\partial_i,X_j]=\delta_{ij},\qquad i,j=1,\ldots,n.
\end{equation}
Note that \(A_n\) comes with a natural action on \(\mathbb R[X]\) given by
\begin{equation}\label{eq:weyl-mu}
\mu\colon A_n\to\End(\mathbb R[X]),\qquad
\mu(X_i)u=X_i u,\quad\mu(\partial_i)u=\partial_{X_i}u.
\end{equation}
In what follows we need a generalization of this construction. Let \(p\in\mathbb
R[X]\) be nonzero. By definition, \(A_n[p^{-1}]\) is the \(\mathbb
R\)-algebra given by the generators
\(X_1,\ldots,X_n,p^{-1}\), \(\partial_1,\ldots,\partial_n\),
relations~\eqref{eq:weyl-strel}, and additional relations
\begin{equation}\label{eq:weyl-addrel}
pp^{-1}=p^{-1}p=1,\quad [X_i,p^{-1}]=0,\quad
[\partial_i,p^{-1}]=-p^{-2}\partial_{X_i}p,\qquad i=1,\ldots,n,
\end{equation}
where \(p^{-s}:=(p^{-1})^s,s\in\mathbb N_0\). Note that there is a natural
inclusion \(A_n\subset A_n[p^{-1}]\). Moreover, the action \(\mu\) extends to an
action of \(A_n[p^{-1}]\) on the localization \(\mathbb R[X,p^{-1}]\) also
denoted by \(\mu\).

For any \(K\in\mathbb\End(S^\bullet T_x^*M)\) the restriction of \(K\) onto
\(S^{\le l}T_x^*M\) is denoted by \(K^l\).
\begin{definition}
An operator \(K\in\mathbb\End(S^\bullet T_x^*M)\) is called
\begin{itemize}
\item\textbf{of bounded shift} if there exists \(s\in\mathbb N_0\) such that
\(\im K^l\subset S^{\le l+s}T_x^*M\) for all \(l\in\mathbb N_0\);
\item\textbf{negligible} if it is of bounded shift and there exist
positive constants \(C,L\) such that \(\dim \im K^l\le Cl^{n-1}\) for all \(l\ge
L\).
\end{itemize}
\end{definition}
It is easy to see that the subspace \(\mathcal{BS} \subset \End(S^\bullet
T_x^*M)\) of operators of bounded shift is a subalgebra in \(\End(S^\bullet
T_x^*M)\) containing all formal differential operators on \(S^\bullet T_x^*M\).
The subspace \(\mathcal N\subset \mathcal{BS}\) of negligible operators is a
two-sided ideal in \(\mathcal{BS}\). In particular, the quotient algebra
\[\mathcal{LR}=\mathcal{BS}/\mathcal N\] of operators \textbf{of large rank} is
well-defined. It follows from~\eqref{eq:jets-dimkersigma} that all right
inverses \(R\) to the elliptic symbol \(\sigma\) of order \(k\) induce the same
equivalence class in \(\mathcal{LR}\). 

\begin{proposition}\label{pr:weyl-phi}
There is a well-defined homomorphism of \(\mathbb R\)-algebras
\[
\phi\colon A_n[p^{-1}]\to \mathcal{LR}
\]
given by
\[
\phi(X_i):=[\mu(\partial_i)],\quad
\phi(\partial_i):=[-\mu(X_i)],\quad
\phi(p^{-1}):=[R],
\]
where \([\cdot]\) stands for the equivalence class in \(\mathcal{LR}\).
\end{proposition}

\begin{proof}
By the universal property of \(\mathbb R\)-algebras given by generators and
relations, it suffices to check that the map \(\phi\) respects the defining
relations~\eqref{eq:weyl-strel} and~\eqref{eq:weyl-addrel} of the algebra
\(A_n[p^{-1}]\). Let us check this, for example, for the relation
\[
[\partial_i,p^{-1}]=-p^{-2}\partial_{X_i}p,\qquad i=1,\ldots,n.
\]
We need to verify that
\begin{equation}\label{eq:weyl-phi0rel}
[\mu(X_i),R]-\partial_{\partial_{X_i}p}R^2\in\mathcal N,\qquad i=1,\ldots,n.
\end{equation}
Fix \(i\in\{1,\ldots,n\}\). We have
\begin{align*}
\sigma([\mu(X_i),R]-\partial_{\partial_{X_i}p}R^2)&=
\sigma\mu(X_i)R-\sigma R\mu(X_i)-
\sigma\partial_{\partial_{X_i}p}R^2\\&=
(\mu(X_i)+\partial_{\partial_{X_i}p}R)-\mu(X_i)-\partial_{\partial_{X_i}p}R=0,
\end{align*}
where we used that \([\sigma,\mu(X_i)]=\partial_{\partial_{X_i}p}\). It
follows that
\[
\im([\mu(X_i),R]-\partial_{\partial_{X_i}p}R^2)\subset\ker\sigma.
\]
By~\eqref{eq:jets-dimkersigma}, \(\dim\ker\sigma^l\) grows as \(l^{n-1}\).
This shows~\eqref{eq:weyl-phi0rel}. For the remaining relations the
verification is similar and easier.
\end{proof}

\begin{remark}
There are at least two different intuitive ways to think about the map \(\phi\).
The first one is the Fourier transform. Indeed, the Fourier transforms of the
operators \(\mu(X_j),\mu(\partial_j),\mu(p),\sigma\) extended to the space of
tempered distributions \(\mathcal S'(\mathbb R^n)\) are
\(i\mu(\partial_j),i\mu(X_j),i^k\sigma,i^k p\) respectively. One can get rid of
the ``\(i\)'' factors by pre-composing with the automorphism of \(A_n^{\mathbb
C}\) given by \(X_j\mapsto iX_j,\partial_j\mapsto -i\partial_j\). The result is
essentially the map \(\phi|_{A_n}\). In particular, this suggests (by analogy
with solving a linear PDE with constant coefficients) to put
\(\phi(p^{-1}):=[R]\) since the latter map is the inverse of \([\sigma]\).

The second way comes from the following standard trick: define a dot product
\((\cdot,\cdot)\) on \(\mathbb R[X]\) by \((u,v):=(\partial_u v)(0)\). This dot
product has the following property: the maps \(\mu(X_j)\) and
\(\mu(\partial_j)\) are adjoint. In particular, the maps \(\sigma\) and
\(\mu(p)\) are adjoint, which means that the hypothetical adjoint to \(R\) is a
kind of ``multiplication by \(p^{-1}\)''. One can consider the map \(\phi\) as
an attempt to formalize this vague idea.
\end{remark}

Define an automorphism
\[
\iota\colon A_n\to A_n,\qquad\iota(X_i)=\partial_i,\quad
\iota(\partial_i)=-X_i.
\]
Note that \(\phi(D)=[\mu\iota(D)]\) for all \(D\in A_n\) for \(\phi,\mu\) from
Proposition~\ref{pr:weyl-phi} and~\eqref{eq:weyl-mu}.

%

\begin{proposition}\label{pr:weyl-phiLR}
\begin{enumerate}[label=(\roman*)]
\item
The map \(\phi\) of Proposition~\ref{pr:weyl-phi} is injective.
\label{it:weyl-phinj}

\item
Let \(\widetilde V\subset\mathcal{BS}\) be a finite dimensional vector subspace
such that \(\widetilde V/(\widetilde V\cap\mathcal N)\subset\im\phi\). Then
there exist constants \(\widetilde C,\widetilde L\) depending only on
\(\widetilde V\) such that for all \(l\ge\widetilde L\) and \(\widetilde
D\in\widetilde V\setminus\mathcal N\) we have \(\rk\widetilde D^l\ge\widetilde
Cl^n\).
\label{it:weyl-CL}
\end{enumerate}
\end{proposition}

\begin{proof}
\ref{it:weyl-phinj} Let \(D\in A_n[p^{-1}]\) be nonzero. Using the
relations~\eqref{eq:weyl-addrel}, we can ``pull out'' all the \(p^{-1}\) factors
in each term of \(D\). In particular, there exists \(s\in\mathbb N_0\) such that
\(p^s D\in A_n\). Put \(D':=\iota(p^s D)\in A_n\). We have
\[
\phi(D)=\phi(p^{-s})\phi(p^s D)=[R^s]\cdot[\mu(D')]=[R^s\mu(D')].
\]
Since \(D\ne 0\), we have \(D'\ne 0\). Using the
relations~\eqref{eq:weyl-strel}, one can write \(D'\) uniquely in the form
\[
D'=\sum_{\alpha,\beta\in\mathbb N_0^n}
a_{\alpha,\beta}X^\alpha\partial^\beta,\qquad a_{\alpha,\beta}\in\mathbb R.
\]
Put
\[
h:=\max\{|\alpha|+|\beta|\colon\alpha,\beta\in\mathbb N_0^n,
a_{\alpha,\beta}\ne 0\},\qquad N:=2h+1,
\]
and for any \(l\in\mathbb N_0\) put
\[
l_0:=\Bigl[\frac lN\Bigr],\qquad
S^l:=\spn\{X^{N\gamma}\mid\gamma\in\mathbb N_0^n,|\gamma|\le l_0\}.
\]
Observe that
\[
\dim S^l=\sum_{j=0}^{l_0}\binom{j+n-1}{n-1}=\binom{l_0+n}n\ge
\frac{l^n}{n!N^n}+O(l^{n-1}).
\]
The space \(S^l\) is constructed so that
\[
\dim\ker\mu(D')|_{S^l}=\#\{\gamma\in\mathbb N_0^n\colon |\gamma|\le l_0,\quad
\mu(D')X^{N\gamma}=0\}.
\]
Since \(D'\ne 0\), there exist \(\alpha,\beta\in\mathbb N_0^n\) such that
\(a_{\alpha,\beta}\ne 0\). Take such a pair \(\alpha,\beta\) and put
\(\delta:=\beta-\alpha\in\mathbb Z^n\). Then for \(\gamma\in\mathbb N_0^n\), the
coefficient of \(X^{N\gamma-\delta}\) in \(\mu(D')X^{N\gamma}\) is a nonzero
polynomial in \(\gamma\). It is easy to see that any such polynomial has
\(O(l_0^{n-1})=O(l^{n-1})\) zeros on the set \(\{\gamma\in\mathbb
N_0^n\colon |\gamma|\le l_0\}\) as \(l,l_0\to+\infty\). Hence,
\[
\rk\mu(D')^l\ge\rk\mu(D')|_{S^l}=\dim S^l-\dim\ker\mu(D')|_{S^l}\ge Cl^n,
\qquad l\ge L
\]
for \(C=\frac 1{2n! N^n}\) and \(L\) sufficiently large. Since \(R^s\) is
injective, we obtain that the operator \(R^s\mu(D')\) is not negligible, hence
\(\phi(D)=[R^s\mu(D')]\ne 0\).

\ref{it:weyl-CL} For each \(\widetilde D\in\widetilde V\) define
\(D:=\phi^{-1}([\widetilde D])\). Let \(s\in\mathbb N_0\) be such that \(p^s
D\in A_n\) for all \(\widetilde D\in\widetilde V\) and let \(\widetilde
D_*:=R^s(\mu\iota(p^s D))\). The same argument as in the proof of
part~\ref{it:weyl-phinj} above shows that there exist constants \(C,L\) such
that
\[
\rk\widetilde D_*^l\ge Cl^n,\qquad\widetilde D\in\widetilde V\setminus\mathcal N,
\quad l\ge L.
\]
Note that the constants \(C,L\) can be chosen to be independent of \(\widetilde
D\) because \(\widetilde V\) is finite dimensional. Since \([\widetilde
D_*]=\phi(D)=[\widetilde D]\), the operator \(\widetilde D_*-\widetilde D\) is
negligible and depends linearly on \(\widetilde D\). Hence there exist constants
\(C',L'\) such that
\[
\rk(\widetilde D_*^l-\widetilde D^l)\le C'l^{n-1},\qquad\widetilde D\in
\widetilde V,\quad l\ge L'.
\]
Again, the constants \(C',L'\) can be chosen to be independent of \(\widetilde
D\) because \(\widetilde V\) is finite dimensional. Thus,
\[
\rk\widetilde D^l\ge\rk\widetilde D_*^l-\rk(\widetilde D_*^l-\widetilde D^l)\ge
Cl^n-C'l^{n-1}\ge\widetilde Cl^n,\qquad\widetilde D\in\widetilde V\setminus
\mathcal N,\quad l\ge\widetilde L
\]
for some new constants \(\widetilde C,\widetilde L\), as desired.
\end{proof}

\subsection{A criterion for the polynomial Wendl condition}\label{sec:wendliff}

Let \(\alpha,\beta\in\mathbb N_0^n\). We have
\[
\Pi(R(X^\alpha Y^\beta f))=\Pi(Y^\beta R(X^\alpha f))=X^\beta R(X^\alpha f)=
\mu(X^\beta)R\mu(X^\alpha)f,\qquad f\in\mathbb R[X].
\]
Taking the quotient by \(N\), we obtain
\[
[W_{X^\alpha Y^\beta}]={(-1)^{\alpha+\beta}} \phi(\partial^\beta
p^{-1}\partial^\alpha).
\]
For each \(q\in\mathbb R[X,p^{-1}]\) homogeneous and nonzero we define {the
linear map}
\[
\omega_q\colon\mathbb R[X,Y]\to A_n[p^{-1}],\qquad
\omega_q(X^\alpha Y^\beta):=(-1)^{\alpha+\beta}\partial^\beta q\partial^\alpha.
\]
It follows that
\begin{equation}\label{eq:wendliff-Weqomega}
[W_B]=\phi(\omega_{p^{-1}}(B)),
\end{equation}
where \(B\in\mathbb R[X,Y]\) is homogeneous of degree \(d\).

%

\begin{proposition}\label{pr:wendliff-omegainj}
The symbol \(\sigma\) satisfies the polynomial Wendl condition if and only if
\(\omega_{p^{-1}}(B)\ne 0\) for every homogeneous nonzero \(B\in\ker\Pi\cap
(S^2\ker\sigma)\).
\end{proposition}

\begin{proof}
Follows directly from~\eqref{eq:wendliff-Weqomega} and
Proposition~\ref{pr:weyl-phiLR}\ref{it:weyl-CL} with
\[
\widetilde V=\{W_B\mid B\in\ker\Pi\cap (S^2\ker\sigma)\text{ homogeneous of
degree }d\},\qquad d\in\mathbb N_0.
\]
\end{proof}

To study the injectivity of \(\omega_{p^{-1}}\), we need a few auxiliary claims.
In the following claim we use the identification \(\sigma=\partial_p\).

\begin{claim}\label{cl:wendliff-Bmin}
Suppose that \(B\in\mathbb R[X,Y]\) is homogeneous of degree \(d>0\) and
satisfies
\begin{equation}\label{eq:wendliff-Bmin}
(\partial_{X_i}+\partial_{Y_i})B=0,\qquad i=1,\ldots,n.
\end{equation}
Then \(B(X,Y)=H(X-Y)\) for some \(H\in\mathbb R[X]\) homogeneous of degree
\(d\). In particular, \(B^\star=(-1)^d B,\,\Pi(B)=0\), and
\(B\in(\ker\sigma)^{\otimes 2}\) if and only if \(H\in\ker\sigma\).
\end{claim}

\begin{proof}
After a linear change of coordinates \(\widetilde X:=X-Y,\,\widetilde Y:=X+Y\),
condition~\eqref{eq:wendliff-Bmin} reads as \(\partial_{\widetilde Y_i}B=0\) for
all \(i=1,\ldots,n\). Therefore \(B(\widetilde X,\widetilde Y)=H(\widetilde
X)=H(X-Y)\) for some \(H\in\mathbb R[X]\), which is automatically homogeneous of
degree \(d\). The remaining assertions are clear.
\end{proof}

\begin{claim}\label{cl:wendliff-omega}
Let \(q\in\mathbb R[X,p^{-1}]\) be homogeneous and nonzero and \(B,H\) be as in
Claim~\ref{cl:wendliff-Bmin}. Then \(\omega_q(B)=\partial_H q\).
\end{claim}

\begin{proof}
Let \(H=\sum\limits_{|\alpha|=d}h_\alpha X^\alpha\). Then
\[
B=\sum_{|\alpha|=d}\sum_{\gamma\le\alpha}
(-1)^{|\gamma|}\binom{\alpha}\gamma h_\alpha X^{\alpha-\gamma}Y^\gamma,
\]
and we calculate as follows:
\[
\omega_q(B)=(-1)^d\sum_{|\alpha|=d}h_\alpha\sum_{\gamma\le\alpha}
(-1)^{|\gamma|}\binom{\alpha}\gamma\partial^\gamma q\partial^{\alpha-\gamma}=
\sum_{|\alpha|=d}h_\alpha\partial_{X^\alpha}q=\partial_H q,
\]
where
\[
\binom{\alpha}\gamma:=\prod_{i=1}^n\binom{\alpha_i}{\gamma_i}.
\]
Indeed, the following identity holds for each \(\alpha\in\mathbb
N_0^n,|\alpha|=d\):
\begin{align*}
h_\alpha\sum_{\gamma\le\alpha}(-1)^{|\gamma|}\binom{\alpha}
\gamma\partial^\gamma q\partial^{\alpha-\gamma}&=
h_\alpha\sum_{\gamma'\le\gamma\le\alpha}
(-1)^{|\gamma|}\binom{\alpha}\gamma\binom{\gamma}{\gamma'}
(\partial_{X^{\gamma'}}q)\partial^{\alpha-\gamma'}\\&=
h_\alpha\sum_{\gamma'\le\gamma\le\alpha}
(-1)^{|\gamma|}\binom{\alpha}{\gamma'}\binom{\alpha-\gamma'}{\gamma-\gamma'}
(\partial_{X^{\gamma'}}q)\partial^{\alpha-\gamma'}\\&=
h_\alpha\sum_{\gamma'\le\alpha}\binom{\alpha}{\gamma'}
\left(\sum_{\gamma'\le\gamma\le\alpha}
(-1)^{|\gamma|}\binom{\alpha-\gamma'}{\gamma-\gamma'} \right)
(\partial_{X^{\gamma'}}q)\partial^{\alpha-\gamma'}\\&=
(-1)^d h_\alpha\partial_{X^\alpha}q,
\end{align*}
since \(\gamma'=\alpha\) is the only \(\gamma'\) for which the sum in the
brackets is nonzero. 
\end{proof}

\begin{proposition}\label{pr:wendliff-main}
Let \(q\in\mathbb R[X,p^{-1}]\) be homogeneous and nonzero. Then there exists
\(B\in\ker\Pi\cap S^2\ker\sigma\) homogeneous and nonzero such that
\(\omega_q(B)=0\) if and only if there exists nonzero \(H\in\ker\sigma\)
homogeneous of positive even degree such that \(\partial_H q=0\).
\end{proposition}

\begin{proof}
We begin with the ``only if'' part. Suppose that there exists \(B\in\ker\Pi\cap
S^2\ker\sigma\) nonzero and homogeneous of degree \(d\) such that
\(\omega_q(B)=0\). Let \(d\) be minimal possible. Note that \(\Pi(B)=0\)
implies that \(d>0\). We have
\begin{equation}\label{eq:wendliff-comm}
0=[X_i,\omega_q(B)]=\omega_q((\partial_{X_i}+\partial_{Y_i})B),\qquad
i=1,\ldots,n.
\end{equation}
Also, it is easy to see that \((\partial_{X_i}+\partial_{Y_i})B\in\ker\Pi\cap
S^2\ker\sigma\). The minimality condition imposed on \(d\)
implies~\eqref{eq:wendliff-Bmin}. By Claim~\ref{cl:wendliff-Bmin},
\(B(X,Y)=H(X-Y)\) for some \(H\in\ker\sigma\). Since \(B\) is symmetric, we have
\(B^\star=B\). On the other hand, \(B^\star=(-1)^d B\) by
Claim~\ref{cl:wendliff-Bmin}, which implies that \(d\) is even. Finally,
\(\partial_H q=\omega_q(B)=0\) by Claim~\eqref{cl:wendliff-omega}.

The ``if'' part follows directly from Claims~\ref{cl:wendliff-Bmin}
and~\ref{cl:wendliff-omega}. Indeed, if there exists \(H\in\ker\sigma\)
homogeneous of even degree \(d>0\) such that \(\partial_H q=0\), then
\(B(X,Y):=H(X-Y)\) lies in \(\ker\Pi\cap S^2\ker\sigma\) by
Claim~\ref{cl:wendliff-Bmin} and satisfies \(\omega_q(B)=\partial_H q=0\) by
Claim~\ref{cl:wendliff-omega}.
\end{proof}

\begin{corollary}\label{cor:wendliff-main}
The symbol \(\sigma\) satisfies the polynomial Wendl condition if and only if
there does not exist a nonzero element \(H\in\ker\sigma\) homogeneous of
positive even degree such that \(\partial_H(p^{-1})=0\).
\end{corollary}

\begin{proof}
Follows directly from Proposition~\ref{pr:wendliff-main} with \(q=p^{-1}\)
and Proposition~\ref{pr:wendliff-omegainj}.
\end{proof}

\begin{example}\label{ex:wendliff-fail}
One can use Corollary~\ref{cor:wendliff-main} to construct an example of a
symbol that \emph{does not} satisfy the polynomial Wendl condition. Let \(m\ge
3,n:=2m,p:=(X_1^2+\ldots+X_n^2)^{m-1}\) (so that
\(\sigma=(\partial_{X_1}^2+\ldots+\partial_{X_n}^2)^{m-1}=\Delta^{m-1}\) is a
power of the Laplacian), and \(H:=X_1^2+\ldots+X_n^2\). We have
\(H\in\ker\sigma\) since \(\deg H<\deg p\) and
\(\partial_H(p^{-1})=\Delta(p^{-1})=0\) since \(p^{-1}\) is the fundamental
solution for the Laplace equation on \(\mathbb R^n\), as \(n=2m\). It follows
from Corollary~\ref{cor:wendliff-main} that \(\sigma\) does not satisfy the
polynomial Wendl condition.
\end{example}

\begin{proposition}\label{pr:wendliff-lap}
Any scalar elliptic symbol of order two satisfies the polynomial Wendl
condition.
\end{proposition}

\begin{proof}
Let \(\sigma\colon S^\bullet T_x^*M\to S^\bullet T_x^*M\) be the formal
differential operator corresponding to a scalar symbol of order two. By an
appropriate choice of coordinates and multiplication of the symbol by \(-1\) if
necessary, we may assume that \(\sigma=\partial_p\) where
\(p=X_1^2+\ldots+X_n^2\), i.e. \(\sigma=\Delta\), the Laplacian. By
Corollary~\ref{cor:wendliff-main} (cf. Example~\ref{ex:wendliff-fail}) it
suffices to check that \(\partial_H(p^{-1})\ne 0\) for every \(H\in\mathbb
R[X]\) homogeneous of even degree \(d>0\) and satisfying \(\Delta H=0\) (i.e.,
\(H\) is harmonic). By Fa\`a di Bruno's formula, for each \(\gamma\in\mathbb
N_0^n\) we have 
\begin{align*}
\partial_{X^\gamma}(p^{-1})&=
\sum_{\pi\in P_\gamma}(y^{-1})^{(|\pi|)}|_{y=p}
\prod_{\gamma'\in\pi}\partial_{X^{\gamma'}}p\\&=
\sum_{\pi\in P_\gamma}(-1)^{|\pi|}|\pi|!\,p^{-|\pi|-1}
\prod_{\gamma'\in\pi}\partial_{X^{\gamma'}}p\\&=
p^{-|\gamma|-1}\left(pf_\gamma+(-1)^{|\gamma|}|\gamma|!
\prod_{i=1}^n(\partial_{X_i}p)^{\gamma_i}\right)\\&=
p^{-|\gamma|-1}(pf_\gamma+(-2)^{|\gamma|}|\gamma|!\,X^\gamma)\qquad
\text{for some }f_\gamma\in\mathbb R[X],
\end{align*}
where \(P_\gamma\) is the multiset consisting of all partitions of the multiset
\[
\{1^{\gamma_1},\ldots,n^{\gamma_n}\}:=
\{\underbrace{1,\ldots,1}_{\gamma_1},\ldots,
\underbrace{n,\ldots,n}_{\gamma_n}\}
\]
counting with multiplicities, and a multiindex \(\gamma'\in\mathbb N_0^n\) is
identified with the multiset \(\{1^{\gamma_1'},\ldots,n^{\gamma_n'}\}\). Writing
\(\partial_H=\sum_{|\gamma|=d}h_\gamma\partial_{X^\gamma}\), it follows
from the above identity that
\begin{equation}\label{eq:wendliff-dibruno}
p^{d+1}\partial_H(p^{-1})\equiv (-2)^d d!\,H\pmod p.
\end{equation}
However, it is well-known (see, e.g.,~\cite[Proposition~5.3]{ABR1992}) that a
nonzero harmonic polynomial cannot be divisible by \(p\). This implies that
\(\partial_H(p^{-1})\ne 0\) whenever \(H\ne 0\), as desired.
\end{proof}

\subsection{Traceless \texorpdfstring{\(\Delta\)}{Δ}-Wendl
conditions}\label{sec:tdwendl}

Let a point \(x\in M\) be fixed and \(l\in\mathbb N\). For every
\((g,\lambda)\in\mathcal M^*\), an element \(B\in S^2 J_x^{l+2}(M)\), and a
right inverse \(R\colon J_x^l(M)\to J_x^{l+2}(M)\) to \(\widetilde{\mathbf
T}_{g,\lambda}^l=\Delta_g^l-\lambda\) define the \textbf{jet traceless
\(\Delta\)-Wendl map} as follows:
\begin{gather*}
W_B^{\perp;l}\colon J_x^{l+1}(S^2 T^*M)\to J_x^{l+1}(S^2 T^*M),\\
W_B^{\perp;l}(h):=
2\Pi_g^{\perp;l+1}(((\d^{l+1}R\div_g^l(.\iprod h))\odot\mathbf 1)B_\d)+
\frac 1n\tr_g^{l+1}(\Pi_\d^{l+1}(B_\d))h,
\end{gather*}
where \(B_\d:=(\d^{l+1}\otimes\d^{l+1})B\).

\begin{definition}\label{def:tdwendl-jetw}
We say that the \textbf{jet traceless \(\Delta\)-Wendl condition} is satisfied
at \((g,\lambda)\in\mathcal M^*\) if for every choice of the map \(R\) and every
\(d\in\mathbb N_0,\tau\in\mathbb N\) there exist constants
\(C=C(d,\tau),L=L(d,\tau)\) such that for all \(l\ge L\) the following holds:
for all \(B\in S^2\ker\widetilde{\mathbf T}_{g,\lambda}^l\) such that
\(B_\d\in\ker\Pi_g^{\perp;l+1}\) is nonzero and satisfies \(\ord B_\d\le d,\rk
B_\d=\tau\), the restriction of the map \(W_B^{\perp;l}\) onto the space
\(\{h\in J_x^{l+1}(S^2 T^*M)\colon\tr_g^{l+1}h=0\}\) has rank at least \(Cl^n\).
\end{definition}

The maps \(\d,\div_g,\tr_g\) induce the corresponding polynomial maps,
\begin{align*}
\hat{\d}&:=\gr_{-1}\d^\infty\colon S^\bullet T_x^*M\to
S^{\bullet-1}T_x^*M\otimes T_x^*M,\\
\widehat{\div}_g&:=\gr_{-1}\div_g^\infty\colon S^\bullet T_x^*M\otimes T_x^*M\to
S^{\bullet-1}T_x^*M,\\
\widehat{\tr}_g&:=\gr_0\tr_g^\infty\colon S^\bullet T_x^*M\otimes S^2 T_x^*M\to
S^\bullet T_x^*M,\\
\widehat\Delta_g&:=\gr_{-2}\Delta_g^\infty\colon S^\bullet T_x^*M\to
S^\bullet T_x^*M.
\end{align*}
Of course, \(\widehat\Delta_g=-\widehat{\div}_g\circ\hat{\d}\). The
\textbf{polynomial traceless \(\Delta\)-Petri map} is defined by
\begin{gather*}
\widehat\Pi_g^\perp\colon S^2(S^\bullet T_x^*M\otimes T_x^*M)\to
S^\bullet T_x^*M\otimes S^2 T_x^*M,\\
\widehat\Pi_g^\perp(\widehat B_\d):=\widehat\Pi_\d(\widehat B_\d)-
\frac{1}n\widehat{\tr}_g\widehat\Pi_\d(\widehat B_\d)\otimes g_x,
\end{gather*}
where
\begin{gather*}
\widehat\Pi_\d\colon S^2(S^\bullet T_x^*M\otimes T_x^*M)\to
S^\bullet T_x^*M\otimes S^2 T_x^*M,\\
\widehat\Pi_\d((u\otimes\theta)\odot
(v\otimes\omega)):=uv\otimes(\theta\odot\omega).
\end{gather*}
For any \(\widehat B\in S^2(S^\bullet T_x^*M)\) homogeneous of degree \(\hat
d+2\) define the \textbf{polynomial traceless \(\Delta\)-Wendl map} by
\begin{gather*}
\widehat W_{\widehat B}^\perp\colon S^\bullet T_x^*M\otimes S^2 T_x^*M\to
S^{\bullet+\hat d}T_x^*M\otimes S^2 T_x^*M,\\
\widehat W_{\widehat B}^\perp(\hat h):=
2\widehat\Pi_g^\perp(((\hat{\d}\,\widehat R\,\widehat{\div}_g(.\iprod\hat h))
\odot\mathbf 1)\widehat B_\d)+\frac{1}n\widehat{\tr}_g(\widehat\Pi_\d(\widehat
B_\d))\otimes\hat h.
\end{gather*}
Put
\[
\kere\widehat\Delta_g:=\{u\in S^\bullet T_x^*M\colon\widehat\Delta_g u=\const\}.
\]

\begin{definition}\label{def:tdwendl-polyw}
We say that the \textbf{polynomial traceless \(\Delta\)-Wendl condition} is
satisfied at \((g,\lambda)\in\mathcal M^*\) if for every choice of the map
\(\widehat R\) and every \(\hat d\in\mathbb N_0\) there exist constants
\(\widehat C=\widehat C(\hat d),\widehat L=\widehat L(\hat d)\) such that for
all \(l\ge\widehat L\) the following holds: for all \(\widehat B\in
S^2\kere\widehat\Delta_g\) homogeneous of degree \(\hat d+2\) and such that
\(\widehat B_\d\in\ker\widehat\Pi_g^\perp\) is nonzero, the restriction of the
map \(\widehat W_{\widehat B}^\perp\) onto the space \(\{\hat h\in S^{\le
l+1}T_x^*M\otimes S^2 T_x^*M\colon\widehat{\tr}_g\,\hat h=0\}\) has rank at
least \(\widehat Cl^n\).
\end{definition}

\begin{claim}\label{cl:tdwendl-tojets}
If the polynomial traceless \(\Delta\)-Wendl condition is satisfied at
\((g,\lambda)\in\mathcal M^*\), then the jet traceless \(\Delta\)-Wendl
condition is satisfied at \((g,\lambda)\).
\end{claim}

\begin{proof}
The proof is similar to the proof of Claim~\ref{cl:polywendl-tojets}. Let \(B\)
be as in Definition~\ref{def:tdwendl-jetw}. Then the lowest order term \(\lt
B_\d\) is homogeneous of degree \(\hat d\le d\) and lies in
\(\im\hat{\d}^{\otimes 2}\cap\ker\widehat\Pi_g^\perp\). In particular, there
exists \(\widehat B\in S^2(S^\bullet T_x^*M)\) homogeneous of degree \(\hat
d+2\) such that \(\widehat B_\d=\lt B_\d\). The crucial point is to show that
one can choose \(\widehat B\in S^2\kere\widehat\Delta_g\), i.e. that \(\lt
B_\d\in S^2(\hat{\d}\kere\widehat\Delta_g)\). This can be seen with the help of
the \textbf{Hodge Laplacian}
\[
\Delta_g^\d\colon C^\infty(T^*M)\to C^\infty(T^*M),\qquad
\Delta_g^\d:=-(\div_g\circ\,\d+\d\circ\div_g).
\]
It is easy to see that \(\d\circ\Delta_g=\Delta_g^\d\circ\d\), which implies
that \(\im\hat{\d}\cap\ker\widehat\Delta_g^\d=\hat{\d}\kere\widehat\Delta_g\)
and \(B_\d\in S^2\ker(\Delta_g^{\d;l-1}-\lambda)\). Hence,
\[
(\widehat\Delta_g^\d\otimes\mathbf 1)\lt B_\d=
(((\Delta_g^{\d;l-1}-\lambda)\otimes\mathbf 1)B_\d)^{\hat d-2}=0,
\]
which implies that
\[
\lt B_\d\in\im\hat{\d}^{\otimes 2}\cap S^2\ker\widehat\Delta_g^\d =
S^2(\im\hat{\d}\cap\ker\widehat\Delta_g^\d)=S^2(\hat{\d}\kere\widehat\Delta_g),
\]
as desired.

The rest of the argument is almost the same as in the proof of
Claim~\ref{cl:polywendl-tojets}. The only additional subtle point is that
Definitions~\ref{def:tdwendl-jetw} and~\ref{def:tdwendl-polyw} deal with
traceless tensors. To address this issue, let
\[
\pi\colon J_x^{l+1}(S^2 T^*M)\to J_x^{l+1}(S^2 T^*M),\qquad
\pi(h):=h-\frac 1n(\tr_g^{l+1}h)j_x^{l+1}g
\]
be the projection onto the space of traceless jets, and \(\hat\pi:=\gr_0\pi\).
Then just replace \(W_B^l\circ i\) and \(\widehat W_{\widehat B}^{l-\hat d}\)
from the proof of Claim~\ref{cl:polywendl-tojets} by
\(W_B^{\perp;l}\circ\pi\circ i\) and \(\widehat W_{\widehat B}^{\perp;l+1-\hat
d}\circ\hat\pi\) respectively (where now \(i\colon J_x^{l+1-\hat d}(S^2 T^*M)\to
J_x^{l+1}(S^2 T^*M)\)).
\end{proof}

Let a \(g\)-orthonormal basis \(X=(X_1,\ldots,X_n)\) of \(T_x^*M\) be chosen.
It is convenient to write an element \(u\in S^\bullet T_x^*M\otimes T_x^*M\)
in the form \(u=\sum_{i=1}^n u_i\,\d X_i\), where \(u_i\in\mathbb R[X]\) and
\(\d X=(\d X_1,\ldots,\d X_n)\) is just a copy of \(X\) (one should think of the
second instance of \(T_x^*M\) in \(S^\bullet T_x^*M\otimes T_x^*M\) as of a
vector bundle \(E\), which accidentally coincides with \(T_x^*M\)). Similarly,
an element \(h\) of \(S^\bullet T_x^*M\otimes S^2 T_x^*M\) can be written in the
form \(\sum_{i,j=1}^n h_{ij}\,\d X_i\otimes\d X_j\), where \(h_{ij}\in\mathbb
R[X]\) and \(h_{ij}=h_{ji}\) for all \(i,j=1,\ldots,n\). It is easy to see that
with this interpretation, \(\hat{\d},\widehat{\div}_g,\widehat{\tr}_g\) coincide
with the standard differential, divergence, and trace respectively. In
particular, they are all independent of \(g\) (more precisely, the dependence is
incorporated into the choice of the orthonormal basis \(X\)), and we omit it in
the sequel. Also, similarly to Sections~\ref{sec:coord}--\ref{sec:wendliff}, we
omit hats and write simply \(\d,\div,\tr,\Pi^\perp\) etc. till the end of this
section.

Let \(B\in S^2\kere\Delta\subset\mathbb R[X,Y]\) and
\(i,i',j,j'\in\{1,\ldots,n\}\) be such that \(i<i'\) and \(j<j'\). Suppose that
\(h=\sum_{k,l=1}^n h_{kl}\d X_k\otimes\d X_l\in S^\bullet T_x^*M\otimes S^2
T_x^*M\) satisfies \(h_{kl}=0\) for all \(\{k,l\}\ne\{i,i'\}\). Then a direct
computation shows that the coefficient of \(\d X_{j'}\otimes\d X_j\) in
\(W_B^\perp(h)\) is equal to \(W_{B;ii';jj'}^\perp(h_{ii'})\), where
\begin{align*}
W_{B;ii';jj'}(f)&:=
\begin{multlined}[t]
w_{B;ii'jj'}^\perp(f)+w_{B;i'ijj'}^\perp(f)+
w_{B;ii'j'j}^\perp(f)+w_{B;i'ij'j}^\perp(f)+\\
\frac 1n\delta_{ij}\delta_{i'j'}\tr(\Pi_\d(B_\d))f,
\end{multlined}\\
w_{u\otimes v;ijkl}^\perp(f)&:=(\partial_{X_k}R\partial_{X_j}
((\partial_{X_i}u)f))\partial_{X_l}v,\qquad f,u,v\in\mathbb R[X].
\end{align*}
(We use the identification \(\mathbb R[X,Y]\cong\mathbb R[X]^{\otimes 2}\).) In
the terminology of Sections~\ref{sec:weyl}-\ref{sec:wendliff}, we have
\[
(\partial_{X_k}R\partial_{X_j}((\partial_{X_i}u)f))\partial_{X_l}v=
(\mu(\partial_{X_l}v)\circ\mu(\partial_k)\circ R\circ\mu(\partial_j)\circ
\mu(\partial_{X_i}u))f.
\]
Taking the quotient by negligible operators, we obtain
\begin{align*}
[w_{u\otimes v;ijkl}^\perp]&=
[\mu(\partial_{X_l}v)]\cdot[\mu(\partial_k)]\cdot [R]\cdot
[\mu(\partial_j)]\cdot [\mu(\partial_{X_i}u)]\\&=
-(-1)^{\deg v-1}\phi(\partial_{\partial_{X_l}v})\cdot\phi(X_k)\cdot
\phi(p^{-1})\cdot\phi(X_j)\cdot
(-1)^{\deg u-1}\phi(\partial_{\partial_{X_i}u})\\&=
-\phi(\omega_{p^{-1}X_j X_k}(\partial_{X_i}u\otimes\partial_{X_l}v)),
\end{align*}
where \(p:=X_1^2+\ldots+X_n^2\) (note that \(\phi(p^{-1})=-[R]\) in this case).
Hence,
\[
[w_{B;ijkl}^\perp]=-\phi(\omega_{p^{-1}X_j X_k}(\partial_{X_i Y_l}B)).
\]
Put
\begin{multline*}
\omega_{ii';jj'}^\perp(B):=-(\omega_{p^{-1}X_{i'} X_j}(\partial_{X_i Y_{j'}}B)+
\omega_{p^{-1}X_i X_j}(\partial_{X_{i'}Y_{j'}}B)+
\omega_{p^{-1}X_{i'}X_{j'}}(\partial_{X_i Y_j}B)\,+\\
\omega_{p^{-1}X_i X_{j'}}(\partial_{X_{i'}Y_j}B))+
\frac{(-1)^d}n\delta_{ij}\delta_{i'j'}\partial_{\tr(\Pi_\d(B_\d))}.
\end{multline*}
Then 
\[
[W_{B;ii';jj'}^\perp]=\phi(\omega_{ii';jj'}^\perp(B)).
\]
Using this and applying Proposition~\ref{pr:weyl-phiLR}\ref{it:weyl-CL} with
\begin{gather*}
\widetilde V=\{W_{B;ii';jj'}^\perp\mid B\in
S^2\kere\Delta\text{ homogeneous of degree }d+2,\,B_\d\in\ker\Pi^\perp\},\\
d\in\mathbb N_0,\qquad i,i',j,j'=1,\ldots,n,\quad i<i',\,j<j',
\end{gather*}
we obtain that the polynomial traceless \(\Delta\)-Wendl condition is satisfied
at every \((g,\lambda)\in\mathcal M^*\) whenever for every homogeneous nonzero
\(B\in S^2\kere\Delta,B_\d\in\ker\Pi^\perp\setminus\{0\}\) there exist
\(i,i',j,j'\in\{1,\ldots,n\}\) satisfying \(i<i',\,j<j'\) and such that
\(\omega_{ii';jj'}^\perp(B)\ne 0\).

We are in a position to prove

\begin{proposition}\label{pr:tdwendl-main}
Let \(n>2\). Then the polynomial traceless \(\Delta\)-Wendl condition is
satisfied at every \((g,\lambda)\in\mathcal M^*\).
\end{proposition}

\begin{proof}
Assume the converse. Then (see the discussion above) there exists an element
\(B\in S^2\kere\Delta,B_\d\in\ker\Pi^\perp\setminus\{0\}\) such that
\begin{equation}\label{eq:tdwendl-iijj}
\omega_{ii';jj'}^\perp(B)=0,\qquad i,i',j,j'=1,\ldots,n,\quad i<i',\,j<j'.
\end{equation}
Assume that \(d:=\deg B-2\) is minimal possible. If \(d=0\), then it follows
from \(\Pi^\perp(B_\d)=0\) that
\[
B_\d=c\sum_{k=1}^n\d X_k\,\d Y_k,\qquad c\in\mathbb R.
\]
We may assume that \(c=1\). Then
\[
B=P(X)+Q(Y)+\sum_{k=1}^n X_k Y_k
\]
for some polynomials \(P,Q\) of degree~2. In particular, \(\partial_{X_i
Y_j}B=\delta_{ij}\) for all \(i,j=1,\ldots,n\), and we have
\[
0=\omega_{13;23}^\perp(B)=-\omega_{p^{-1}X_1 X_2}(1)=-p^{-1}X_1 X_2\ne 0.
\]
This contradiction shows that \(d>0\). Then, similarly
to~\eqref{eq:wendliff-comm}, we have
\begin{gather*}
0=[X_k,\omega_{ii';jj'}^\perp(B)]=
\omega_{ii';jj'}^\perp((\partial_{X_k}+\partial_{Y_k})B),\\
i,i',j,j',k=1,\ldots,n,\quad i<i',\,j<j'.
\end{gather*}
Also, it is easy to see that \((\partial_{X_k}+\partial_{Y_k})B\in
S^2\kere\Delta\) and \((\partial_{X_k}+\partial_{Y_k})B_\d\in\ker\Pi^\perp\).
The minimality condition imposed on \(d\) implies that
\[
(\partial_{X_k}+\partial_{Y_k})B_\d=
\d^{\otimes 2}((\partial_{X_k}+\partial_{Y_k})B)=0,\qquad k=1,\ldots,n,
\]
or, equivalently,
\[ 
(\partial_{X_k}+\partial_{Y_k})B_{ij}=0,\qquad i,j,k=1,\ldots,n,
\]
where \(B_{ij}:=\partial_{X_i}\partial_{Y_j}B\). Since \(d>0\), by
Claim~\ref{cl:wendliff-Bmin} we have \(B_{ij}(X,Y)=H_{ij}(X-Y)\) for some
\(H_{ij}\in\ker\Delta\). In particular, \(\Pi(B_{ij})=0\), which implies that
\(\Pi_\d(B_\d)=0\) and \(\tr(\Pi_\d(B_\d))=0\). Observe also that
\begin{gather*}
(\partial_k H_{ij})(X-Y)=\partial_{X_k}(H_{ij}(X-Y))=
\partial_{X_k}(B_{ij}(X,Y))=\partial_{X_k X_i Y_j}(B(X,Y)),\\
i,j,k=1,\ldots,n,
\end{gather*}
where \(\partial_k\) denotes the partial derivative in \(k\)-th variable. It
follows that \(\partial_k H_{ij}=\partial_i H_{kj}\). A similar argument shows
that \(\partial_k H_{ij}=\partial_j H_{ik}\). Hence,
\[
\partial_k H_{ij}=\partial_i H_{kj}=\partial_j H_{ki}=\partial_k H_{ji},
\]
i.e. \(H_{ij}-H_{ji}=\mathrm{const}\). Since all \(H_{ij}\) are homogeneous, we
have \(H_{ij}=H_{ji}\) for all \(i,j=1,\ldots,n\).

By Claim~\ref{cl:wendliff-omega}, we have
\begin{gather}
\begin{multlined}[t]
\omega_{ii';jj'}^\perp(B)=-(\partial_{H_{ij'}}(p^{-1}X_{i'}X_j)+
\partial_{H_{i'j'}}(p^{-1} X_i X_j)+
\partial_{H_{ij}}(p^{-1}X_{i'}X_{j'})+\\
\partial_{H_{i'j}}(p^{-1}X_i X_{j'})),
\end{multlined}\label{eq:tdwendl-iijj2}\\
i,i',j,j'=1,\ldots,n,\quad i<i',\,j<j'.\notag
\end{gather}
We have
\begin{multline*}
\partial_{H_{ij}}(p^{-1}X_{i'}X_{j'})=
X_{i'}X_{j'}\partial_{H_{ij}}(p^{-1})+
X_{i'}\partial_{\partial_{X_{j'}}H_{ij}}(p^{-1})+
X_{j'}\partial_{\partial_{X_{i'}}H_{ij}}(p^{-1})\,+\\
\partial_{\partial_{X_{i'}X_{j'}}H_{ij}}(p^{-1}).
\end{multline*}
Hence~\eqref{eq:wendliff-dibruno} implies that
\begin{equation}\label{eq:tdwendl-modp}
p^{d+1}\partial_{H_{ij}}(p^{-1}X_{i'}X_{j'})\equiv
(-2)^d d!\,X_{i'}X_{j'}H_{ij}\pmod p,
\end{equation}
and similarly for the other terms. Hence it follows
from~\eqref{eq:tdwendl-iijj}, \eqref{eq:tdwendl-iijj2}, \eqref{eq:tdwendl-modp}
that
\begin{gather*}
X_{i'}X_j H_{ij'}+X_i X_j H_{i'j'}+
X_{i'}X_{j'}H_{ij}+X_i X_{j'}H_{i'j}\equiv 0\pmod p,\\
i,i',j,j'=1,\ldots,n,\quad i<i',\,j<j'.
\end{gather*}
Considering all options with \(\{i,i',j,j'\}\subset\{1,2,3\}\), we obtain
\[
\begin{pmatrix}
X_2^2 & X_1^2 & 0 & 2X_1X_2 & 0 & 0 \\
X_3^2 & 0 & X_1^2 & 0 & 0 & 2X_1X_3 \\
0 & X_3^2 & X_2^2 & 0 & 2X_2X_3 & 0 \\
0 & X_1X_3 & 0 & X_2X_3 & X_1X_2 & X_2^2 \\
X_2X_3 & 0 & 0 & X_1X_3 & X_1^2 & X_1X_2 \\
0 & 0 & X_1X_2 & X_3^2 & X_1X_3 & X_2X_3
\end{pmatrix}\begin{pmatrix}
H_{11} \\ H_{22} \\ H_{33} \\ H_{12} \\ H_{23} \\ H_{13}
\end{pmatrix}\equiv 0\pmod p.
\]
The determinant of the matrix on the left-hand side equals \(16X_1^4 X_2^4
X_3^4\). Since it is not divisible by \(p\), we obtain (e.g. from Cramer's
rule) that \(p\mid H_{ij}\) for all \(i,j\in\{1,2,3\}\) and since \(H_{ij}\) are
harmonic, they all vanish. Similarly, \(H_{ij}=0\) for all
\(i,j\in\{1,\ldots,n\}\) and \(B_\d=0\). This contradiction concludes the proof.
\end{proof}

\appendix
\section{Proof for the Coulomb gas on \texorpdfstring{\(\mathbb R^3\)}{R³}}
\label{sec:coulomb}

In this section we prove Theorem~\ref{th:intro-co} extending Theorem
\ref{th:intro-sm} to the case of Schr\"o\-din\-ger operators for the Coulomb
gas on \(\mathbb R^3\) in a confining potential.

\begin{proof}[Proof of Theorem~\ref{th:intro-co}]

First, we will prove a slightly more general statement, allowing arbitrary
smooth \(V\colon M\to\mathbb R\) invariant under the \(G=S_N\)-action such that
\(V(x)\to\infty\) as \(|x|\to\infty\), instead of \(V(x)=\sum_{i=1}^N v(x_i)\)
as in the formulation of the theorem. 

Consider the operator \(L_f=\Delta+V+C+f\), as in the statement of the theorem. 

It is clearly \(G\)-equivariant. This operator satisfies the strong unique
continuation property at every point away from the multi-diagonal
\[
\Sigma=\bigcup_{1\le i<j\le N}\{x\in M=(\mathbb R^3)^N\mid x_i=x_j\}
\]
as there the potential is smooth. Note that \(M\setminus\Sigma\) is connected.
Moreover, one can find a regular point for the \(G\) action away from \(\Sigma\)
in the interior of \(B_S.\) We will prove that for all \(f,\,L_f\) defines a
self-adjoint Fredholm operator (of index zero) in suitable Hilbert spaces. With
this in mind, the arguments in the proof of Theorem \ref{th:intro-sm}, including
Floer \(\varepsilon\)-spaces, go through without essential changes to establish
Theorem \ref{th:intro-co}. For instance, to prove that \(\Pi(B)\) vanishes near
a selected regular point \(x\) in the interior of \(B_S\) and away from
\(\Sigma\), we can look at \(\chi\Pi(B)\) for a non-negative \(G\)-invariant
cutoff function strictly positive near \(x\) and zero outside a \(G\)-invariant
neighborhood of the orbit \(G\cdot x\) of \(x\).

We note that writing \(C=\sum\chi(|x_i-x_j|)C_{ij}+\sum
(1-\chi(|x_i-x_j|))C_{ij}\) for \(C_{ij}=C(|x_i-x_j|)\) and a cutoff function
\(\chi\colon [0,\infty)\to [0,1]\) with \(\chi=1\) on \([0,\delta)\) and \(\chi=
0\) on \([2\delta,\infty)\) for \(\delta\) sufficiently small, and absorbing the
second sum into \(V\), and replacing \(C\) by \(\widetilde C=\sum
\chi(|x_i-x_j|)C_{ij}\) we may assume that \(C=0\) outside the \(2\delta\)
neighborhood of \(\Sigma\). Furthermore, each
\[
\widetilde C_{ij}=\chi(|x_i-x_j|)C_{ij}
\]
satisfies
\[
\widetilde{C}_{ij}\le E\frac 1{|x_i-x_j|}
\]
for a constant \(E>0\) and hence \(\widetilde C\le E\sum\frac 1{|x_i-x_j|}\). 

\subsection{Estimate on compacts}\label{sec:hardy}

Let us first show that \(L_f\) defines an operator \(H^2_\loc(M)\to
L^2_\loc(M)\). By standard arguments, it is sufficient to check this for the
multiplication operator by \(\widetilde{C}\), and hence in turn by
\(\widetilde{C}_{ij}\) for \(i<j\). We may suppose that \(i=1,\,j=2\).

Recall the following special case of Hardy's inequality~\cite{RS2019}: for a
compactly supported smooth function \(u\) on \(\mathbb R^3\),
\begin{equation}\label{eq:hardy-subj}
\left\|\frac ur\right\|_{L^2(\mathbb R^3)}\le 2\|\nabla u\|_{L^2(\mathbb R^3)},
\end{equation}
where on the left \(r\) denotes the function \(r(x)=|x|\) for the Euclidean
norm on \(\mathbb R^3\).

Let \(y=x_1-x_2,\,z=(x_1+x_2,x_3,\ldots,x_N)\). Let \(u\in C_c^\infty(M)\). It
suffices to see that \(\|\widetilde{C}_{12}u\|_{L^2(M)}\) is finite. In the new
coordinates \((y,z)\) this amounts to estimating as follows
using~\eqref{eq:hardy-subj}, where the constant \(A\) takes into account the
distortion of the metric by our change of coordinates:
\begin{multline*}
\|\widetilde C_{12}u\|^2_{L^2(M)}\le\frac E8
\int_{\mathbb R^3\times (\mathbb R^3)^{N-1}}\frac{|u(y,z)|^2}{|y|^2}\,dy\,dz\le
\frac E2\int_{\mathbb R^3\times (\mathbb R^3)^{N-1}}|\nabla_y u|^2 dy\,dz\le\\
\frac{EA}2\|\nabla u\|^2_{L^2(M)}<\infty.
\end{multline*}

Moreover, the proof above yields the following estimate. 

\begin{corollary}\label{cor:hardy-estimate}
There exists a constant \(K\) depending only on the dimension of \(M\), such
that for all \(u\in H^1(M)\),
\[
\|Cu\|_{L^2(M)}\le K\|\nabla u\|_{L^2(M)}.
\]
\end{corollary}

\subsection{Self-adjoint extension with smooth potential}\label{sec:sa}

Next we show the following for the potential \(W=V+f\).

\begin{proposition}\label{pr:sa-main}
The operator \(H_0=\Delta+W\) on \(C^\infty_c(M)\) defines an
essentially self-adjoint operator on \(L^2(M)\) whose closure, which we denote
by \(H\), is its unique self-adjoint extension and has domain
\[
D=\{u\in L^2(M)\cap H^2_{\loc}(M)\mid\Delta u+Wu\in L^2(M)\}.
\]
\end{proposition}

Note that in view of elliptic regularity, \(D\) coincides with the space of \(u
\in L^2\) such that \(\Delta u+Wu\in L^2(M)\) as distributions.

Set \(H_0=\Delta+W\). Since \(V\) is continuous and tends to \(+\infty\) as
\(|x|\to\infty\), and \(f\) has compact support, \(W\) has a finite lower
bound \(m\). Choose \(c\in\mathbb R\) so that
\[
Z:=W+c\ge 1,\qquad A_0:=H_0+cI\ge I.
\]
Adding a real constant preserves essential self-adjointness and the operator
domain. It is therefore enough to prove the theorem for
\[
A_0=\Delta+Z,
\]
whose closure we denote by \(A\). We proceed in two steps. First we prove an
abstract criterion in terms of the kernel of the adjoint, and then provide a
direct cutoff calculation showing that \(\ker A_0^*=\{0\}\).

\begin{lemma}\label{lem:sa-criterion}
Let \(S_0\) be a densely defined symmetric operator on a Hilbert space
\(\mathcal H\), and suppose
\[
\langle S_0u,u\rangle\ge\|u\|^2,\qquad u\in\Dom(S_0).
\]
Then \(S_0\) is essentially self-adjoint if and only if \(\ker S_0^*=\{0\}\).
\end{lemma}

\begin{proof}
A densely defined symmetric operator is closable. Set \(S=\overline{S_0}\).
Then \(S\) is closed and symmetric, \(S\ge I\), and \(S^*=S_0^*\).
Assume \(\ker S_0^*=\{0\}\).

\emph{First, \(\im S\) is closed.} Cauchy--Schwarz gives
\begin{equation}\label{eq:sa-coercive}
\|Su\|\,\|u\|\ge\langle Su,u\rangle\ge\|u\|^2,
\qquad\text{hence}\qquad\|Su\|\ge\|u\|.
\end{equation}
If \(Su_j\to f\) in \(\mathcal H\), then
\[
\|u_j-u_k\|\le\|Su_j-Su_k\|,
\]
so \((u_j)\) is Cauchy and converges to some \(u\in\mathcal H\). Closedness of
\(S\) now implies \(u\in\Dom(S)\) and \(Su=f\).

\emph{Second, \(\im S\) is dense.} For any densely defined operator,
\[
(\im S)^\perp=\ker S^*.
\]
This follows directly from the definition of the adjoint: \(v\) is orthogonal to
every \(Su\) precisely when \(v\in\Dom(S^*)\) and \(S^*v=0\). Our hypothesis
therefore makes \(\im S\) dense. Since it is also closed,
\begin{equation}\label{eq:sa-onto}
\im S=\mathcal H.
\end{equation}

\emph{Third, \(S\) is self-adjoint.} Take \(v\in\Dom(S^*)\). By
\eqref{eq:sa-onto}, choose \(u\in\Dom(S)\) with \(Su=S^*v\). Symmetry gives
\(S\subset S^*\), hence
\[
S^*(v-u)=S^*v-Su=0.
\]
The triviality of \(\ker S^*\) implies \(v=u\in\Dom(S)\). Thus
\(\Dom(S^*)\subseteq\Dom(S)\); symmetry gives the reverse inclusion, so
\(S=S^*\).

Conversely, if \(S_0\) is essentially self-adjoint, then \(S_0^*=S\ge I\).
Estimate~\eqref{eq:sa-coercive} implies \(\ker S_0^*=\{0\}\).
\end{proof}

\begin{remark}
More generally, if \(S_0\ge aI\) and \(\lambda<a\) is real, then
\[
S_0\text{ is essentially self-adjoint}\quad\Longleftrightarrow\quad
\ker(S_0^*-\lambda I)=\{0\}
\]
by an application of Lemma~\ref{lem:sa-criterion} to \((S_0-\lambda
I)/(a-\lambda)\).
\end{remark}

Now we apply Lemma~\ref{lem:sa-criterion} to prove Proposition~\ref{pr:sa-main}.
Suppose that \(u\in\ker A_0^*\). By the distributional description of the
adjoint,
\begin{equation}\label{eq:sa-weak}
(\Delta+Z)u=0\quad\text{on }M,\qquad u\in L^2(M).
\end{equation}

Since \(Z\) is smooth, \(Zu\in L^2_\loc\). Interior elliptic regularity
therefore gives \(u\in H^2_\loc( M)\).

Choose real cutoffs \(\chi_R\in C_c^\infty( M)\) satisfying
\[
0\le\chi_R\le 1,\quad\chi_R=1\text{ on }B_R,\quad
\operatorname{supp}\chi_R\subset B_{2R},\quad |\nabla\chi_R|\le C/R.
\]
Multiply~\eqref{eq:sa-weak} by \(\chi_R^2\overline u\), integrate, and take real
parts. Integration by parts is justified by the local regularity and compact
support of the test function. It yields
\[
0=\int\chi_R^2|\nabla u|^2
+2\operatorname{Re}\int\chi_R\overline u\,\nabla u\cdot\nabla\chi_R
+\int Z\chi_R^2|u|^2.
\]
Expanding \(|\nabla(\chi_Ru)|^2\) gives the exact identity
\[
\int_{M}\Bigl(|\nabla(\chi_Ru)|^2+Z|\chi_Ru|^2\Bigr)
=\int_{M}|\nabla\chi_R|^2|u|^2
\le\frac{C^2}{R^2}\|u\|_2^2.
\]
Because \(Z\ge1\), we obtain
\[
\int_{B_R}|u|^2\le\|\chi_Ru\|_2^2\le\frac{C^2}{R^2}\|u\|_2^2.
\]
Letting \(R\to\infty\) proves \(u=0\). Lemma~\ref{lem:sa-criterion} proves
essential self-adjointness of \(A_0\), and hence Proposition~\ref{pr:sa-main}
for \(H_0\).

\begin{remark}
In particular, \textbf{a smooth real potential bounded below already suffices};
the growth to \(+\infty\) is stronger than needed for this conclusion.
\end{remark}

\subsection{Kato's boundedness of the singular multiplication operator}

We now look at \(L_f=\widetilde H=H+C\) as a perturbation of \(H\) and prove
that it is self-adjoint with the same domain as \(H\) by means of Kato's
perturbation theorem.

\begin{theorem}
The operator \(\widetilde H=H+C\) is self-adjoint with domain exactly
\(\Dom(H)\), and the graph norms of \(H\) and \(\widetilde H\) are equivalent.
The operator \(\Delta+V+f+C\) on \(C^\infty_c(M)\) is essentially self-adjoint,
with closure \(\widetilde H\).
\end{theorem}

As before, it is enough to consider the operator \(\widetilde A=A+C\). Note that
for \(u\in\Dom(A)\),
\begin{equation}\label{eq:kato-energyA}
\|\nabla u\|_2^2\le\|Au\|_2\|u\|_2.
\end{equation}
Indeed, by calculating on \(u \in C^\infty_c(M)\) and then approximating in the
graph norm, we get that for \(u\in\Dom(A)\),
\begin{equation}\label{eq:kato-formineq}
\|\nabla u\|_2^2+\|Z^{1/2}u\|_2^2\le\langle Au,u\rangle\le\|Au\|_2\|u\|_2.
\end{equation}
In particular, Corollary~\ref{cor:hardy-estimate} and~\eqref{eq:kato-energyA}
imply that for every \(u\in\Dom(A)\) and every \(\varepsilon>0\),
\[
\|Cu\|_2\le K\|Au\|_2^{1/2}\|u\|_2^{1/2}
\le\varepsilon\|Au\|_2+\frac{K^2}{4\varepsilon}\|u\|_2.
\]
In particular, \(\Dom(A)\subset\Dom(C)\), and moreover \(C\) is relatively
\(A\)-bounded with coefficient \(\epsilon<1\). The Kato--Rellich theorem says
that a symmetric perturbation of relative bound less than one preserves
self-adjointness on the original domain and preserves operator cores: see for
example~\cite[Theorem 1.1]{Kat1966} for closability and the domains coinciding,
and~\cite[Theorem 6.4, p.~135]{Teschl} and~\cite[Theorem 7.3]{Simon} for some
more precise references. Applying it to \(A\) and \(C\) proves the theorem and
gives
\[
\widetilde A=A+C\text{ self-adjoint on }\Dom(A),\qquad
\widetilde H:=\widetilde A-cI=H+C\text{ on }\Dom(H).
\]


\subsection{The Fredholm property}

It follows immediately from Section~\ref{sec:sa} that \(L_f=H+C\colon D\to
L^2(M)\) is a well-defined operator which is bounded on \(D\) with the graph
norm
\[
\|u\|_D:=\|u\|_{L^2(M)}+\|Au\|_{L^2(M)}.
\]
We claim that this bounded operator is Fredholm. 

First, by an argument similar to the proof of Rellich's embedding theorem, as
\(V(x)\to\infty\) for \(|x|\to\infty\), the inclusion \(i\colon D \to L^2(M)\)
is a compact operator. Indeed, for \(R>0\) let \(B_R\) be the ball of radius
\(R\) in \(M\). Suppose \(\|u\|_D\le Q\) for some \(Q>0\).
Then~\eqref{eq:kato-formineq} implies that for all \(R>0\)
\[
\|u\|_{H^1(B_R)}\le Q,
\]
while
\[
\|u\|_{L^2(M\setminus B_R)}\le\frac Q{\inf_{M\setminus B_R}\sqrt Z}.
\]
Now for a sequence \(u_i\in D\) with \(\|u_i\|_D\le Q\) choose a sequence
\(R_k\to\infty\), \(k\ge 1\). By Rellich's embedding, there is a subsequence,
say \(v_j\), with Cauchy restriction to \(L^2(B_{R_k})\) for all \(k\). As
\(Z(x)\to\infty\) as \(|x|\to\infty\),
\[
\|v_j\|_{L^2(M\setminus B_{R_k})}\to 0\quad\text{as}\quad k\to\infty
\]
uniformly in \(j\), we obtain that \(v_j\) is Cauchy in \(L^2(M)\), finishing
the proof. 

Furthermore, \(T=H+C+(c+d)I\), where we pick \(c\) as before and \(d\) so that
\(C+d\ge 0\), is a self-adjoint operator on \(L^2(M)\) with domain \(D\). It is
moreover a bounded operator \(D\to L^2(M)\). Furthermore 
\[
\langle Tu,u\rangle=\langle Au,u\rangle+\langle (C+d)u,u\rangle\ge\|u\|_2^2.
\]
Hence \(T\ge I\) and therefore \(T\colon D\to L^2(M)\) is bijective and
therefore an isomorphism by the bounded inverse theorem. Finally,
\(L_f=\Delta+V+C+f=H+C=T-(c+d)I\) is Fredholm as a compact perturbation of an
isomorphism.
\end{proof}

\bibliographystyle{alpha}
\bibliography{mybib.bib}

\end{document}